\documentclass[11pt,letterpaper]{article}
\usepackage{graphicx,epstopdf,subfigure,mathtools,mathrsfs, amsmath, amssymb} 
\usepackage[font=small,labelfont=bf]{caption}
\usepackage{float, cancel, amssymb}
\usepackage[normalem]{ulem}
\usepackage[title]{appendix}
\PassOptionsToPackage{usenames,dvipsnames}{xcolor}
\usepackage[usenames,dvipsnames]{xcolor}
 
\usepackage[margin=1 in]{geometry}
\usepackage{cite}
\usepackage{caption}
\usepackage[font=small,labelfont=bf]{caption}

\usepackage{amsfonts, soul, authblk, enumitem}

\newcommand{\rev}[1]{\color{black}#1\normalcolor}

\newcommand{\widebar}[1]{%
   \hbox{%
     \vbox{%
       \hrule height 0.5pt 
       \kern0.5ex
       \hbox{%
         \kern-0.1em
         \ensuremath{#1}%
         \kern-0.1em
       }%
     }%
   }%
} 

\graphicspath{{../}}
\newfloat{algorithm}{!t}{loa}
\providecommand{\algorithmname}{Algorithm}
\floatname{algorithm}{\protect\algorithmname}

\newcommand{\V}[1]{\boldsymbol{#1}}                 
\newcommand{\M}[1]{\boldsymbol{#1}}
\newcommand{\Lop}[1]{\boldsymbol{\mathcal{#1}}}

\newcommand{\ind}[2]{{#1}^{(#2)}}

\newcommand{\tdisc}[2]{#1_{#2}}

\global\long\def\dt#1{\partial_t#1} 
\global\long\def\ds#1{\partial_s#1} 
\global\long\def\Xmp{\V{X}_\text{MP}}
\global\long\def\Ump{\V{U}_\text{MP}}
\global\long\def\X{\M{\Lop{X}}}

\global\long\def\Xuni{\widebar{\V{X}}}
\global\long\def\Xu{\V{X}^{(u)}}
\global\long\def\Bind#1{{\{#1\}}}
\global\long\def\Pot{\mathcal{E}}
\global\long\def\Nhat{\widehat{\M{N}}}
\global\long\def\eqdi{\stackrel{d}{=}}

\global\long\def\D#1{\Delta#1}

\global\long\def\norm#1{\left\Vert #1\right\Vert }

\global\long\def\kon{k_\text{on}}
\global\long\def\koff{k_\text{off}}

\global\long\def\rmax{r_\text{max}}
\global\long\def\Nseg{N_\text{seg}}
\global\long\def\Lseg{\ell_s}
\global\long\def\Slet#1{\M{S}\left(#1\right)}
\global\long\def\Dlet#1{\M{D}\left(#1\right)}

\global\long\def\Xs{\V{\tau}}
\global\long\def\EPMI{\frac{1}{8\pi\mu}}
\DeclareMathOperator*{\argmin}{arg\,min}
\global\long\def\ML{ \widetilde{\M{M}}^\text{L}}

\global\long\def\MNL{\widetilde{\M{M}}^\text{NL}}

\global\long\def\koff{k_\text{off}}
\global\long\def\kon{k_\text{on}}
\global\long\def\konb{k_\text{on,s}}
\global\long\def\koffb{k_\text{off,s}}

\global\long\def\eps{{\hat{a}}}
\global\long\def\epsh{{\hat{a}}^*}
\global\long\def\epsRS{\hat{\epsilon}}
\global\long\def\epsRSh{\hat{\epsilon}^*}
\global\long\def\rc{a}
\global\long\def\epsc{\epsilon}
\global\long\def\dt#1{\partial_t#1} 

\global\long\def\Xsbar{\bar{\Xs}}

\global\long\def\Wproc{\V{\mathcal{W}}}
\global\long\def\gauss{\V{\eta}}
\global\long\def\Mfor{\widetilde{\M{M}}}
\global\long\def\div{\partial}
\global\long\def\Wt{\widetilde{\M{W}}}
\global\long\def\Mref{\widetilde{\M{M}}_\text{os}}
\global\long\def\Msqs{\widetilde{\M{M}}_\text{SQS}}
\global\long\def\fext{\V{f}^{\text{(ext)}}}

\global\long\def\Fext{\V{F}^{\text{(ext)}}}
\global\long\def\FCL{\V{F}^{\text{(CL)}}}
\global\long\def\Rhat{\widehat{\V R}}
\global\long\def\XPoly{\mathbb{X}}
\global\long\def\TauPoly{\mathbb{T}}

\usepackage{hyperref}
\hypersetup{
    colorlinks=false,
    pdfborder={0 0 0},
}

\title{A simulation platform for slender, semiflexible, and inextensible fibers with Brownian hydrodynamics and steric repulsion}

\date{}

\author[1,2,$*$]{Ondrej Maxian}
\author[3,$\dagger$]{Aleksandar Donev}
\affil[1]{\small Department of Molecular Genetics and Cell Biology, University of Chicago, Chicago, IL 60637}
\affil[2]{\small Institute for Biophysical Dynamics, University of Chicago, Chicago, IL 60637}
\affil[3]{\small Courant Institute of Mathematical Sciences, New York University, New York, NY 10012}
\affil[$\dagger$]{d.\ November 2, 2023}
\affil[*]{\small Corresponding author: ondrej@uchicago.edu}

\begin{document}
\maketitle

\vspace{-1.2 cm}

\begin{abstract}
The last few years have witnessed an explosion of new numerical methods for filament hydrodynamics. Aside from their ubiquity in biology, physics, and engineering, filaments present unique challenges from an applied-mathematical point of view. Their slenderness, inextensibility, semiflexibility, and meso-scale nature all require numerical methods that can handle multiple lengthscales in the presence of constraints. Accounting for Brownian motion while keeping the dynamics in detailed balance and on the constraint is difficult, as is including a background solvent, which couples the dynamics of multiple filaments together in a suspension. In this paper, we present a simulation platform for deterministic and Brownian inextensible filament dynamics which includes nonlocal fluid dynamics and steric repulsion. \rev{For nonlocal hydrodynamics}, we define the mobility on a single filament using line integrals of Rotne-Prager-Yamakawa regularized singularities, and numerically preserve the symmetric positive definite property by using a thicker regularization width for the nonlocal integrals than for the self term. For steric repulsion, we introduce a soft local repulsive potential defined as a double-integral over two filaments, then present a scheme to identify and evaluate the nonzero components of the integrand. Using a temporal integrator developed in previous work, we demonstrate that Langevin dynamics sample from the equilibrium distribution of free filament shapes, and that the modeling error in using the thicker regularization is small. We conclude with two examples, sedimenting filaments and cross-linked fiber networks, in which nonlocal hydrodynamics does and does not generate long-range flow fields, respectively. In the latter case, we show that the effect of hydrodynamics can be accounted for through steric repulsion.
\end{abstract}

\section{Introduction}

Fibers, fibers, fibers! Perhaps this is a recency illusion, but it feels as if the applied mathematics literature of the last few years has teemed with new efficient methods for computing the dynamics of fibers in flow \cite{ehssan17,walker2020efficient, keavRPY,inexRS}. The driving force behind this flurry of activity is likely a result of modern computing power \cite{krishna2022petascale,yan2022toward,koshakji2023robust} meeting decades-old problems that involve filaments, such as the motion of flagellated bacteria \cite{blum1979biophysics,lauga2009hydrodynamics}, rheology of filament suspensions \cite{fibsusps,rahnama1995effect,mackfibs}, and dynamics of the cytoskeleton \cite{theriot1991actin,nazockdast2017cytoplasmic}. Another reason might be the unique challenges filaments pose to simulate: they are \emph{slender}, which endows the problem with multiple lengthscales, \emph{inextensible}, which adds a constraint to the dynamics, and \emph{semiflexible}, which often makes Brownian motion important in setting their shapes. All of the methods presented over the last few years have incorporated novel ways of addressing these challenges. Yet, despite all of this progress, a method that can deal with all of them at once remains elusive. 

Even prior to the recent flurry of activity, the literature in the area of filament dynamics was quite vast, and the reader will pardon any over-simplifications on our part. Painting with a broad brush, there seem to be two main sub-areas: first, there is the question of how to resolve the \emph{deterministic} interactions of a slender filament with the fluid that surrounds it. This area is primarily occupied by applied mathematicians, and has seen the application of immersed boundary methods \cite{peskin1972flow} and slender body theories \cite{batch70,krub} to the filament problem. Then, on the opposite end, there is the question of how to formulate and integrate the equations of \emph{Brownian} dynamics for filaments, regardless of the medium they are immersed in. This seems to be the purview of engineers and physicists \cite{morse2003theory,butler2005brownian, ottinger2012stochastic, freedman2017versatile}, and methods in this area have typically used approximations to the hydrodynamic interactions based on local drag.

Beginning with the first part of the literature, there is by now a long history of using immersed boundary methods to simulate the dynamics of filaments. In the classical immersed boundary method as formulated by Peskin \cite{peskin1972flow,peskin2002acta}, the fiber has a Lagrangian representation, on which the internal forces can be computed and ``spread'' to a background fluid grid. Solving the fluid equations on the grid and ``interpolating'' the result back to the filament yields a velocity of the structure, allowing for an update in its position. The cornerstone of the immersed boundary method is the regularized delta function which is used to spread the force onto the grid, thereby imparting a thickness \rev{(called the hydrodynamic radius)} on which the fluid ``feels'' the fiber \cite{ttbring08}. In the classical IB method, the width of the regularized delta (or ``blob'') function is tied to the width of the Eulerian fluid grid. Other extensions based on the same idea eliminate the fluid grid entirely by using a clever choice of delta function. For example, the force coupling method \cite{fcm03,fcm10,keavRPY} uses a Gaussian blob function, the Rotne-Prager-Yamakawa kernel uses a spherical surface delta function \cite{rpyOG, wajnryb2013generalization,keavRPY}, and the method of regularized Stokeslets uses other carefully-chosen functions depending on the dimension \cite{cortez2001method,cortez2005method} (notably omitting the inerpolation step). Even though these methods lack a fluid grid, they retain the physics of the blob function representing the thickness of the fiber. As such, modeling slender fibers requires a blob function with width $\epsc$ (the fiber aspect ratio). Even though fast methods exist to compute the velocity due to many blobs in linear time \cite{griffith2007adaptive,rostami2016kernel,PSRPY, su2024accelerating}, a naive representation of the filament would require $1/\epsc$ blobs, making even the most optimized methods impractical.

Because of this, filament representations based on line-integrals of regularized singularities have become more popular \cite{cortez2012slender,ohm2021remarks,zhao2021regularized}. In the case when the line integrals are discretized with direct quadrature, this of course has no advantage, but in the cases when special quadrature (or analytical) techniques are available to integrate the kernel, the spatial discretization can be decoupled from the width of the regularization function. For example, if the integral can be computed analytically for straight lines, a fiber can be discretized into a series of rigid segments \cite{cortez2018regularized,hall2019efficient,huang2021numerical,koshakji2023robust}, the number of which are a function of the expected fiber shape rather than its aspect ratio. Consequently, because only the local interactions are a strong function of aspect ratio, interactions between distinct fiber segments (potentially on distinct fibers) can be integrated directly with $\mathcal{O}(1)$ points, and the method has cost independent of aspect ratio. 

A related approach to filament hydrodynamics, which is also based on integral equations, is to use 
slender body theory (SBT) \cite{batch70, krub,johnson,gotz2001interactions,SBTRev} to asymptotically remove the smallest lengthscale in the problem. The original result, which was based on singularities and matched asymptotic expansions \cite{krub}, can be recast as an asymptotic reduction of the three-dimensional boundary integral equations on the filament surface, whereby the forces and velocities are expanded in a Fourier representation (over the cross section), and only the constant mode is kept \cite{koens2018boundary}. The result of this analysis is a local drag mobility matrix which is $\mathcal{O}(\log(\epsc))$ and expresses the velocity due to force concentrated $\mathcal{O}(\epsc)$ from the point of interest, plus an $\mathcal{O}(1)$ integral that expresses the velocity from the rest of the filament. Because the $\mathcal{O}(\epsc)$ parts are concentrated in the local drag term, the smallest dimension of the problem is integrated out in a sense, and it becomes attractive to use SBT in numerical simulations \cite{ts04,young2009hydrodynamic,ehssan17}. Yet methods based on SBT are plagued by a number of issues, foremost among which is the nonlocal finite part integral, which causes the mobility to have negative eigenvalues when the filament is over-resolved \cite{mori2020accuracy}. Approaches to deal with this include regularizing the integrand \cite{ts04, andersson2021integral} and putting discretization nodes on the cross sections to avoid singularities \cite{koshakji2023robust}. The tricky nature of regularization and near-singular integration often lead the practitioner away from SBT and towards the more implementation-friendly regularized schemes. This is especially the case for Brownian dynamics, where the square root of the (symmetric positive definite) mobility is required.

Once the method for the hydrodynamics is specified, it remains to deal with the inextensibility constraint. Most immersed boundary methods use penalty terms for the purposes of treating constraints \cite{lim2008dynamics, liu2019efficient,lee2021novel}. To eliminate the potentially stiff timescales associated with these forces, recent studies have explicitly solved for the forces required to keep the dynamics on the constraint \cite{keavRPY, inexRS}. Yet, unless the solver is ultimately nonlinear \cite{keavRPY}, the dynamics will always drift off the constraint numerically, and including penalty parameters \cite{ts04,ehssan17} to correct for this seems to defeat the purpose of implementing a constrained method in the first place. Because the method of \cite{keavRPY} uses discrete regularized singularities with direct integration, it seems there is no method in the literature which harmonizes exact treatment of the inextensibility constraint with efficient quadrature for the mobility in three dimensions (see \cite{hall2019efficient} for two dimensions).

Returning to our second broad area of the literature, there are a number of methods for Brownian filament fluctuations, but these do not treat the hydrodynamics between filament pieces or distinct filaments. Typically, these methods take the form of a set of beads connected by springs, which are governed by stretching and bending energies and are therefore constraint free \cite{freedman2017versatile, lin2014combined,cunha2022settling}. Similar versions of these models fix the distance between beads, which confines the dynamics to a certain constrained manifold \cite{butler2005brownian}. This has an effect on the (Ito) Langevin equation that governs equilibrium Brownian dynamics, as constraints give rise to new stochastic drift terms which can be quite complex \cite{morse2003theory, ottinger2012stochastic}. These terms can in principle be handled using Fixman's method \cite{fixman1978simulation,butler2005brownian,delong2015brownian}, but this requires a (potentially costly) additional resistance solve every time step.

Separate from the chain representation is the matter of Brownian hydrodynamics. Here the slender body theory problem rears its ugly head once again. If the mobility is not guaranteed to be symmetric positive definite (as is the case also in regularized Stokeslets because of a lack of symmetry in the spread/interpolate operations), it is impossible to define its square root, which is necessary for fluctuation-dissipation balance. Thus, the available methods that attempt Brownian motion with a faithful representation of hydrodynamics either approximate the hydrodynamic interactions (not using the full SBT) \cite{butler2005brownian,manikantan2016effect}, or use regularized singularities (IB methods) \cite{liu2019efficient}. Indeed, while the spread/interpolate symmetry in regularized singularity methods gives an automatically SPD mobility, \rev{for slender filaments the number of regularized singularities required to accurately resolve the hydrodynamics is still prohibitively large.} 

All of this literature gently points us to a path that harmonizes regularized singularities (with their SPD properties) with slender body theory (to remove the smallest lengthscale), combined with constrained dynamics (to remove stiff penalty terms) and efficient handling of stochastic drift terms (to integrate the equations of Brownian dynamics) \cite{delong2015brownian}. This is more or less the path we have followed over the past few years. First, we showed that slender body theory is asymptotically equivalent to a line integral of Rotne-Prager-Yamakawa (RPY) singularities \cite{rpyOG,maxian2021integral}. Second, we developed an efficient quadrature scheme for the self RPY integral, fully decoupling the degrees of freedom from the fiber aspect ratio \cite{maxian2022hydrodynamics}. And third, we showed how to implement the special quadrature scheme in a constrained Brownian framework for a single filament \cite{maxian2023bending}. Until now, however, we were not able to formulate a method for Brownian dynamics of multiple filaments; that is, a method with Brownian fluctuations \emph{and} inter-fiber hydrodynamic interactions that has cost independent of the fiber aspect ratio. The primary purpose of this paper is to present such a method, while collecting, in one readable place, all of the main equations and numerical methods that were developed in previous work and are preserved in this final implementation. Finally, to complete the simulation framework, in this paper we present a novel method for steric repulsion (to keep the fibers well separated). While this method is based on penalty forces rather than newer constraint-based ideas \cite{broms2023barrier, yan2019computing, yan2022toward, ferguson2021intersection}, our tests show that it is effective at keeping the fibers apart while reducing the required time step size by at most a factor of ten. 

The paper is therefore laid out as follows: in Section \ref{sec:space}, we introduce the continuum and discrete equations of motion for deterministic dynamics. Here we leave the mobility $\Mfor$ (force-velocity relationship) general, and focus more on the inextensibility constraint. Once a discrete evolution equation is obtained, it becomes straightforward to introduce Brownian motion in a manner consistent with detailed balance \cite{makino2004brownian,delong2015brownian}. The formulation of Brownian motion in Section \ref{sec:BM} demonstrates what we need from the mobility $\Mfor$, in particular its symmetric-positive-definiteness. This paves the way for Section \ref{sec:Mob}, where we introduce the continuum mobility, its discretization, and the key step of ``fattening'' the nonlocal mobility which allows for computing hydrodynamic interactions in an SPD manner with cost independent of $\epsc$. In Section \ref{sec:Sterics}, we introduce a novel algorithm for steric repulsion which keeps the fibers from passing through each other. Importantly, both the mobility and steric interaction algorithms are based on access to a continuum representation of the fiber $\XPoly(s)$, which is available in our discretization. After presenting temporal integrators for deterministic and Brownian motion in Section \ref{sec:TInt}, we show results of numerical tests and large-scale simulations in Sections \ref{sec:NT} and \ref{sec:Bundling}, where we simulate fibers in their equilibrium state, fibers sedimenting under gravity, and cross-linked actin networks. The sedimentation example in particular reveals the benefits of the new ``fat-corrected'' mobility formulation, while the cross linking fiber example shows the importance of resolving steric interactions.
 
\section{Equations of motion \label{sec:space}}
In this section, we lay out the governing equations for the discrete spatial variables $\V{X}$ which describe the filament centerlines. This material is entirely a review of previous publications \cite{maxian2021integral, maxian2022hydrodynamics, maxian2023bending}, but it appears here for completeness, and in an effort to collect all of the key information in one place. In that spirit, we begin with a continuum formulation, which only makes sense in the deterministic context. After the continuum formulation, we present the discretization in space, which uses a spectral collocation method and carefully handles the nonlinearities by using two different grids for the tangent vectors and collocation nodes. The Langevin equation for Brownian dynamics follows from the deterministic dynamics and detailed balance.

\subsection{Continuum formulation for deterministic dynamics}
Let us begin with a continuous curve $\XPoly(s)$ which describes the centerline of a filament. In this paper, the filaments are inextensible with constant length $L$, and so $s \in [0,L]$ is an arclength parameterization. The filament shapes have an associated bending energy density
\begin{equation}
\label{eq:BendEnd}
\Pot^{(\kappa)}\left[\XPoly\right] = \frac{\kappa}{2} \int_0^L  \ds^2 \XPoly(s) \cdot \ds^2 \XPoly(s) \, ds.
\end{equation}
\rev{Using the principle of virtual work, the pointwise force density on the filament can be found by taking a variational derivative of the energy \eqref{eq:BendEnd}
\begin{equation}
\V{f}^{(\kappa)}=-\frac{\delta  \mathcal{E}^{(\kappa)}}{\delta \XPoly}=-\kappa \ds^4 \XPoly.
\end{equation}
The same procedure also gives the free-fiber boundary conditions \cite[Sec.~2.1]{li2013sedimentation}
\begin{equation}
\label{eq:FrBCs}
\ds^2 \XPoly(s=0,L)=\V 0 \qquad \ds^3 \XPoly(s=0,L)=\V 0.
\end{equation}
Since these ``natural'' boundary conditions come from differentiating the energy functional in continuum, they can be implemented in the discrete setting by discretizing the energy \eqref{eq:BendEnd} directly and differentiating the result (see Section \ref{sec:BForce}).}

To obtain velocity from forces, we introduce the linear mobility operator $\Lop{M}\left[\XPoly \right]$, and set the velocity of the filament centerline $\V{U} = \Lop{M}\V{f}$. Note that because the power dissipated in the fluid is always positive,
\begin{equation}
\label{eq:Pmob}
\langle \V{U}, \V{f} \rangle_{L^2} = \langle \Lop{M}\V{f}, \V{f} \rangle_{L^2} > 0,
\end{equation}
the operator $\Lop{M}$ is symmetric positive definite with respect to the $L^2$ inner product. Denoting any external force density (e.g., gravity) by $\fext$, the equation of motion for the fiber centerline is 
\begin{equation}
\label{eq:Motion}
\dt \XPoly = \Lop{M}\left[\XPoly \right] \left(\V{f}^{(\kappa)}+\V{\lambda}+\fext\right). 
\end{equation}
Here $\V{\lambda}(s)$ is a constraint force which enforces the inextensibility constraint,
\begin{equation}
\label{eq:Constr}
\TauPoly(s) \cdot \TauPoly(s) = 1,
\end{equation} 
where $\TauPoly=\ds \XPoly$. 

\subsubsection{Closing the equations: the kinematic operator $\Lop{K}$}
To close the equation of motion, we differentiate the constraint \eqref{eq:Constr} with respect to time to yield 
\begin{equation}
\label{eq:TauCross}
\dt \TauPoly \cdot \TauPoly = 0 \rightarrow \dt \TauPoly = \V{\Omega} \times \TauPoly,
\end{equation}
where $\V{\Omega}(s)$ denotes the rotation rate of the tangent vectors on the unit sphere. Integrating \eqref{eq:TauCross}, we obtain a representation of all inextensible motions
\begin{equation}
\label{eq:Kdef}
\dt \XPoly(s)= \Ump + \int_{L/2}^s \left(\V{\Omega}\left(s^\prime \right) \times \TauPoly\left(s^\prime\right)\right) \, ds' :=\left(\Lop{K}\V{\alpha}\right)(s),
\end{equation}
in terms of the degrees of freedom $\V{\alpha}=\left(\V{\Omega}, \Ump\right)$, which are the tangent vector rotation rates $\V{\Omega}$ \rev{and velocity of the fiber midpoint $\Ump$}. Equation \eqref{eq:Kdef} defines the \emph{kinematic} operator $\Lop{K}\left[\XPoly\right]$.

While the representation of inextensible motions \eqref{eq:Kdef} reduces the number of degrees of freedom (removing rotations parallel to the tangent vector), it does not close the dynamics \eqref{eq:Motion} because it does not constrain the forces $\V{\lambda}$. This needs to be done by imposing the principle of virtual work, which states that the constraint forces $\V{\lambda}$ dissipate no power in the fluid with respect to all possible inextensible motions \cite{varibook, maxian2021integral}
\begin{equation}
\label{eq:Vwork}
\mathcal{P}_\lambda = \langle \V{\lambda}, \V{U} \rangle_{L^2} =\langle \V{\lambda}, \Lop{K}\V{\alpha} \rangle_{L^2}=\langle \Lop{K}^* \V{\lambda}, \V{\alpha} \rangle =0 \rightarrow \Lop{K}^* \V{\lambda}=\V 0.
\end{equation}
The last equality defines $\Lop{K}^*[\XPoly]$ as the $L^2$ adjoint of $\Lop{K}$. In numerical methods, this definition is sufficient, as we will form $\M{K}$ as a matrix, then apply an $L^2$ weights matrix to $\M{K}^T$ to obtain a representation of $\Lop{K}^*$. 

In previous formulations of inextensibility based on Euler beam theory \cite{ts04}, the equations of motion are closed by differentiating the constraint with respect to time to obtain $\ds \left(\dt \XPoly \right) \cdot \TauPoly =0$, then substituting $\V{\lambda} = \ds \left(T \TauPoly\right)$ to obtain a \emph{line tension equation} for the tensions $T$, which is solved simultaneously with the equation of motion \eqref{eq:Motion}. We choose to close the formulation differently for two main reasons: first, solving for the rotation rates $\V{\Omega}$ is advantageous in numerical methods, since we can rotate the tangent vectors by those rates to keep the dynamics exactly on the constraint. Second, the kinematic constraint $\Lop{K}^*\V{\lambda}=0$ allows us to eliminate $\V{\lambda}$ and write a closed-form equation for $\XPoly$, which will allow us to write the Langevin equation for thermal fluctuations. Nevertheless, it is instructive to insert the definition of $\Lop{K}$ in \eqref{eq:Kdef} into the virtual work constraint \eqref{eq:Vwork} to show that the two formulations are equivalent,
\begin{align}
\nonumber
\mathcal{P}_\lambda&=\int_0^L \V{\lambda}(s) \cdot \Ump + \int_0^L \V{\lambda}(s) \cdot \left(\int_0^s \left(\V{\Omega}\left(s^\prime \right) \times \TauPoly\left(s^\prime\right)\right) \, ds'\right)\\
\nonumber
&= \int_0^L \V{\lambda}(s) \cdot \Ump + \int_0^L \left(\V{\Omega}\left(s^\prime \right) \times \TauPoly\left(s^\prime\right)\right) \cdot \left(\int_{s'}^L \V{\lambda}(s) \, ds\right) \, ds'\\
\label{eq:Powerz}
&= \int_0^L \V{\lambda}(s) \cdot \Ump + \int_0^L \V{\Omega}\left(s^\prime \right) \cdot \left( \TauPoly\left(s^\prime\right) \times \left(\int_{s'}^L \V{\lambda}(s) \, ds\right)\right) \, ds'=0.
\end{align}
Since \eqref{eq:Powerz} must hold for all $\Ump$ and $\V{\Omega}$, we have that 
\begin{equation}
\Lop{K}^* \V{\lambda} = \begin{pmatrix} \displaystyle{\TauPoly\left(s\right) \times \left(\int_{s}^L \V{\lambda}(s') \, ds'\right)}  \\[10pt]  \displaystyle{\int_0^L \V{\lambda}(s)} \end{pmatrix}=\V 0.
\end{equation}
If the fiber is not straight, the first condition implies that $\V{\lambda}(s)= \ds\left(T(s) \TauPoly(s)\right)$ for some scalar function $T(s)$ with $T(L)=0$, while the second equation gives $T(0)=0$. In the case of a straight fiber, we have $\V{\lambda}(s)=\left(\ds T(s)\right) \TauPoly$, and we can only recover $T(L)=T(0)$, which is logical since tension is only defined up to a constant if the filament is straight. This demonstrates that our kinematic formulation is equivalent to the more standard line tension equation with free fiber boundary conditions $T(0)=T(L)=0$ \cite[Sec.~3]{maxian2021integral}.

\subsubsection{Continuum saddle point system}
\rev{To conclude the continuum formulation, we insert} the kinematic equations \eqref{eq:Kdef} and \eqref{eq:Vwork} into the mobility equation \eqref{eq:Motion} to obtain a saddle point system in the constraint forces $\V{\lambda}$ and degrees of freedom $\V{\alpha}$ 
\begin{gather}
\label{eq:SPC}
\begin{pmatrix} - \Lop{M} & \Lop{K} \\ 
\Lop{K}^*& \M{0} \end{pmatrix}
\begin{pmatrix} \V{\lambda} \\ \V{\alpha} \end{pmatrix} 
= \begin{pmatrix} \Lop{M}\left(-\delta \mathcal{E}^{(\kappa)}/\delta \XPoly +\fext \right)\\ \V{0} \end{pmatrix}.
\end{gather}
This system is written over a single fiber, where $\Lop{K}$ is the kinematic operator \eqref{eq:Kdef}, and $\V{\lambda}$ and $\V{\alpha}$ are the constraint forces and motions over that fiber. But it can also be extended to multiple fibers by slight abuse of notation, in which case $\Lop{K}$ is a block diagonal operator whose diagonal entries are the $\Lop{K}$ operators for each fiber, and $\V{\lambda}$ and $\V{\alpha}$ are the list of all constraint forces and kinematic motions. The only distinction between these two systems comes when we include nonlocal hydrodynamics. In that case the mobility operator $\Lop{M}$ is \emph{not} block diagonal, and couples the fibers together.

\subsection{Spatial discretization \label{sec:SpatDisc}}
\begin{figure}
\centering
\includegraphics[width=0.6\textwidth]{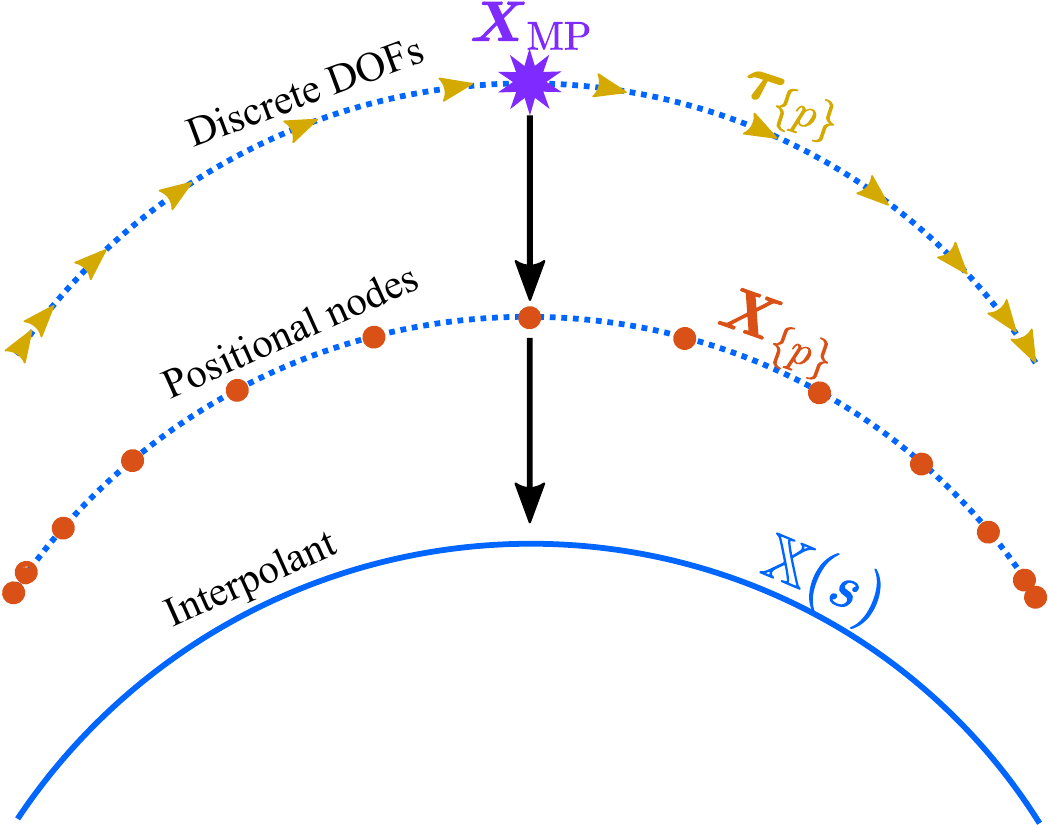}
\caption{\label{fig:Disc}Spectral discretization. The discrete degrees of freedom $\Xs$ and $\Xmp$ define a set of Chebyshev points $\V{X}$ and continuum interpolant $\XPoly(s)$.}
\end{figure}

The first decision that has to be made when discretizing in space is the representation of the fiber centerline $\XPoly(s)$. For slender filaments, the mobility $\Lop{M}$ is typically written as an integral over the fiber centerline. This leads us down the path of a spectral discretization \cite{ehssan17}, whereby we can use a set of collocation nodes to define a representation of $\XPoly(s)$ and develop quadrature schemes for integral equations. The smoothness of the fiber shapes we are interested in reinforces this idea, as spectral methods allow us to represent fibers with only a few modes, removing the need to resolve (irrelevant) high-order bending modes. However, it leaves inextensibility unclear: in spectral methods, the fiber shapes are ultimately polynomials, and it is impossible for the polynomial $\ds \XPoly \cdot \ds \XPoly-1$ to be zero \emph{everywhere} without being identically zero. And, if the space of polynomials to work in is unclear, how will we write an overdamped Langevin equation over a continuous constraint we do not understand?

Our discretization philosophy is to take advantage of the dual nature of a spectral collocation discretization: on the one hand, we have a set of discrete collocation points and tangent vectors which are the degrees of freedom in the simulation. But on the other, these points give us access to a continuum interpolant $\XPoly(s)$ which we can use to design efficient quadrature schemes for the mobility of slender filaments. Concretely, as shown in Fig.\ \ref{fig:Disc}, the fiber is defined by a set of interior tangent vectors $\Xs$ on a \emph{type 1} (not including the endpoints) Chebyshev grid of size $N$, together with the fiber midpoint $\Xmp$. From this, the positions of the fiber $\V{X}$ can be obtained \emph{exactly} (no integration error) on a grid of size $N_x=N+1$ via integration
\begin{equation}
\label{eq:Xmap}
\V{X} = \X \begin{pmatrix} \Xs \\ \Xmp \end{pmatrix}.
\end{equation}
The exact expression for $\X$ in terms of the constituent matrices is given in \cite[Sec.~2.1.1]{maxian2023bending}; here it suffices to note that $\X$ is an invertible matrix which maps the constrained degrees of freedom $\Xs$ to the fiber position (with the additional constant of the fiber midpoint). 

\subsubsection{$L^2$ inner products and the mobility}
While the continuum equations \eqref{eq:SPC} are formulated in terms of force densities $\V f$, in discrete variables we always need to work with \emph{forces} $\V F$, and the matrix $\Mfor$ which maps forces to velocity. To understand how to relate force to force density, we consider the $L^2$ inner product for power dissipated in the fluid 
\begin{equation}
\label{eq:Pdisc}
\mathcal{P} = \langle \V{f}, \V{U} \rangle_{L^2} = \V{f}^T \M{E}_{N_x \rightarrow 2N_x}^T \M{W}_{2N_x}  \M{E}_{N_x \rightarrow 2N_x}\V{U}:=\V{f}^T \Wt \V{U}:=\V{F}^T \V{U}
\end{equation}
In the case when $\V{f}$ and $\V{U}$ are polynomials defined on a Chebyshev grid of size $N_x$, the computation \eqref{eq:Pdisc} gives the power \emph{exactly} by first upsampling both $\V{f}$ and $\V{U}$ to a grid of size $2N_x$. \rev{The upsampling, which is denoted by the extension matrix $\M{E}_{N_x \rightarrow 2N_x}$, transforms the data from a coarse to fine grid by using the coarse grid values to compute a Chebyshev polynomial, which is evaluated pointwise on the finer grid to give the upsampled representation. After upsampling, integration is performed on the finer grid }by applying the integration weights (diagonal matrix $\M{W}_{2N_x}$). The last two equalities define the (symmetric positive definite) $L^2$ inner product weights matrix $\Wt=\Wt^T$, which we use to define a force $\V{F}=\Wt \V{f}$ from force density.

Let us now reconsider the power equation with the mobility \eqref{eq:Pmob}. Letting $\M{M}$ be the discretization of the matrix which maps force \emph{densities} to forces, the inner product \eqref{eq:Pmob} can be written in terms of forces as
\begin{equation}
\mathcal{P} = \langle \V{f}, \M{M}\V{f} \rangle_{L^2} =\langle \Wt^{-1} \V{F}, \M{M}\Wt^{-1} \V{F} \rangle_{L^2} = \V{F}^T \Wt^{-1}\Wt \M{M}\Wt^{-1}\V{F}:=\V{F}^T \Mfor \V{F}>0,
\end{equation}
which defines the SPD matrix $\Mfor=\M{M}\Wt^{-1}$ that acts on \emph{forces} to give velocity. 

\subsubsection{Bending force \label{sec:BForce}}
Once we define the $L^2$ inner product, discretization of the bending energy \eqref{eq:BendEnd} is straightforward. We compute the bending energy exactly via the $L^2$ inner product 
\begin{equation}
\label{eq:Edisc}
\mathcal{E}^{(\kappa)}\left[\V{X}\right]= \frac{\kappa}{2} \V{X}^T \left(\M{D}_{N_x}^2\right)^T \Wt \M{D}_{N_x}^2 \V{X}:=\frac{1}{2}\V{X}^T\M{L}\V{X}.
\end{equation}
The bending \emph{force} is now the derivative of the energy with respect to $\V{X}$, 
\begin{equation}
\label{eq:FKappa}
\V{F}^{(\kappa)} = -\frac{\partial \mathcal{E}^{(\kappa)}}{\partial \V{X}} = -\M{L}\V{X}.
\end{equation}
This type of discretization, which is common in immersed boundary methods \cite{peskin2002acta, zhu2002simulation}, enforces the boundary conditions naturally, in essence because minimizing the energy forces the boundary terms to be zero.

\subsubsection{Kinematic matrices}
The discretization of the kinematic operator $\Lop{K}$ is straightforward once we have a representation for the filament and $L^2$ inner product. We recall the definition of $\Lop{K}$ in \eqref{eq:Kdef} and the mapping from the tangent vectors to the nodes \eqref{eq:Xmap}. Letting $\M{C}\left[\Xs\right]$ be a matrix encoding $\M{C}\V{\Omega}=\Xs \times \V{\Omega}$ at each point, the product $\M{K}\V{\alpha}$ is discretized as 
\begin{equation}
\label{eq:KNp1}
\M{K}\V{\alpha}=\X \begin{pmatrix}-\M{C} & \M{0} \\ \M{0} & \M{I} \end{pmatrix} \begin{pmatrix} \V{\Omega}\\ \Ump \end{pmatrix},
\end{equation}
where now $\V{\Omega}$ represents a $3N$ vector of rotation rates for each point. The matrix $\M{K}$ is a square $3N_x \times 3(N+1)$ matrix, with a null space of size $N$ which contains the tangent vectors at each of the $N$ tangent vector nodes.

The constraint on the \emph{forces} $\V{\Lambda}:=\Wt \V{\lambda}$ follows from the virtual work principle \eqref{eq:Vwork}
\begin{equation}
P_\lambda = \langle \V{\lambda}, \M{K}\V{\alpha}\rangle_{L^2} = \V{\lambda}^T \Wt \M{K}\V{\alpha} = \V{\alpha}^T \M{K}^T \Wt \V{\lambda}=\V{\alpha}^T \M{K}^T \V{\Lambda}=\V 0.
\end{equation}
Since this must hold for all $\V{\alpha}$, we have the simple equation $\M{K}^T \V{\Lambda}=\V 0$.

\subsubsection{Discrete saddle point system}
Summing up, the discrete saddle point system is the continuum saddle point system \eqref{eq:SPC}, but rewritten in terms of \emph{forces} which contain the $L^2$ inner product weights, 
\begin{gather}
\label{eq:SPD}
\begin{pmatrix} - \Mfor & \M{K} \\ 
\M{K}^T& \M{0} \end{pmatrix}
\begin{pmatrix} \V{\Lambda} \\ \V{\alpha} \end{pmatrix} 
= \begin{pmatrix} \Mfor\left(-\M{L}\V{X}+\Fext \right)\\ \V{0} \end{pmatrix},
\end{gather}
and the deterministic dynamics can be obtained by eliminating $\V{\Lambda}$ to obtain 
\begin{equation}
\label{eq:Nhat}
\dt{\V{X}} = -\Nhat \M{L}\V{X}, \qquad \Nhat =\M{K}\left(\M{K}^T \Mfor^{-1} \M{K}\right)^\dagger \M{K}^T.
\end{equation}

\subsection{Langevin equation for Brownian fluctuations \label{sec:BM}}
Once the deterministic dynamics are known, we can formulate an overdamped Langevin equation which is in detailed balance with respect to the energy \eqref{eq:Edisc}, constrained on the inextensibility of the tangent vectors. This requires two pieces: first, in order to satisfy fluctuation-dissipation balance, the covariance of the noise must be proportional to $\Nhat$. In other words, coefficients of the noise should be proportional to $\Nhat^{1/2}$, where $\Nhat^{1/2}$ (which is not necessarily unique) satisfies 
\begin{equation}
\Nhat^{1/2} \left(\Nhat^{1/2} \right)^T = \Nhat. 
\end{equation}
Secondly, because the mobility $\Nhat$ depends on the fiber positions/tangent vectors, the noise is multiplicative. Thus, when the SDE for the motion of the filament collocation points is formulated in an Ito sense, stochastic drift terms are required to ensure detailed equilibrium \cite[Sec.~II.B.2]{delong2015brownian}. The Ito Langevin equation governing the evolution of the fiber centerline is given by \cite[Eq.~(35)]{maxian2023bending}
\begin{subequations}
\label{eq:Lang}
\begin{align}
\label{eq:ItoX}
\dt{ \V{X}} &= -\Nhat \M{L}\V{X} + k_B T \left(\div_{\V{X}} \cdot \Nhat\right) + \sqrt{2 k_B T}\Nhat^{1/2} \Wproc \\[4 pt]
\label{eq:kinetic}
& \eqdi-\Nhat \M{L}\V{X} +  \sqrt{2 k_B T}\Nhat \circ \Nhat^{-1/2}\Wproc,
\end{align}
\end{subequations}
where $\Wproc$ is a collection of standard i.i.d.\ Gaussian white noise processes (the formal derivatives of Brownian motion), the divergence of a matrix is defined as
\begin{equation}
\label{eq:DefDiv}
\left(\div_{\V{X}} \cdot \M{A}\right)_j:=\sum_k \frac{\partial  \M{A}_{jk}}{\partial \V X_k}
\end{equation}
and $\eqdi$ denotes equality in distributions of trajectories. The second equation \eqref{eq:kinetic} is the split Stratonovich-Ito \cite{delong2015brownian} or kinetic \cite{KineticStochasticIntegral_Ottinger} form, where the terms before the $\circ$ are evaluated
at the mid-point of the time interval (Stratanovich interpretation), while those after the $\circ$ are evaluated at the beginning of the time interval (Ito interpreation). When we present our temporal integrator for \eqref{eq:kinetic}, which was developed in \cite{maxian2023bending} and shown to be consistent with \eqref{eq:ItoX} therein, it will be helpful to keep \eqref{eq:kinetic} in mind.

At first glance, computing $\Nhat^{1/2}$ is quite difficult, since \eqref{eq:Nhat} suggests it requires the square root of a resistance solve. In fact, however, solving the saddle point system
\begin{gather}
\label{eq:SPB1}
\begin{pmatrix} - \Mfor & \M{K} \\ 
\M{K}^T& \M{0} \end{pmatrix}
\begin{pmatrix} \V{\Lambda} \\ \V{\alpha} \end{pmatrix} 
= \begin{pmatrix} \Mfor^{1/2} \V{\eta} \\ \V{0} \end{pmatrix},
\end{gather}
where $\gauss$ is an i.i.d.\ vector of standard normal random variables, gives 
\begin{align}
\nonumber
\dt{\V{X}}= \M{K}\V{\alpha} &=\Nhat \Mfor^{-1/2}\gauss=\Nhat^{1/2} \V{\eta},
\end{align}
where the last equality holds because 
\begin{align*}
\Nhat \Mfor^{-1/2}\left(\Nhat \Mfor^{-1/2}\right)^T &= \Nhat \Mfor^{-1} \Nhat =\M{K}\left(\M{K}^T \Mfor^{-1} \M{K}\right)^\dagger \M{K}^T\Mfor^{-1}\M{K}\left(\M{K}^T \Mfor^{-1} \M{K}\right)^\dagger \M{K}^T = \Nhat.
\end{align*}
Thus, solving a saddle point system with right hand side $\Mfor^{1/2}\gauss$ generates the noise $\Nhat^{1/2}\gauss$, and only a single resistance solve is necessary \cite[Sec.~II(B)]{delong2015brownian}.

In \cite{maxian2023bending}, we showed that the Langevin equation \eqref{eq:Lang} samples from the equilibrium Gibbs-Boltzmann distribution
\begin{align}
\label{eq:GBDist}
P_\text{eq}\left(\Xsbar\right)&= Z^{-1} \exp{\left(-\mathcal{E}^{(\kappa)}(\Xsbar)/ k_B T\right)} \prod_{p=1}^N \delta \left(\Xs_\Bind{p}^T \Xs_\Bind{p} - 1\right),
\end{align}
where the discrete energy $\mathcal{E}^{(\kappa)}$ is defined in \eqref{eq:Edisc}, $Z$ is a normalization constant, and the product of $\delta$ functions encodes the constraint that the tangent vectors $\Xs_\Bind{p}$ are independently uniformly distributed on the unit sphere for $p=1, \dots N$. This distribution is actually a postulate more than a fact; for worm-like chains, the tangent vectors are all equally spaced and thus independent. But for the spectral discretization shown in Fig.\ \ref{fig:Disc}, the tangent vectors near the fiber endpoints are quite close together, and it is doubtful that they are truly independent. Nevertheless, we showed in \cite{maxian2023bending} that \rev{Markov Chain Monte Carlo} (MCMC) samples from this equilibrium distribution converge to the theoretical end-to-end distribution of a free filaments as $N$ increases.

\section{Mobility \label{sec:Mob}}
This section discusses how to compute the action of the mobility $\Mfor$ and its square root $\Mfor^{1/2}$, both of which are necessary in Brownian dynamics simulations. To do this, we first define the velocity on filament $i$ as the sum of regularized Rotne-Prager-Yamakawa singularities over all other filaments. Then, we introduce a simple SPD ``reference'' mobility matrix which is based on global oversampling to compute the integrals accurately. We demonstrate through numerical experiments that the accuracy of the reference mobility is limited by the self term; some $1/\epsc$ points are required to resolve it accurately. It therefore becomes impractical for many filaments. As a result, the third part of this section seeks an SPD mobility where the self term is separated from the nonlocal terms, i.e., where the mobility can be written as 
\begin{equation}
\label{eq:MforSpl}
\Mfor = \Mfor_1+\Mfor_2,
\end{equation}
where both $\Mfor_1$ and $\Mfor_2$ are SPD matrices, and $\Mfor_1$ is a local operation while $\Mfor_2$ is a global operation. When the mobility is split in this way, the expectation of the covariance of the noise $\Mfor_1^{1/2} \Wproc_1 + \Mfor_2^{1/2} \Wproc_2$ (where $\Wproc_1$ and $\Wproc_2$ are independent white noise processes) is equal to $\Mfor$, as required to satisfy fluctuation-dissipation balance \cite{bao2018fluctuating, PSRPY}. We demonstrate that a combination of special quadrature for the local mobility, plus a ``fattening'' the nonlocal part of the mobility, gives a splitting with the desired properties and a cost independent of the fiber slenderness.

\subsection{The definition of the mobility}
In \cite{maxian2022hydrodynamics}, we motivated our choice of mobility (for a single fiber) by appealing to the immersed boundary or regularized singularity literature \cite{peskin1972flow, cortez2001method, fcm03, rpyOG,wajnryb2013generalization} . In our approach, the centerline of the fiber is modeled as a series of infinitely-many regularized delta functions (``blobs''), which in the Rotne-Prager-Yamakawa (RPY) case are surface delta functions on spheres of radius $\eps$. For two blobs, a regularized kernel $\M{M}_\text{RPY}$ is obtained by solving the Stokes equations with forcing given by the regularized delta function centered at a point $\V{y}$, then averaging the resulting fluid velocity field at another surface delta function centered at a point $\V{x}$ \cite{rpyOG, wajnryb2013generalization}
\begin{gather}
\label{eq:MbttRPY}
\M{M}_\text{RPY}
\left(\V{x},\V{y}; \eps \right)= \frac{1}{8 \pi \mu}\begin{cases}
 \left(\dfrac{\M{I}+\Rhat \Rhat}{R}+\dfrac{2\eps^2}{3} \dfrac{\M{I}-3\Rhat \Rhat}{R^3}\right)& R > 2\eps \\[8 pt]
\left(\left(\dfrac{4}{3\eps}-\dfrac{3R}{8\eps^2}\right)\M{I}+\dfrac{R}{8\eps^2} \Rhat \Rhat \right) & R \leq 2\eps
\end{cases}.
\end{gather}
When the blobs are well separated ($R > 2\eps$), the kernel is the sum of a Stokeslet and an $\mathcal{O}(\eps^2)$ multiple of the doublet, while when the blobs overlap there is a nonsingular correction which tends to the classical Stokes drag mobility $1/(6 \pi \eps \mu)$ as $\V{x} \rightarrow \V y$. Because of the symmetric nature of the ``spreading'' and ``averaging'' of force, the grand mobility matrix for a series of blobs is symmetric positive definite, which is a vital property for thermal fluctuations.

Once the regularized kernel for two blobs is defined, the velocity on fiber $i$ is defined as an integral of the regularized kernel
\begin{gather}
\label{eq:Ucont}
\ind{\V{U}}{i}(s) = \int_0^L \M{M}_\text{RPY}\left(\ind{\XPoly}{i}(s),\ind{\XPoly}{i}(s'); \eps \right) \ind{\V{f}}{i}(s')\, ds'\\ \nonumber +\sum_{j \neq i } \int_0^L \M{M}_\text{RPY}\left(\ind{\XPoly}{i}(s),\ind{\XPoly}{j}(s'); \epsh \right) \ind{\V{f}}{j}(s')\, ds',
\end{gather}
where \rev{the superscript $(i)$ indexes the fibers}, and we have used a distinct radius ($\eps$) for the self term than for the nonlocal terms ($\epsh$) for the sake of generality. While we will typically assume $\eps=\epsh$, we will later see that it is numerically useful to set $\epsh > \eps$. Indeed, for fibers that are well separated, the radius enters only in the doublet part of the regularized singularity \eqref{eq:MbttRPY}, and as such has a negligible effect. 

In the case of a single fiber (first line of \eqref{eq:Ucont}) \rev{with radius $\rc$}, we previously showed \cite[Appendix~A]{maxian2021integral} through asymptotic analysis that the choice of regularized radius
\begin{equation}
\label{eq:RadR}
\eps = \frac{e^{3/2}}{4}\rc 
\end{equation}
gives the same mobility as slender body theory \cite{krub, johnson,gotz2001interactions,ts04} to $\mathcal{O}(\eps)$. Unlike slender body theory, however, which suffers from ill-posedness on lengthscales less than $\eps$ (leading to negative eigenvalues), the RPY formulation has the advantage of giving a well-posed SPD mobility with a well-defined square root. Thus, using the RPY integral with the radius \eqref{eq:RadR} imparts all the benefits of slender body theory (in terms of asymptotic accuracy to the true three-dimensional Stokes equations \cite{mori2018theoretical}), without the associated ill-posedness. In our notation, we will often switch back and forth between the \emph{true} fiber radius $\rc$ and aspect ratio $\epsc=\rc/L$, and the RPY radius $\eps$ and aspect ratio $\epsRS=\eps/L$. 

\subsection{Upsampled mobility}
We first define a naive, obviously SPD, discretization of the mobility which will facilitate comparisons with more efficient methods. Let us set $\epsh=\eps$ in \eqref{eq:Ucont} and consider a simple way to \rev{compute} the velocity at a point on fiber $i$. Given a set of forces $\V{F}$ defined at Chebyshev grid points (either on one or many fibers), we form a reference mobility by applying the following steps
\begin{enumerate}
\item Obtain the force density $\V{f}=\Wt^{-1}\V{F}$. 
\item Extend the force density to an upsampled grid with $N_u$ points by resampling the Chebyshev interpolant, $\V{f}_u = \M{E}_u \V{f}$. 
\item Convert these force densities to forces on the upsampled grid  $\V{F}_u = \M{W}_u \V{f}_u$
\item Apply the grand RPY mobility matrix to obtain velocities on the upsampled grid $\V{U}_u = \Mfor_{\text{RPY}, u} \V{F}_u$. If there are multiple fibers, this is the step which encodes the nonlocal interactions.
\item Perform a least-squares projection in $L^2$ to obtain the velocity on the original $N$ point grid, $\V{U}=\Wt^{-1}\M{E}_u^T \M{W}_u \V{U}_u$. \rev{Here the $L^2$ weights matrix $\Wt$, which is computed in \eqref{eq:Pdisc} by going through a grid of size $2N_x$, is the same when computed on a grid of larger size, since any grid of size $N_u \geq 2N_x$ eliminates aliasing errors. Since least-squares projection in $L^2$ would use the matrix computed on a grid of size $N_u$, this step implies the restriction $N_u \geq 2N_x$.}
\end{enumerate}
Putting these steps together, we arrive at a \emph{reference oversampled mobility} matrix 
\begin{gather}
\label{eq:Mref}
\Mfor = \Mref:=\Wt^{-1} \M{E}_u^T \M{W}_u \Mfor_{\text{RPY}, u} \M{W}_u \M{E}_u \Wt^{-1},
\end{gather}
which has the obvious Cholesky factor
\begin{equation}
\label{eq:MrefCh}
\Mref^{1/2}=\Wt^{-1} \M{E}_u^T \M{W}_u \Mfor_{\text{RPY}, u} ^{1/2}.
\end{equation}
The most expensive part of applying the mobility is the application $\Mfor_{\text{RPY}, u}$ on the upsampled grid, which can be done in linear time on periodic domains using the positively-split Ewald method \cite{PSRPY} (as we do here), or in free space using the fast multipole method \cite{guan2018rpyfmm, yan2021kernel}. Likewise, we apply the square root of the mobility \eqref{eq:MrefCh} by using the positively-split Ewald method implemented in the UAMMD GPU library \cite{pelaez2025universally}. This applies $ \Mfor_{\text{RPY}, u} ^{1/2}$ on the \emph{upsampled} grid, with the premultipliers in \eqref{eq:MrefCh} downsampling to the original Chebyshev grid. 

The mobility \eqref{eq:Mref} is a robust, SPD mobility which can be applied in linear time. However, as we demonstrate in Appendix \ref{sec:fatNL}, the number of oversampling points required to resolve the \emph{self} interaction to a given accuracy scales as $1/\epsc$, while the number required to resolve the \emph{nonlocal} interactions is $\mathcal{O}(1)$. Consequently, making the mobility \eqref{eq:Mref} practical for slender filament suspensions requires replacing the block diagonal entries with a more efficient method.

\subsection{Special quadrature for the self term}
We previously developed \cite{maxian2022hydrodynamics} a special quadrature scheme for the self-RPY integral, which is summarized in Appendix \ref{sec:SQ}. This scheme takes the force density $\V{f}$ on the Chebyshev grid as input, and computes a velocity $\V{U}=\M{M}_\text{SQ}\V{f}=\M{M}_\text{SQ}\Wt^{-1} \V{F}$. Thus, the special quadrature mobility matrix can be written as $\M{M}_\text{SQ}\Wt^{-1}$, but this matrix is not even symmetric, much less positive definite. However, because it is localized to a single fiber, we can perform dense linear algebra operations on it, such as symmetrizing it and truncating its eigenvalues, thus resulting in 
\begin{equation}
\label{eq:Msqs}
\Msqs = \left[\frac{1}{2}\left(\M{M}\Wt^{-1}+\Wt^{-1} \M{M}^T\right)\right]_{\lambda^*}
\end{equation}
where $\left[\cdot\right]_{\lambda^*}$ denotes the operation of computing an eigenvalue decomposition, and setting all eigenvalues less than ${\lambda^*}$ as equal to ${\lambda^*}$. As shown in Appendix \ref{sec:fatNL}, truncating the eigenvalues causes a loss of accuracy, but only at low tolerances. While we previously used \eqref{eq:Mref} as a basis to set the eigenvalue threshold \cite{maxian2023bending}, here we find smaller tolerance to give higher accuracy without compromising robustness. We therefore use $\lambda^*=10^{-3}$ throughout this paper.

As shown in Appendix \ref{sec:fatNL}, the special quadrature scheme can correctly resolve the self velocity with a cost independent of $\epsc$. In fact, for moderate numbers ($N_x \approx 25$) of collocation points, we find that the special quadrature scheme is as accurate as using $N_u=2/\epsc$ oversampling points in the reference mobility \eqref{eq:Mref}. Thus, replacing the block diagonal terms of \eqref{eq:Mref} with the special quadrature matrix \eqref{eq:Msqs} allows for $\mathcal{O}(1)$ oversampling points \rev{(which are now only required for nonlocal terms)} as $\epsc$ decreases.

\subsection{Nonlocal mobility}
Therefore, a first guess for an $\epsc$-independent nonlocal mobility is a replacement of the block diagonal parts of \eqref{eq:Mref} with special quadrature
\begin{equation}
\label{eq:MFirst}
\Mfor = \Mref-\text{BDiag}\left\{\Mref\right\}+\text{BDiag}\left\{\Msqs\right\}.
\end{equation}
Because the oversampled mobility $\Mref$ is only used for the nonlocal velocity and not the self term, a given accuracy could be obtained using $N_u = \mathcal{O}(1)$ points (with respect to $\epsRS$). The special quadrature, which has cost independent of $\epsRS$, then gives the self velocity contributions. Thus the mobility \eqref{eq:MFirst} has cost independent of $\epsRS$. 

\begin{figure}
\centering
\includegraphics[width=0.49\textwidth]{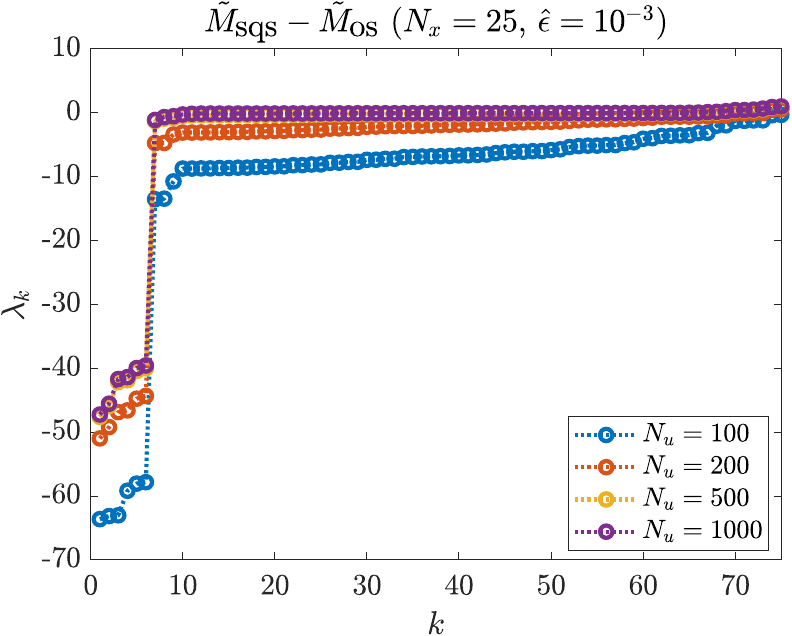}
\includegraphics[width=0.49\textwidth]{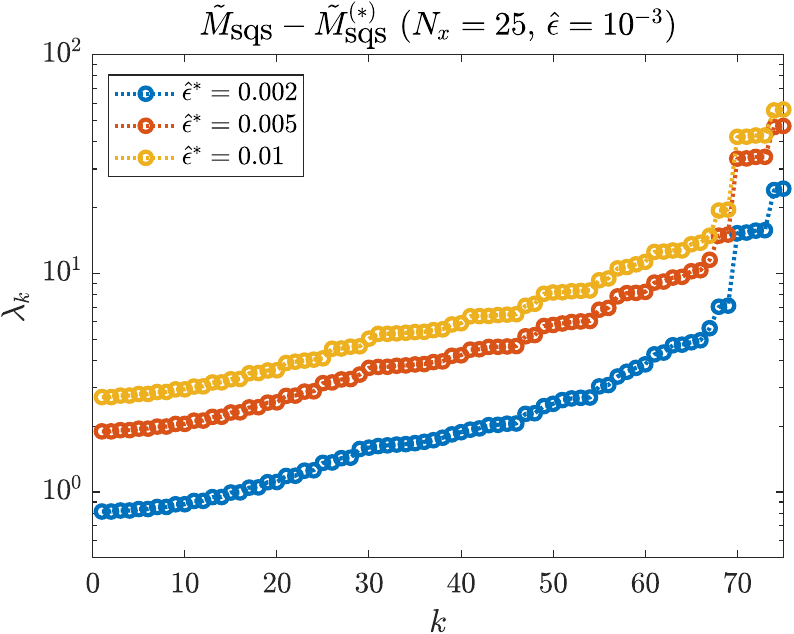}
\caption{\label{fig:EigValsSQ}Eigenvalues of possible corrections to the nonlocal mobility for $\epsRS=10^{-3}$. Left panel: correcting the mobility via \eqref{eq:MFirst}, by subtracting oversampled quadrature (with number of oversampled points $N_u$) and adding the special quadrature, gives a correction matrix with negative eigenvalues. Right panel: using a fatter radius $\eps^*$ for the nonlocal mobility via\ \eqref{eq:MNew}, then using the correction matrix $\Msqs-\Msqs^{(*)}$, always gives a symmetric positive definite correction. For these plots, $\Msqs$ is symmetrized, but the eigenvalues are not truncated (truncating them does not change the plots by eye). These plots are generated with a specific fiber ($N_x=25$ and $\epsRS=10^{-3}$), but are unchanged when we vary the shape.}
\end{figure}

While this is an appealing property, for the purposes of fluctuations we need to split the mobility into a sum of two SPD matrices (see \eqref{eq:MforSpl} and the discussion there). While $\Mref$ is SPD by construction, the block diagonal correction $\text{BDiag}\left\{\Msqs\right\}-\text{BDiag}\left\{\Mref\right\}$ is not SPD; in fact, the left panel of Fig.\ \ref{fig:EigValsSQ} shows that it has mostly negative eigenvalues. If we conceptualize the special quadrature as the limit $N_u \rightarrow \infty$, then we can think of this result as a consequence of how the fluid ``sees'' the fiber. In the extreme limit $N_u=1$, the mobility becomes that of a sphere, $\Mref \sim 1/\epsRS$, while in the limit $N_u \rightarrow \infty$, the mobility becomes filamentous, $\Msqs \sim \log{\epsRS}$. Thus, the difference matrix $\text{BDiag}\left\{\Msqs\right\}-\text{BDiag}\left\{\Mref\right\}$ will have almost exclusively negative eigenvalues (an example for a filament with $N_x=25$ is shown in Fig.\ \ref{fig:EigValsSQ}).

We propose counteracting this drop in the eigenvalues via a smoothing of the oversampled mobility matrix. Returning to the generalized mobility \eqref{eq:Ucont}, we set a larger \rev{(``fatter'')} RPY radius $\epsh > \eps$, which is independent of the true radius $\eps$. We then propose changing the \emph{nonlocal} kernel to use this larger radius, denoting this by $\Mref^{(*)}$. That is, $\Mref^{(*)}$ is the reference mobility \eqref{eq:Mref}, but with the larger radius $\epsh$. We then set the total mobility 
\begin{equation}
\label{eq:MNew}
\Mfor = \Mref^{(*)}-\text{BDiag}\left\{\Msqs^{(*)}\right\}+\text{BDiag}\left\{\Msqs\right\},
\end{equation}
i.e., we compute the oversampled fattened mobility for the nonlocal velocity, subtract special quadrature on the fattened mobility (to remove the incorrect self velocity), then add the special quadrature with the correct value of $\eps$ (to add the correct self velocity). This formulation relies on three key ideas:
\begin{enumerate}
\item The correct aspect ratio is only important for the self velocity, and not the nonlocal term. As discussed previously, for well-separated fibers, the nonlocal integrals take the form of a Stokeslet + $\mathcal{O}(\eps^2) \times$ the doublet. Thus, for small aspect ratios, the thickness of the fiber is secondary to the line of Stokeslets, and artificially thickening the fibers in the nonlocal velocity calculation makes little difference for the error (see Fig.\ \ref{fig:NLStar} in Appendix\ \ref{sec:fatNL}). This was recognized previously in methods that only used lines of Stokeslets for the nonlocal velocity \cite{ehssan17, koshakji2023robust}. Our approach introduces an $\mathcal{O}(\left(\epsRSh\right)^2)$ asymptotic error for well-separated filaments, where $\epsRSh=\epsh/L$ is the thickened aspect ratio.
\item  Because $\epsRSh$ will be relatively large, we can afford a larger number of oversampled points to compute the fattened mobility, $N_u=1/\epsRSh$. As a result, special quadrature on the fattened kernel will be a good approximation to the block diagonal terms of $\Mref^{(*)}$, and the difference $\text{BDiag}\left\{\Mref^{(*)}\right\} -\text{BDiag}\left\{\Msqs^{(*)}\right\}$ is near zero. As shown in Fig.\ \ref{fig:SelfErFat}, the addition and subtraction of two different numerical representations of the self mobility on the fattened fibers contributes an error of at most 1\%, and can be controlled by increasing the number of upsampled points.
\item Since the difference $\text{BDiag}\left\{\Msqs\right\} - \text{BDiag}\left\{\Msqs^{(*)}\right\}$ represents the same numerical scheme with two different values of $\eps$, the eigenvalues of the correction matrix should be positive if the correct radius $\eps < \epsh$. This is demonstrated in the right panel of Fig.\ \ref{fig:EigValsSQ}, where we plot the eigenvalues of the correction matrix $\Msqs-\Msqs^{(*)}$ for a single fiber with difference choices of $\epsh$ (they approach zero as $\epsh \rightarrow \eps$).
\end{enumerate}
Putting these observations together, the mobility \eqref{eq:MNew} gives the self mobility to numerical precision, while approximating the nonlocal mobility in a controllable way (approximations can be checked by decreasing $\epsh$ towards $\eps$), and, most importantly, is the sum of two SPD components. Thus, we can generate $\Mfor^{1/2} \Wproc$ by taking (c.f.\ \eqref{eq:MforSpl}) 
\begin{equation}
\left( \Mref^{(*)}\right)^{1/2}\Wproc_1+\left(\text{BDiag}\left\{\Msqs\right\}-\text{BDiag}\left\{\Msqs^{(*)}\right\}\right)^{1/2}\Wproc_2.
\end{equation}
The first square root is applied using \eqref{eq:MrefCh} and the positively-split Ewald method, while the second is localized to each fiber and is computed by dense linear algebra. 

In practice, we find the correction matrix $\Msqs-\Msqs^{(*)}$ can occasionally have negative eigenvalues when the time step size does not correctly resolve thermal fluctuations. While we should never be running in this regime, we prevent our code from crashing in this case by truncating the eigenvalues at $\lambda^*=10^{-3}$, as discussed in \eqref{eq:Msqs}.

\section{Steric repulsion \label{sec:Sterics}}
The change in the RPY mobility \eqref{eq:MbttRPY} on $R \leq 2\eps$ is unphysical, since fibers immersed in a common medium should never overlap one another. Of course, this happens quite a lot in our simulations, since the filament geometries are slender and underresolved. This section offers a way to keep filaments well-separated by using a soft Gaussian force between nearly-touching fibers. This approach is in no way superior to the recently-developed (implicit) methods that guarantee no overlap between particles \cite{broms2023barrier, yan2019computing, yan2022toward, ferguson2021intersection}, but it is more practical for our purposes, since it uses an explicit method (one \rev{force calculation} per time step) and does not become singular when fibers overlap (which will happen randomly due to Brownian motion). Simulations in Section \ref{sec:Bundling} show that our approach is effective in keeping the filaments apart without introducing a strong restriction on the time step size. The focus then becomes how to efficiently determine near-contacts and evaluate forces.

As in \cite{popov2016medyan}, we propose a double-integral steric interaction energy between two fibers $i$ and $j$,
\begin{gather}
\label{eq:EDCInt}
\ind{\Pot}{ij}=\int_0^{L} \int_0^{L} \hat{\Pot}\left(r\left(\ind{s}{i},\ind{s}{j}\right)\right) \, d\ind{s}{i} \, d\ind{s}{j},\\ \nonumber
r\left(\ind{s}{i},\ind{s}{j}\right)=\norm{\ind{\XPoly}{i}\left(\ind{s}{i}\right)-\ind{\XPoly}{j}\left(\ind{s}{j}\right)},
\end{gather}
where $\hat{\Pot}$ is the potential density function (units energy/area), which we assume is compactly supported on $0 \leq r \leq \rmax$. For the potential density function, we use an error function so that the force will ultimately be a Gaussian with standard deviation $\delta$
\begin{gather}
\label{eq:PotSter}
\hat{\Pot}(r) = \frac{{\Pot}_0}{\rc}\text{erf}\left(r/(\delta \sqrt{2})\right)\\[3 pt] \nonumber
\frac{d\hat{\Pot}}{dr} = \frac{{\Pot}_0}{\rc \delta}\sqrt{\frac{2}{\pi}} \exp{\left(-r^2/\left(2\delta^2\right)\right)}
\end{gather}
In this potential density, we have two free parameters: $\delta$, which is the scale on which the force decays, and ${\Pot}_0$, which controls the magnitude of the forces (and therefore the temporal stiffness) and has units of energy per length (force). We set $\delta=\rc$ to make the steric repulsion short ranged, and truncate the Gaussian at four standard deviations, $\rmax=4\delta$. Because the force $d\hat{\Pot}/dr$ is typically multiplied by weights on both interacting filaments with size $\mathcal{O}(\rc)$, choosing ${\Pot}_0$ constant with $\delta \sim \rc$ gives a force which is $\mathcal{O}(1)$ with respect to the fiber aspect ratio, thus ensuring a relatively constant required time step size with decreasing aspect ratio. Based on these considerations, we set $\Pot_0=B k_B T$, with $B$ a simulation-dependent constant which controls the magnitude of the force.

\rev{While it is possible to differentiate the energy\ \eqref{eq:EDCInt} in continuum first, this results in a force density on fiber $i$ that is a \emph{single integral} over fiber $j$. In this case, we could only detect contacts close to the chosen Chebyshev discretization points on fiber $i$. Rather, we will discretize the energy first, thus ensuring we actually compute a double integral, detecting contacts between all pairs of fiber regions. We then differentiate the discrete energy expression to get force. The details are as follows:} if the Chebyshev points used to discretize $\XPoly(s)$ are given by $\V{X}$, then the upsampled points will be denoted by $\Xu = \M{E}\V{X}$. The double integral can then be evaluated and differentiated via
\begin{align}
\Pot =& \sum_k \sum_j \hat{\Pot}\left(\norm{\Xu_{\Bind{k}}-\Xu_{\Bind{j}}}\right) w_k w_j\\ \nonumber
=& \sum_k \sum_j \hat{\Pot}\left(\norm{\M{E}_{kp}\V{X}_\Bind{p}-\M{E}_{jq}\V{X}_{\Bind{q}}}\right) w_k w_j\\
\label{eq:sumForceUp}
\ind{\V{F}}{i}_{\Bind{a}} = -\frac{\partial \Pot}{\partial \ind{\V{X}}{i}_{\Bind{a}}}=-&\sum_k \sum_j \frac{\partial \hat{\Pot}}{\partial r}\left(r_{kj}\right) \widehat{\V{r}}_{kj} \M{E}_{ka} w_k w_j,
\end{align}
where $\V{r}_{kj} = \Xu_{\Bind{k}}-\Xu_{\Bind{j}}$. The last equation gives the force at Chebyshev node $a$ on fiber $i$, and is a function of the integration weights $w_k$ and $w_j$ of points $k$ and $j$ on the upsampled grid. 

In our algorithm, we will \emph{always} use\ \eqref{eq:sumForceUp} to compute forces. The freedom we have is how to choose the oversampled points $\Xu$ (and therefore the matrix entry $\M{E}_{ka}$ in\ \eqref{eq:sumForceUp}). We present two possibilities for doing this: global upsampling (guaranteed accuracy, but inefficient in the limit $\epsc \rightarrow 0$), as well as a second algorithm based on selective upsampling for pieces of the fibers that are nearly in contact. This second algorithm achieves our goal of a cost independent of $\epsc$. In Appendix \ref{sec:StericAppen}, we analyze the differences between the two algorithms, demonstrating that the accuracy of the segment-based algorithm is comparable to the uniform-point-based algorithm with $1/\epsc$ points.

\subsection{Global uniform point resampling \label{sec:UnifRS}}
Let us suppose that we resample the fiber $\XPoly(s)$ at $N_u$ uniform points to form the vector $\V{X}^{(u)}$. We will typically use $N_u=1/\epsc=1/\delta$, but more upsampling might be necessary to achieve higher accuracy. Thus the points we choose are at arclength coordinates $\V s=0, \Delta s_u, \dots, L$, where $\Delta s_u =L/(N_u-1)$ is the spacing. The corresponding weights are $\V w=\Delta s_u/2, \Delta s_u, \dots \Delta s_u, \Delta s_u/2$, so that the first and last point have a weight of 1/2, in accordance with the trapezoid rule. Then a simple algorithm to evaluate\ \eqref{eq:sumForceUp} is as follows:
\begin{enumerate}
\item Resample every fiber at $N_u$ points.
\item Perform a (linear-cost) neighbor search to determine pairs of points for which the potential function $\partial \hat{\Pot}/\partial r$ is nonzero to within a certain tolerance, i.e., find pairs of points a distance $\rmax$ or less apart \rev{(excluding points on the same fiber that are closer than $\rmax$ when the fiber is straight)}. 
\item For each pair of points with $\norm{\Xu_\Bind{k}-\Xu_\Bind{j}} < \rmax$ (determined in step 2), compute the corresponding term in the sum\ \eqref{eq:sumForceUp} and add it to obtain the forces on the Chebyshev nodes. 
\end{enumerate}
This is obviously an extremely simple algorithm, but it can become costly as the fibers get more slender, since resolving all of the potential contacts requires a large number uniform points (the quadratures are also only second-order accurate, which could limit accuracy). In the next section, we aim for an algorithm which is independent of $\epsc$. 

\subsection{Segment-based algorithm \label{sec:SterSeg}}
\begin{figure}
\centering
\includegraphics[width=\textwidth]{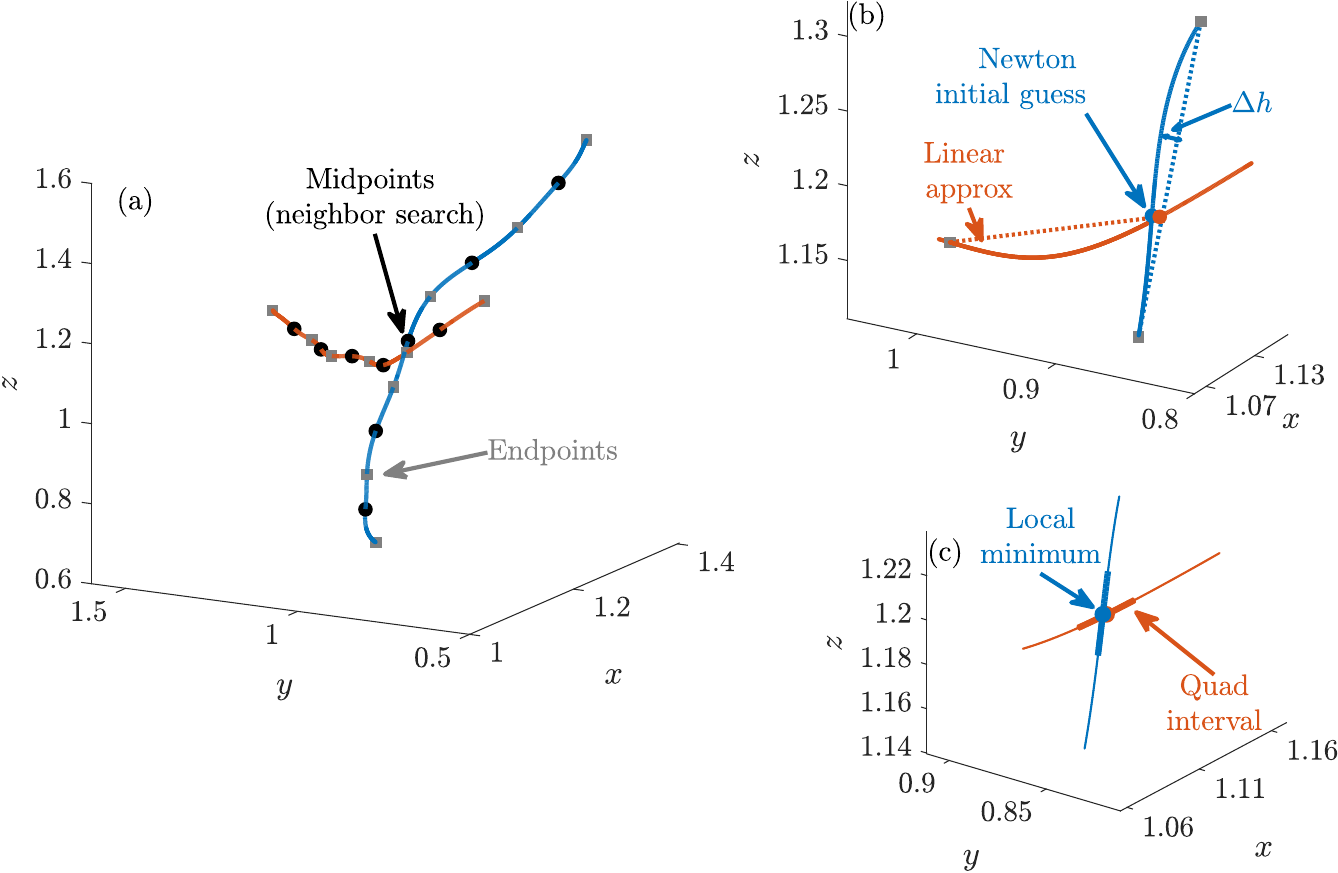}
\caption[Schematic of segment-based steric interaction algorithm, which has cost independent of $\epsc$]{\label{fig:SegSch}Schematic of segment-based steric interaction algorithm. (a) We divide the fiber into pieces and use the midpoint of each piece (black circles) to identify possible close segments. (b) We then approximate each piece using a straight segment (dashed lines), find the minimum distance between the two segments, and use that as an initial guess for Newton's method. (c) Newton's method finds the local minimizer of the distance on the two fibers, around which we draw quadrature intervals (thick lines).}
\end{figure} 

The following algorithm theoretically achieves a cost that is linear in the number of fibers (at constant density) and independent of $\delta \sim \rc$. It is based on dividing the fiber into $\Nseg$ (curved) fiber pieces of length $\Lseg$, then approximating those pieces as straight segments for the purposes of a cheap distance estimation. If the segments are close enough, we perform a nonlinear Newton solve to obtain the closest points on the corresponding fibers. The number of fiber pieces in this algorithm is chosen to minimize computational cost, and should be independent of $\epsc$. A schematic of the algorithm is shown in Fig.\ \ref{fig:SegSch}, and a more detailed explanation of the steps is as follows.

\begin{enumerate}
\item We first sample each fiber at $\Nseg$ points, corresponding to the midpoints of the fiber pieces (circles in Fig.\ \ref{fig:SegSch}(a)). Then, we perform a neighbor search over all the midpoints using the cutoff distance $r_\text{cut}=\Lseg+\rmax$. This gives a list of potentially interacting segments, from which we exclude the self-interaction (there is no steric force between a segment and itself). 
\item For each pair of fiber pieces, we approximate curved fiber pieces as straight segments (dashed lines in Fig.\ \ref{fig:SegSch}(b)) and solve a quadratic equation to determine the closest point of interaction between the straight segments. Let us denote this distance $\ind{d}{ij}_\text{segs}$, with the corresponding closest points on the fibers given by $\ind{s}{i}_\text{seg}$ and $\ind{s}{j}_\text{seg}$.
\item We next determine if $\ind{d}{ij}_\text{segs} < \rmax+\ind{\Delta h}{i}+\ind{\Delta h}{j}$, where $\ind{\Delta h}{i}$ represents the distance the fiber piece $\ind{\XPoly}{i}(s)$ curves from the straight segment we use to approximate it (see Section\ \ref{sec:ClosePts}). If $d_\text{segs} \geq \rmax+\ind{\Delta h}{i}+\ind{\Delta h}{j}$, then we exit, as the corresponding fiber pieces could never be close enough to interact. 
\item We now proceed to the case when the fiber pieces could be close enough to interact. In this case, we solve the minimization problem
\begin{equation}
\label{eq:NewtonMin}
\left(\ind{s}{i}_*,\ind{s}{j}_*\right)= \argmin_{\begin{matrix} 0 \leq \ind{s}{i} \leq L \\ 0 \leq \ind{s}{j} \leq L \end{matrix}} \norm{\ind{\XPoly}{i}\left(\ind{s}{i}\right)-\ind{\XPoly}{j}\left(\ind{s}{j}\right)}^2,
\end{equation}
with a modified projected Newton's method, with initial guess $\ind{s}{i}=\ind{s}{i}_\text{seg}$ and $\ind{s}{j}=\ind{s}{j}_\text{seg}$. Our methodology is described in Appendix\ \ref{sec:Newton}, and is based on the projected Newton methods discussed in \cite{bertsekas1982projected, kim2010tackling}.
\item We then approximate the fibers as quadratic curves around the local minimum, and solve $\norm{\ind{\XPoly}{i}\left(\ind{s}{i}_*+\ind{\D s}{i}\right)-\ind{\XPoly}{j}\left(\ind{s}{j}_*+\ind{\D s}{j}\right)}^2=\rmax^2$ to obtain intervals $\ind{S}{i}$ and $\ind{S}{j}$ on which we need to integrate the energy\ \eqref{eq:EDCInt}. 
\item Because there are multiple segments on each fiber, it is possible that the interval associated with a particular minimum might overlap with another interval for the same fiber pair. Because of that, we assemble a list of fiber pairs $i$ and $j$ and intervals $\ind{S}{i}$ and $\ind{S}{j}$ from step 5, over which we do the integral\ \eqref{eq:sumForceUp} by taking a kind of ``union.'' Because we are in two dimensions, this procedure is not entirely straightforward, and so we describe details of it in Appendix\ \ref{sec:Unions}.
\item{\label{step7}} Finally, we put a grid of Gauss-Legendre points over each interval. Letting $L_\text{int}$ be the length of an interval, the number of points we use is
\begin{equation}
\label{eq:NGL}
N_\text{GL} = \left\lfloor N_{\delta}\frac{L_\text{int}}{\delta}\right\rfloor+1,
\end{equation}
where $\delta$ is the standard deviation of the Gaussian potential\ \eqref{eq:PotSter} and $N_\delta$ is a constant which gives the number of Gauss-Legendre points we use per standard deviation. We then apply the formula\ \eqref{eq:sumForceUp} over the two Gauss-Legendre grids. To implement this step efficiently, we precompute a maximum number of Gauss-Legendre points by setting $L_\text{int}=L$ in\ \eqref{eq:NGL}. We then pretabulate all possible points and weights on $[-1,1]$ for $N_\text{GL}=1,\dots N_\text{GL}^\text{(max)}$. Then, for a given $L_\text{int}$, we compute $N_\text{GL}$ using\ \eqref{eq:NGL}, and look up the points and weights in the table, scaling appropriately by $L_\text{int}$.
\end{enumerate}

To make the cost of the algorithm formally independent of $\epsc$, the neighbor search in step 1 must be on a number of segments that is independent of $\epsc$. The cost of step 2 is then independent of $\epsc$ because we are minimizing segment distances. Steps 3--5 are the most difficult, but have cost independent of $\epsc$ because we use nonlinear solves to obtain the closest points. And because $\rmax \sim \delta \sim \rc$, the grid in steps 5--7 will have a thickness that scales with $\epsc$, so the cost of the quadrature in step 7 is independent of $\epsc$. 

Ultimately, the efficiency of the algorithm (and its favorability over the one presented in Section\ \ref{sec:UnifRS}) comes down to the expense of the neighbor search vs.\ root finding and quadrature. If the root finding is expensive, then it makes sense to just use $1/\epsc$ fiber pieces (which become spheres), in which case the algorithm becomes equivalent to Section\ \ref{sec:UnifRS}. If the root finding and quadrature are cheaper \rev{(which is typically the case as $\epsc \rightarrow 0$)}, then we can use less points in the neighbor search and this segment-based algorithm might be more efficient than the uniform-point algorithm.

\section{Temporal integration \label{sec:TInt}}
Now that we have reviewed the spatial equations, it remains to discretize in time. Because the applications for these methods are mainly biological, the focus of temporal integration is on robustness, and not necessarily accuracy. For this reason, the temporal integrators we present here are first-order accurate (see \cite{maxian2021integral} for a second-order version in the deterministic case). The deterministic integrator is straightforward, and amounts to a first-order backward Euler discretization of \eqref{eq:SPD}, with the caveat that the final update to the fiber is a (nonlinear) rotation of the tangent vectors to keep the dynamics on the constraint, followed by integration to obtain the new fiber positions. The integrator for Brownian fluctuations is based on using the Brownian fluctuations to generate a trial step to the midpoint, then solving the equations using the midpoint configuration. This procedure captures a key part of the drift term \eqref{eq:ItoX}, precisely because it is consistent with the midpoint interpretation in \eqref{eq:kinetic}.

\subsection{Deterministic methods}
Considering the linear system \eqref{eq:SPD}, which must be solved at every time step, let us introduce the time step index $n$ and approximate
\begin{equation}
\label{eq:Xnp1St}
\tdisc{\V{X}}{n+1,*} = \tdisc{\V{X}}{n}+\Delta t \tdisc{\M{K}}{n}\tdisc{\V{\alpha}}{n}. 
\end{equation}
Because bending forces are stiff, we treat them implicitly using the linearized backward Euler method. This gives the linear system
\begin{gather}
\tdisc{\begin{pmatrix} - \Mfor & \M{K} \\ 
\M{K}^T& \M{0} \end{pmatrix}}{n}
\tdisc{\begin{pmatrix} \V{\Lambda} \\ \V{\alpha} \end{pmatrix}}{n}
= \begin{pmatrix} \tdisc{\Mfor}{n}\left(-\M{L}\tdisc{\V{X}}{n+1,*}+\tdisc{\Fext}{n} \right)\\ \V{0} \end{pmatrix},
\end{gather}
Substituting \eqref{eq:Xnp1St} gives a linear system for $\V \alpha$ and $\V{\Lambda}$ in terms of $\tdisc{\V{X}}{n}$. 
\begin{gather}
\label{eq:SPDBE}
\tdisc{\begin{pmatrix} - \Mfor & \M{K}+\Delta t \Mfor \M{L} \M{K} \\ 
\M{K}^T& \M{0} \end{pmatrix}}{n}
\tdisc{\begin{pmatrix} \V{\Lambda} \\ \V{\alpha} \end{pmatrix}}{n}
= \begin{pmatrix} \tdisc{\Mfor}{n}\left(-\M{L}\tdisc{\V{X}}{n}+\tdisc{\Fext}{n} \right)\\ \V{0} \end{pmatrix}.
\end{gather}
Now, depending on the form of the mobility, the solution strategies for this system vary. If $\Mfor$ is a block diagonal matrix (localized to each fiber without any inter-fiber interactions), this system can be solved using dense linear algebra (i.e., by forming $\Mfor$ as a dense matrix and inverting it directly). If, however, $\Mfor$ is a dense matrix, encoding interactions of all fibers with all other fibers, iterative methods are required. In this case, we split the mobility matrix into a local and nonlocal part, 
\begin{align}
\label{eq:Msplit}
\Mfor = \ML+&\MNL \\
\label{eq:LocDiag}
= \text{BDiag}\left\{\Mfor \right\} +& \left(\Mfor -  \text{BDiag}\left\{\Mfor \right\}\right).
\end{align}
We then solve the system \eqref{eq:SPDBE} using a block diagonal preconditioner,
\begin{equation}
\label{eq:Precond}
\M{P}=
\begin{pmatrix}
    -\ML & \M{K}+\D t \ML \M{L}\M{K} \\[2 pt]
    \M{K}^T & \M{0}
    \end{pmatrix},
\end{equation}
This preconditioner ought to be effective since $\ML$ typically dominates $\MNL$ for slender filaments. Note that while other options for the mobility splitting are possible, including making $\ML$ itself diagonal, here we will restrict to the case of\ \eqref{eq:LocDiag}, where the local mobility on fiber $i$ is just the $i$th diagonal block of the matrix $\Mfor$. 

\subsubsection{Time-lagging the nonlocal forces}
For most fiber suspensions, the local behavior dominates the dynamics, and so we can actually obtain similar stability behavior by \emph{time-lagging} the non-local parts of the mobility. In particular, we can only treat the bending and constraint forces implicitly for the \emph{local} parts of the mobility, which gives rise to the linear system 
\begin{gather}
\begin{small}
\label{eq:BEnL}
\tdisc{\begin{pmatrix} - \ML & \M{K}+\Delta t \ML \M{L}\M{K} \\ 
\M{K}^T& \M{0} \end{pmatrix}}{n}
\tdisc{\begin{pmatrix} \V{\Lambda} \\ \V{\alpha} \end{pmatrix}}{n}
= \begin{pmatrix} \left(\tdisc{\ML}{n}+\tdisc{\MNL}{n}\right)\left(-\M{L}\tdisc{\V{X}}{n}+\tdisc{\Fext}{n}\right)+\tdisc{\MNL}{n}\left(\tdisc{\V{\Lambda}}{n-1}\right)\\ \V{0} \end{pmatrix}.
\end{small}
\end{gather}
To the extent that this approach is stable, it significantly reduces the cost, since only one nonlocal hydrodynamic evaluation is required per time step. When the suspension is dense enough, the dynamics are unstable, but we can use the solutions above, denoted by $\tdisc{\widetilde{\V{\Lambda}}}{n}$ and $\tdisc{\widetilde{\V{\alpha}}}{n}$, to obtain a modified system of equations. The key is to rewrite the locally implicit system\ \eqref{eq:BEnL} as
\begin{equation*}
\tdisc{\M{K}}{n} \tdisc{\widetilde{\V \alpha}}{n} = \left(\tdisc{\ML}{n}+\tdisc{\MNL}{n}\right) \left(-\M{L}\tdisc{\V{X}}{n}-\Delta t \M{L} \tdisc{\M{K}}{n}\tdisc{\widetilde{\V \alpha}}{n}+\Fext+\tdisc{\widetilde{\V{\Lambda}}}{n}\right)+\tdisc{\MNL}{n} \left(\Delta t \M{L} \tdisc{\M{K}}{n}\tdisc{\widetilde{\V \alpha}}{n}-\tdisc{\widetilde{\V{\Lambda}}}{n}+\tdisc{\V{\Lambda}}{n-1}\right),
\end{equation*}
then subtract this from the fully implicit system\ \eqref{eq:SPDBE}, which gives the residual form of the saddle-point system
\begin{gather}
\label{eq:residBE}
    \begin{pmatrix}
    -\Mfor& \M{K}+\D t \Mfor \M{L}\M{K} \\[2 pt]
    \M{K}^T & \M{0}
    \end{pmatrix}_{n}
    \begin{pmatrix} 
    \Delta \V{\Lambda}\\[2 pt]
    \Delta \V{\alpha}
    \end{pmatrix}_n =  
    \begin{pmatrix} 
    \tdisc{\MNL}{n}\left(-\Delta t \M{L}\tdisc{\M{K}}{n}\tdisc{\widetilde{\V{\alpha}}}{n}+\tdisc{\widetilde{\V{\Lambda}}}{n}-\tdisc{\V{\Lambda}}{n-1}\right)\\[2 pt]
   \V{0}
   \end{pmatrix}
\end{gather}
to be solved using GMRES for the perturbations $\Delta \tdisc{\V{\Lambda}}{n} = \tdisc{\V{\Lambda}}{n} - \tdisc{\widetilde{\V{\Lambda}}}{n}$ and $\Delta \tdisc{\V{\alpha}}{n} = \tdisc{\V{\alpha}}{n} - \tdisc{\widetilde{\V{\alpha}}}{n}$. 

There are two caveats here. First, obviously the time-lagging procedure does not work at $t=0$, so we must solve the fully implicit system \eqref{eq:SPDBE} there. Second, the solution of the residual system \eqref{eq:residBE} is only required for \emph{stability}, and not accuracy. As such, we do not have to solve the system to a low tolerance. Instead, we empirically set a maximum number of GMRES iterations, increasing this number until we obtain stability.

\subsubsection{Nonlinear update}
No matter the precise linear system being solved, once we solve for $\tdisc{\V{\alpha}}{n}=\left(\tdisc{\V{\Omega}}{n},\tdisc{{\Ump}}{,n}\right)$, we have access to the discrete tangent vector rotation rates $\tdisc{\V{\Omega}}{n}$. As such, we update the fibers by rotating the tangent vectors
\begin{equation}
\label{eq:RotUpdate}
\tdisc{\Xs}{n+1}=\text{rotate}\left(\tdisc{\Xs}{n}, \tdisc{\V{\Omega}}{n} \D t\right),
\end{equation}
using the Rodrigues rotation formula \cite[Eq.~(111)]{maxian2021integral}. We then update the midpoint $\Xmp$ via $\tdisc{{\Xmp}}{,n+1}=\tdisc{{\Xmp}}{,n}+\D t \tdisc{{\Ump}}{,n}$, and obtain $\tdisc{\V{X}}{n+1}$ by applying the matrix $\X$ defined in \eqref{eq:Xmap}. We will denote this nonlinear update, which preserves the length of the tangent vectors exactly, as $\text{RotateAndIntegrate}\left(\tdisc{\V{X}}{n},\tdisc{\V{\alpha}}{n}\Delta t \right)$.

\subsection{Brownian fluctuations}
The temporal integrator for Brownian fluctuations requires additional terms in the saddle point system \eqref{eq:SPD} to account for the Brownian velocity and drift terms in the Langevin equation \eqref{eq:ItoX}. In \cite{maxian2023bending}, we broke the drift term into three pieces, showing that one piece can be obtained by rotating the tangent vectors at the end of the time step, while the second piece can be obtained by solving the saddle point system \eqref{eq:SPD} at a kind of ``midpoint'' configutation (consistent with the kinetic formulation \eqref{eq:kinetic}). The third piece is obtained by adding a drift of the mobility $\Mfor$ to the right hand side, computed using a random finite difference \cite{delong2014brownian}. 

The full temporal integration scheme for Brownian filaments is as follows \cite{maxian2023bending}
\begin{enumerate}
\item Generate a Brownian velocity increment
\begin{equation}
\label{eq:UB}
\tdisc{\V{U}^{(B)}}{n}=\sqrt{\frac{2 k_B T}{\D t}}\Mfor^{1/2}\tdisc{\gauss}{n},
\end{equation}
where $\tdisc{\gauss}{n}$ is a vector of i.i.d.\ standard Gaussian $\mathcal{N}(0,1)$ random variables.
\item Use the Brownian velocity increment to propose a rotation rate for the tangent vectors 
\begin{equation}
\label{eq:OmegaTilde}
\tdisc{\V{\alpha}}{n,*} =\tdisc{\M{K}}{n}^\dagger \tdisc{\V{U}^{(B)}}{n}
\end{equation}
\item Use the proposed rotation rates to generate a ``midpoint'' configuration of the filaments 
\begin{equation}
\tdisc{\V{X}}{n+1/2,*}=\text{RotateAndIntegrate}\left(\tdisc{\V{X}}{n},\tdisc{\V{\alpha}}{n.*}\Delta t/2 \right)
\end{equation}

\item Compute an additional drift velocity (``mobility drift,'' MD) in one of two ways. If the mobility is local and can be formed as a dense matrix, set 
\begin{equation}
\label{eq:UMD}
\tdisc{\V{U}}{n}^\text{(MD)}=\sqrt{ \frac{2k_B T}{\D t} }\left(\tdisc{\Mfor}{n+1/2,*}-\tdisc{\Mfor}{n}\right) \tdisc{\Mfor}{n}^{-T/2} \tdisc{\gauss}{n}.
\end{equation}
This term might be impractical for large systems because it is based on solving a resistance problem to obtain $\Mfor^{-1/2}$. An alternative expression \rev{(used for nonlocal hydrodynamics)} that generates the same velocity in expectation is to use the random finite difference \cite{delong2014brownian}
\begin{gather}
\nonumber
\V{X}^\text{(RFD)} = \text{RotateAndIntegrate}\left(\tdisc{\V{X}}{n},\delta L  \tdisc{\M{K}}{n}^\dagger \gauss^\text{(RFD)} \right)\\ 
\label{eq:MRFD}
\tdisc{\V{U}}{n}^\text{(MD)} =\frac{k_B T}{\delta L} \left(\Mfor \left(\V{X}^\text{(RFD)}\right)-\tdisc{\Mfor}{n}\right) \gauss^\text{(RFD)},
\end{gather}
where $\gauss^\text{(RFD)}$ is a vector of i.i.d.\ standard Gaussian random variables (independent of $\tdisc{\gauss}{n}$), and $\delta =10^{-5}$ is a small number.
\item Compute an additional velocity to modify the backward Euler method for more accurate estimation of the covariance at a finite time step size \cite[Sec.~3.1]{maxian2023bending}
\begin{equation}
\tdisc{\V{U}^\text{(MBE)}}{n}=\sqrt{k_B T}\tdisc{\Mfor}{n}\M{L}^{1/2} \tdisc{\widetilde{\gauss}}{n}.
\end{equation}
\item Solve the saddle point system 
\begin{gather}
\label{eq:SPmidp}
\tdisc{\begin{pmatrix} -\Mfor & \M{K}+\D t \Mfor \M{L}\M{K}\\ \M{K}^T & \V 0 \end{pmatrix}}{n+1/2,*}
\tdisc{\begin{pmatrix} \V{\Lambda} \\ \V{\alpha} \end{pmatrix}}{n+1/2}=
\begin{pmatrix} -\tdisc{\Mfor}{n+1/2,*} \M{L}\tdisc{\V{X}}{n}+\tdisc{\V{U}^\text{(ext)}}{n} \\  \V{0}\end{pmatrix}\\ \nonumber
\text{where} \qquad \tdisc{\V{U}^\text{(ext)}}{n}=\tdisc{\V{U}^{(B)}}{n} +\tdisc{\V{U}^\text{(MBE)}}{n}+\tdisc{\V{U}}{n}^\text{(MD)},
\end{gather}
for $\tdisc{\V{\Lambda}}{n+1/2}$ and $\tdisc{\V{\alpha}}{n+1/2}$.
\item Update the fiber via
\begin{equation}
\tdisc{\V{X}}{n+1}=\text{RotateAndIntegrate}\left(\tdisc{\V{X}}{n},\tdisc{\V{\alpha}}{n+1/2}\Delta t \right)
\end{equation}
\end{enumerate}
The temporal integrator has a certain synergy with the split Stratonovich-Ito form \eqref{eq:kinetic}, in the sense that the Brownian velocity is evaluated at time $n$ (Ito interpretation), while the saddle point solve (which generates the action of $\Nhat$) is done at the midpoint (Stratanovich intepretation). 

\section{Numerical tests \label{sec:NT}}
In this section, we examine the mobility \eqref{eq:MNew} in both steady state and dynamic numerical tests, and demonstrate how changing between local and nonlocal mobilities can yield insights into the physics of deterministic and Brownian filament suspensions. We first consider the equilibrium bending fluctuations of filaments, where we compare the steady state end-to-end distance distribution to that obtained by MCMC using the energy \eqref{eq:GBDist}. We then use a simple example of filaments under gravity to demonstrate how the mobility \eqref{eq:MNew} can accurately resolve nonlocal hydrodynamics while treating the local mobility in an efficient way. After testing the accuracy of the mobility in a deterministic context, we perform simulations of two and four hundred Brownian filaments under gravity, all interacting hydrodynamically. Our simulations demonstrate anew the role of hydrodynamic interactions in constricting Brownian fluctuations \cite{manikantan2016effect, cunha2022settling}.

\subsection{Equilibrium fluctuations \label{sec:eqfl}}
We first test how the mobility \eqref{eq:MNew} changes the time step size required to successfully sample from the end-to-end distance distribution of free fibers at equilibrium. To do this, we consider a suspension of $F$ filaments with length $L=1$ in a periodic box of size $L_d=2$. We initialize the filaments as straight rods, then run the dynamics over a time $10\tau_\text{fund}$, where 
\begin{equation}
\label{eq:TauSm}
\tau_\text{fund} = 0.003 \left(\frac{4 \pi \mu L^4}{\kappa \ln{\left(\epsRS^{-1}\right)}}\right),
\end{equation}
is the longest bending timescale in the system \cite[Eq.~(55)]{maxian2023bending}. We remove the first $\tau_\text{fund}$ of the trajectory, then measure the end-to-end distance distribution of all filaments on the remaining $9\tau_\text{fund}$. We repeat this experiment twice for each parameter set to generate errror bars. The goal is to determine the time step size $\Delta t$, and the corresponding number of GMRES iterations at each time step, necessary to simulate the equilibrium distribution accurately.

\begin{figure}
\centering
\includegraphics[width=0.45\textwidth]{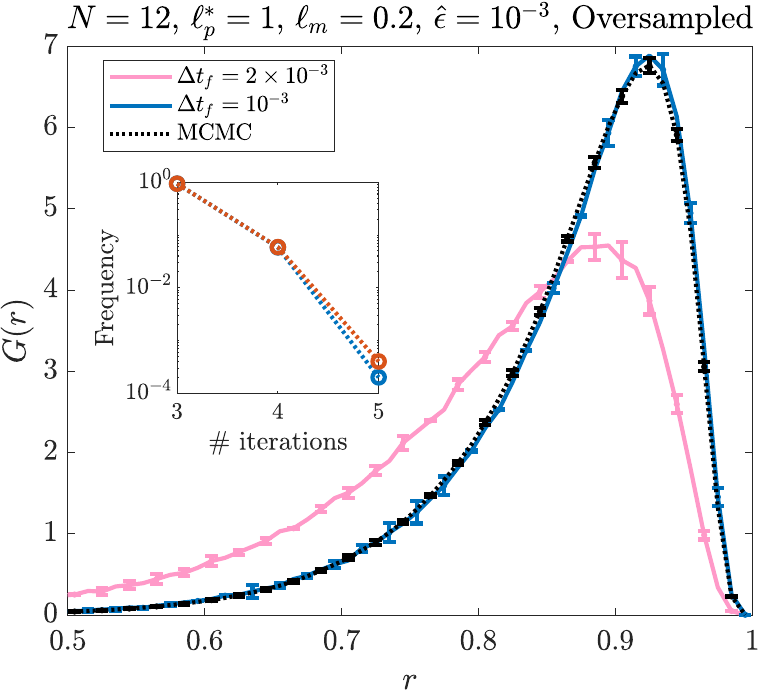}
\includegraphics[width=0.45\textwidth]{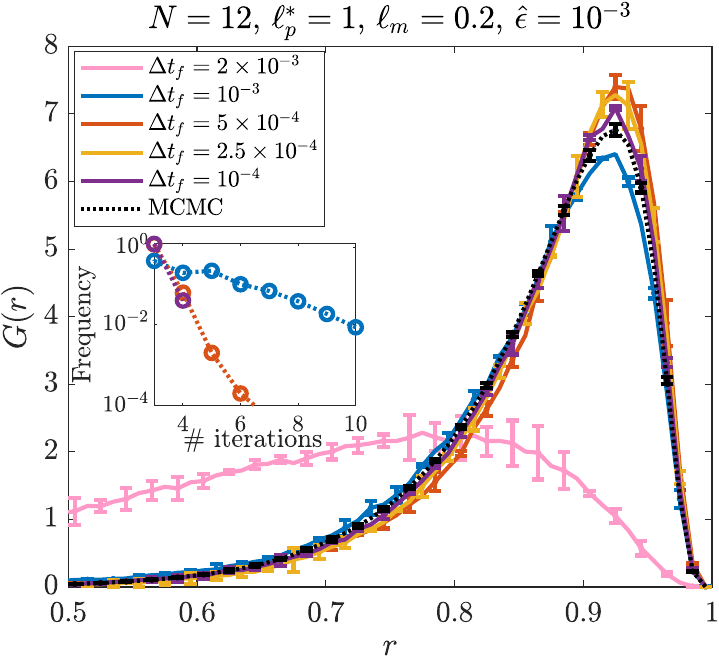}
\includegraphics[width=0.45\textwidth]{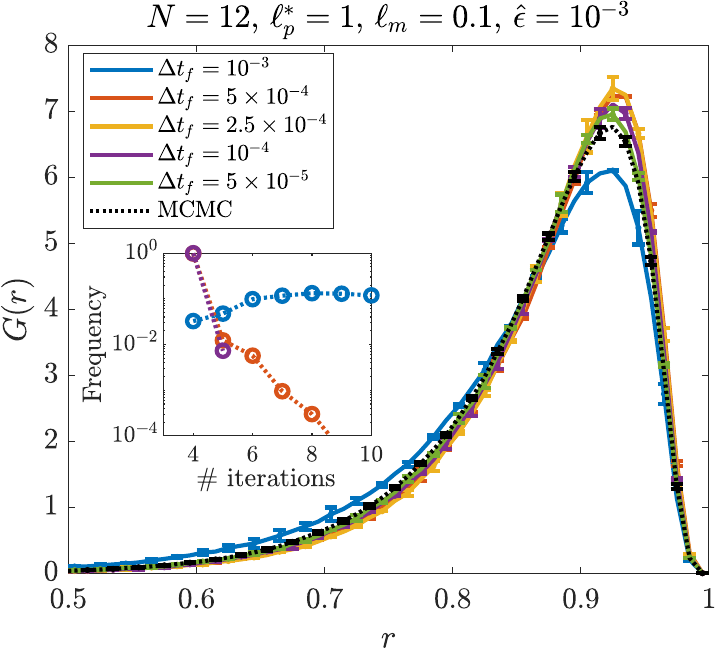}
\includegraphics[width=0.45\textwidth]{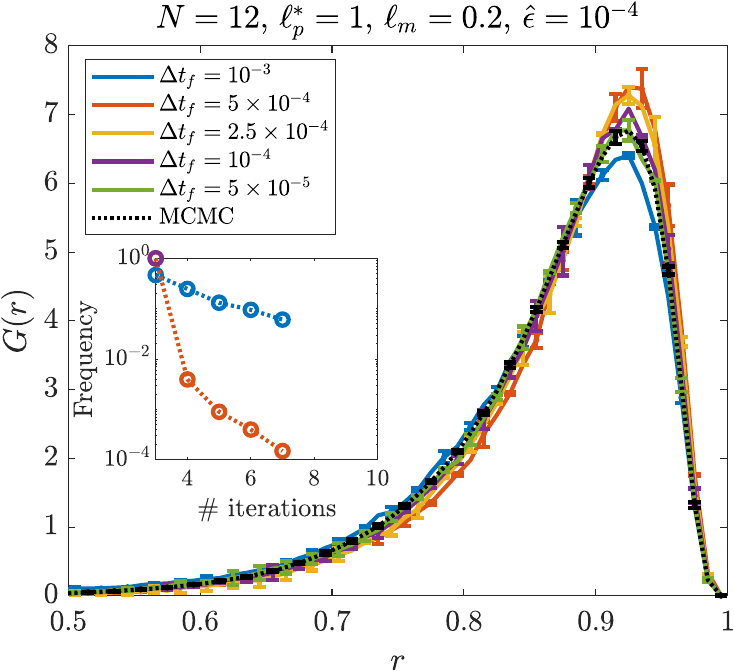}
\includegraphics[width=0.45\textwidth]{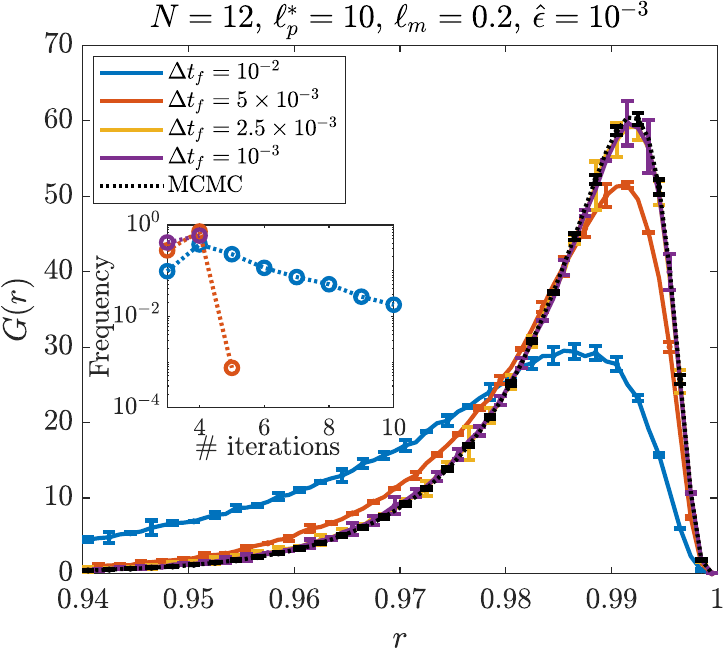}
\includegraphics[width=0.45\textwidth]{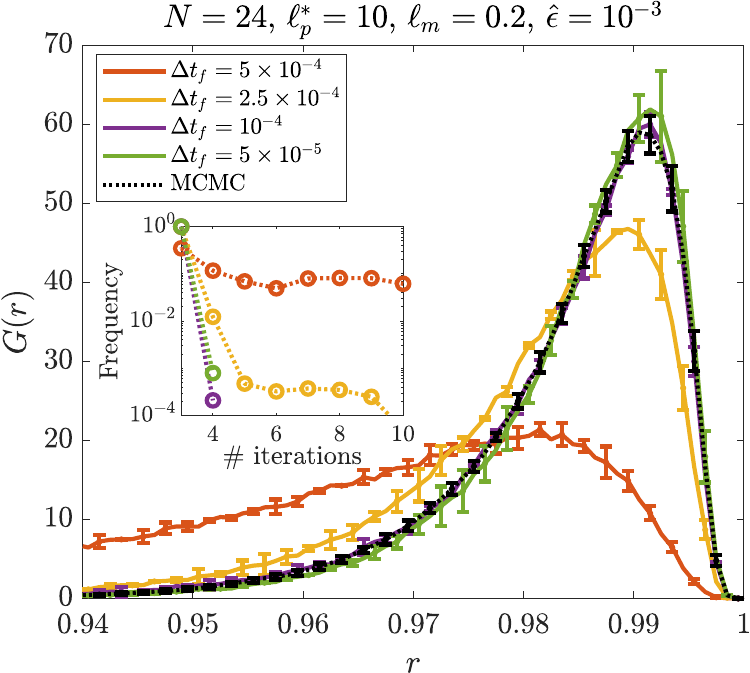}
\caption{\label{fig:ChkMCMC} End-to-end distance distribution with nonlocal hydrodynamics. In all cases we compare the results with Langevin dynamics with various time step sizes $\Delta t_f$ (in units of the slowest bending timescale) to MCMC, thereby determining the time step required to simulate the equilibrium distribution, and how it scales with $\ell_m$ (no change), $\hat \epsilon$ (no change), $\ell_p^*=\ell_p/L$ (relative time scale scales with $\ell_p^*$), and $N$ (time scale goes roughly as $N^{-4}$). The insets in each case show a probability histogram of the number of GMRES iterations required for convergence (tolerance $10^{-3}$).}
\end{figure}

Figure \ref{fig:ChkMCMC} shows how the required time step size and number of iterations vary with the system parameters. To obtain a point of reference, we begin at the top left with the \emph{oversampled} mobility \eqref{eq:Mref} for $\epsRS=10^{-3}$, which uses $N_u=500$ oversampling points globally (this is still an under-resolved configuration). For a persistence length $\ell_p^*=\ell_p/L=1$, mesh size $\ell_m=\sqrt{L_d^3/(FL)}=0.2$ (there are $F=200$ filaments), and $N=12$ tangent vectors, we see that $\Delta t_f=\Delta t/\tau_\text{fund}=10^{-3}$ is sufficient to match the distribution from Langevin dynamics (blue) with that from MCMC (dotted black). Interestingly, the number of GMRES iterations (inset) \emph{is not} a function of the fiber shapes; we almost always need at most four GMRES iterations for convergence, regardless of what the end-to-end distance statistics show. Our previous work \cite[Fig.~6]{maxian2023bending} also revealed a time step size of order $10^{-3}$ for local hydrodynamics only, so the nonlocal hydrodynamics with oversampling does not appear to change the required time step size.

We now consider what happens when we switch to the fat-corrected mobility \eqref{eq:MNew}. Here we use the same parameters as in the oversampled case, but set $\epsRSh=10^{-2}$, so that the nonlocal hydrodynamics can be resolved using only $N_u=100$ points. This makes each time step significantly cheaper. But is the time step size comparable to the oversampled case? The top right panel of Fig.\ \ref{fig:ChkMCMC} shows that the corrective part of the mobility \eqref{eq:MNew} introduces a small error when $\Delta t_f=10^{-3}$ (blue), which appears to be on the order of at most 10\% in the overall probability density function. This error takes a significantly smaller time step size to reduce; at $\Delta t_f=10^{-4}$ (purple) we finally match the distribution obtained from MCMC. Nevertheless, the intermediate distributions we obtain for $\Delta t_f=5 \times 10^{-4}$ (red) and $\Delta t_f=2.5 \times 10^{-4}$ (yellow) are acceptable for our purposes, meaning that the fat-corrected mobility requires a time step size two-fold lower than the oversampled one.

The remaining four panels in Fig.\ \ref{fig:ChkMCMC} show how the required time step size with the fat-corrected mobility changes with the parameters. We observe that halving the mesh size (the middle left plot has $\ell_m=0.1$, or $F=800$ filaments) does not have a strong effect on the required time step size ($\Delta t_f=5 \times 10^{-4}$ is still the first acceptable time step size). The same holds for decreasing $\epsRS$ by a factor of 10 (middle right plot), as $\Delta t_f=10^{-3}$ appears acceptable in that case. Increasing the persistence length by a factor of 10 (bottom left plot) gives an increase in the \emph{relative} allowable time step size; here $\Delta t = 2.5 \times 10^{-3}$ matches the (much tighter) equilibrium distribution. The most severe restrictions on the time step are when we increase the number of tangent vectors. Doubling to $N=24$ while keeping $\ell_p^*=10$ fixed yields an acceptable time step size of $\Delta t_f=10^{-4}$, roughly a factor of 16 lower than for $N=12$. Still, this extreme restriction (doubling the number of tangent vectors decreases the required time step as $N^{-4}$), is a consequence of the problem physics: if we want to resolve the end-to-end distance distribution correctly, we need to resolve every bending mode $k$, the timescale of which scales as $1/k^4$ \cite[Eq.~(24)]{poelert2012analytical}. Indeed, the restriction $N^{-4}$ on the required time step size is unchanged from when we considered local hydrodynamics only \cite[Fig.~6]{maxian2023bending}, and might be alleviated by developing a stochastic exponential integrator \cite{keavRPY}.

In addition to the required time step size, another important part of the mobility cost is the number of GMRES iterations required to solve \eqref{eq:SPmidp} (tolerance $10^{-3}$) for each time step size. Here we have a very convenient result: as shown in the insets of Fig.\ \ref{fig:ChkMCMC}, when we correctly resolve the equilibrium distribution, the number of GMRES iterations is three with a frequency of 99\%. Exceptions to this include larger persistence length, for which there are at most four iterations with frequency 99.9\%, and a smaller mesh size, in which the fewest number of iterations is four (this is unsurprising because nonlocal hydrodynamics makes a larger contribution). Likewise, when the fibers are more slender, so that nonlocal hydrodynamics makes a relatively smaller contribution, we see less iterations (compare the top right and middle right plots in Fig.\ \ref{fig:ChkMCMC}). Overall, the iteration counts give us a clear litmus test for the equilibrium distribution: when the number of iterations exceeds five with frequency larger than 1\%, we are not correctly simulating the equilibrium distribution. Furthermore, the iteration counts that we obtain when converged to the equilibrium distribution (three or four iterations) are the same as those obtained for oversampled quadrature, in which each iteration can be five to ten times more expensive. Thus, the only possible additional cost in the fat-corrected mobility is the small change in required time step size.

\subsection{Fibers in gravity \label{sec:grav}}
Because hydrodynamic interactions influence the behavior of both single filaments \cite{li2013sedimentation, cunha2022settling} and multiple interacting filaments \cite{llopis2007sedimentation, gustavsson2009gravity, saggiorato2015conformations, bukowicki2019sedimenting, makanga2023obstacle} under gravity, settling filaments provide an important context in which we can test our algorithms and various mobility approximations. We begin by considering a single filament under gravity, where we confirm the appearance of two different dominant bending modes, depending on the strength of the gravitational forces relative to the elastic forces \cite{saggiorato2015conformations, cunha2022settling}. In this context, we choose a parameter set with moderate gravitational forces, and examine how different mobilities affect the deterministic behavior of two sedimenting filaments, finding that our fat-corrected mobility gives errors on the same order as the non-SPD mobility \eqref{eq:MFirst}. We then expand our study to Brownian filaments, showing that entrainment flows can effectively confine translational and rotational diffusion which would otherwise prevent the fibers from approaching each other when they fluctuate. We conclude by demonstrating our algorithm's capability to simulate a larger suspension of hundreds of filaments, where the fibers segregate into clumps \cite{gustavsson2009gravity}.

Throughout this section, we fix the fiber aspect ratio $\epsc=10^{-3}$, which is thirty to one hundred times smaller than previous studies which were confined to bead-link representations of the filaments \cite{saggiorato2015conformations,cunha2022settling, makanga2023obstacle}. When the aspect ratio is fixed, the dynamics of deterministic chains are governed by the elasto-gravitational number and elastic timescale \cite{bukowicki2019sedimenting}
\begin{equation}
\label{eq:EGNum}
\beta = \frac{g L^3}{\kappa} \qquad \bar t = \frac{L^4 \mu}{\kappa},
\end{equation}
where $g$ is the gravitational force \emph{density} on the filaments ($gL$ is the total force on each filament). \rev{The simulations here do \emph{not} use a bottom boundary; the single filament is isolated in free space (falling in a domain of infinite extent), while simulations with multiple filaments use a triply periodic domain (see the conclusion for further discussion on a bottom boundary).}

\subsubsection{One deterministic filament}

\begin{figure}
\centering
\includegraphics[width=\textwidth]{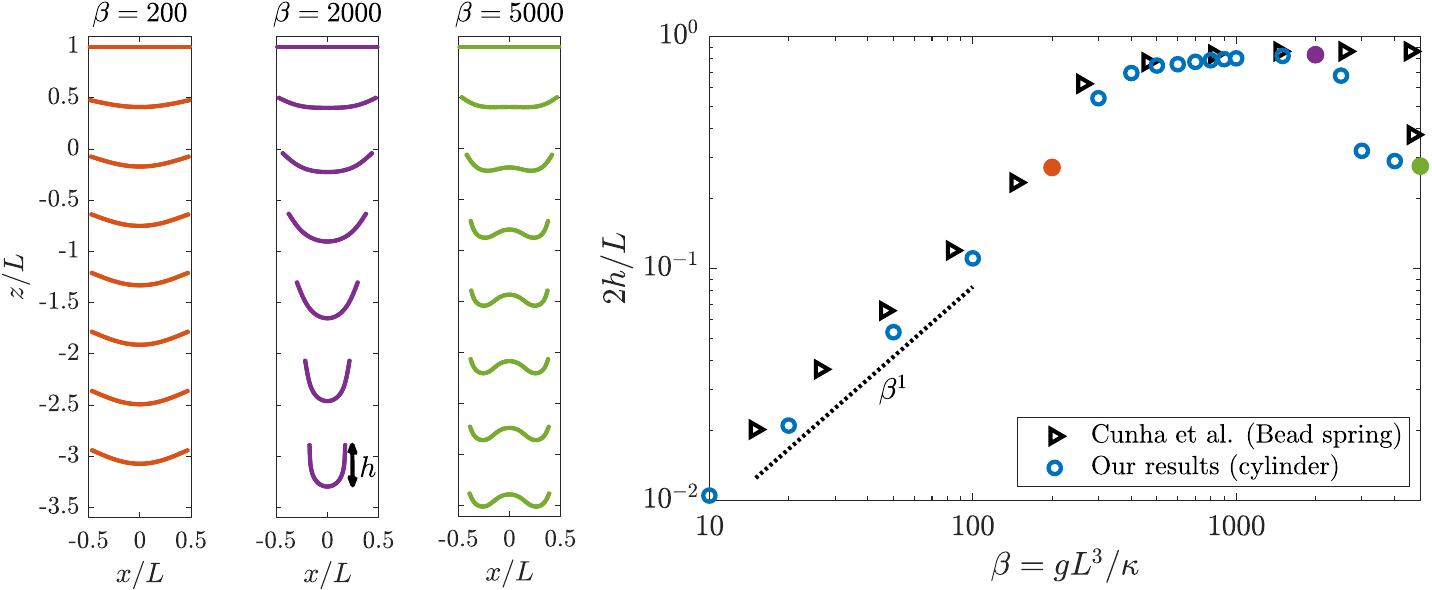}
\caption{\label{fig:Single} Deterministic dynamics of a single falling filament. The left panels show a series of snapshots of the sedimenting filament for a variety of elasto-gravitational numbers $\beta=gL^3/\kappa$. The right panel shows a summary of the vertical extension of the filament $h$ as a function of $\beta$, \rev{compared to a previous bead-spring model with slightly different model geometry \cite[Fig.~4]{cunha2022settling}}. Colored dots on the right plot correspond to the colored fiber snapshots shown in the left panels.  }
\end{figure}

We begin with a single falling filament in isolation, using the special quadrature scheme \eqref{eq:Msqs} for the mobility, and the deterministic time-stepping algorithm \eqref{eq:SPDBE} with $\Mfor=\ML$ (the discretization parameters are $N=16$ and $\Delta t = 2 \times 10^{-4} \left(\bar{t}/g\right)$). As shown in Fig.\ \ref{fig:Single}, we begin with the fiber horizontally-oriented, then run forward in time until it assumes a steady state shape. In agreement with previous observations \cite[Fig.~4]{cunha2022settling}, there are three regimes of behavior: for small $\beta$, the fibers deflect only slightly, with a vertical extension $h$ proportional to $\beta$. For intermediate $\beta$, the fibers assume a maximum extension which is roughly $0.8(L/2)$, or 80\% of the maximum radius. For larger $\beta$, a new bending mode emerges where the fiber assumes a ``W'' instead of ``U'' shape. \rev{We observe small quantitative discrepancies from previous results \cite[Fig.~4]{cunha2022settling} because of differences in fiber geometry; our filaments are continuous rods, whereas the simulations in \cite{cunha2022settling} considered a discrete set of beads of uniform radius $a$ connected by springs of radius $3a$. Nevertheless, good quantitative agreement still results because the shape changes in both cases are driven by nonlocal hydrodynamic interactions, as opposed to non-uniformity of filament thickness, which dictates the corresponding shapes of spheroidal filaments \cite{li2013sedimentation}}.

\subsubsection{Two deterministic filaments}
\begin{figure}
\centering
\includegraphics[width=\textwidth]{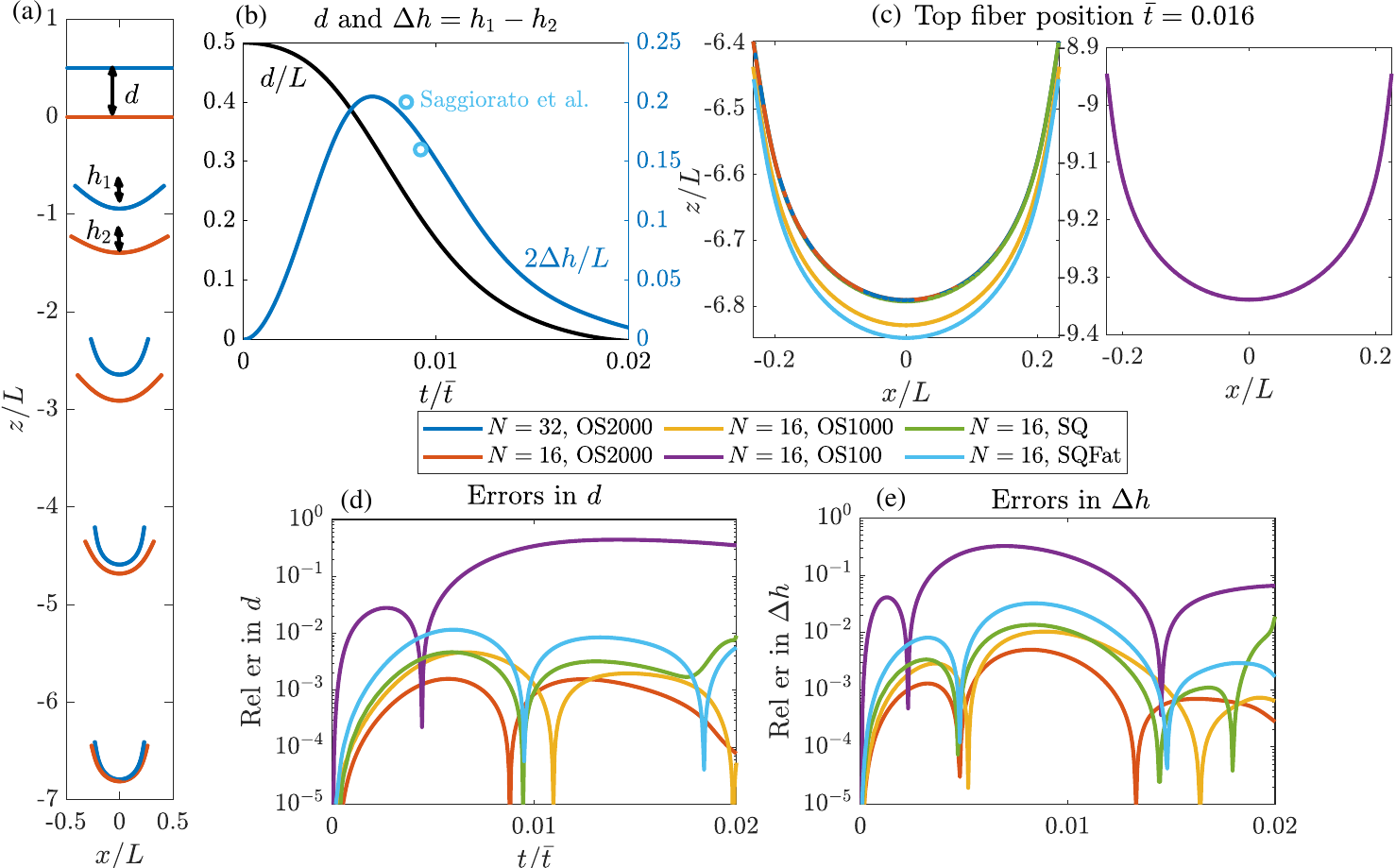}
\caption{\label{fig:TwoFallDet} Using sedimentation of two deterministic filaments to study the accuracy of different nonlocal mobilities. (a) The set-up and dynamics. We consider two filaments in the $xz$ plane, with a series of snapshots shown every $0.004\bar{t}$. (b) The key statistics are the distance between the two fiber midpoints, $d=\ind{\V{X}}{2}(L/2) -\ind{\V{X}}{1}(L/2)$ and the difference in the two radii of curvature $\Delta h = h_1-h_2$. \rev{For comparison, we pull values of $(d,\Delta h)$ from the top two panels of \cite[Fig.~5]{saggiorato2015conformations}, align $t$ to match the value of $d$, and plot the corresponding value of $\Delta h$ with a light circle.} (c) The position $\ind{\V{X}}{1}(\bar t=0.016)$, (d) the errors in $d$ over time, and (e) the errors in $\Delta h$ over time, for the different mobility options. In (c,d,e), we show $N=32$ with 2000 oversampled points in blue (this is the reference curve), $N=16$ with 2000 oversampled points in red, $N=16$ with 1000 oversampled points in yellow, $N=16$ with 100 oversampled points in purple, $N=16$ with special quadrature \eqref{eq:MFirst} in green, and the fat-corrected mobility \eqref{eq:MNew} with $\epsRSh=10^{-2}$ in cyan. }
\end{figure}

Having established the regimes of behavior for a single filament, we now consider the interaction of two deterministic filaments. Here we focus on the \emph{numerical} aspects of the problem, in particular how the different mobility choices perform as the fibers approach each other. As shown in Fig.\ \ref{fig:TwoFallDet}(a), we begin with two fibers parallel in the plane, spaced $d/L=0.5$ apart. As these filaments sediment, the top fiber becomes more curved than the bottom one, and falls faster due to the flows created by the second fiber. As a result, the fibers approach each other, eventually assuming the same shape and falling together  \cite{llopis2007sedimentation,saggiorato2015conformations}. 

The importance of hydrodynamic interactions makes this problem an ideal set-up to examine the accuracy of different nonlocal mobilities, in particular our three candidates \eqref{eq:Mref}, \eqref{eq:MFirst} (which is only possible in deterministic simulations), and \eqref{eq:MNew}. To do this, we fix $\beta=500$, then choose a small time step size $\Delta t = (2.5\times 10^{-4})(\bar t/g)$ to eliminate temporal errors (the final time is $t_f=0.1\bar t/g$). We use the deterministic time-stepping algorithm with time-lagged nonlocal forces \eqref{eq:BEnL}, and use a ``periodic'' domain with $L_d/L=10$, having verified this is sufficient to eliminate periodic artifacts. 

We establish a reference trajectory of the filaments by using the oversampled mobility \eqref{eq:Mref} with $N=32$ collocation points and $N_u=2000$ upsampling points, and collecting the distance between the fiber midpoints $d$, and the difference in the vertical extension of the two fibers ($\Delta h=h_1-h_2$), \rev{which we compare (Fig.\ \ref{fig:TwoFallDet}(b)) to two reference points $(d,\Delta h)$ obtained from \cite[Fig.~5]{saggiorato2015conformations}}. We then systematically change the mobility, first repeating the oversampled mobility \eqref{eq:Mref}, but with $N=16$ collocation points and $N_u=2000$ upsampling points. The smaller number of collocation points has no effect on the dynamics to three digits of accuracy, and the trajectories of the fibers are identical. We then change the number of upsampling points to $N_u=1000$, where we observe noticeably larger errors of order 0.01 in the fiber positions, extensions, and separation. Decreasing the number of oversampling points further to $N_u=100$ yields disastrous consequences, as the self-RPY integrals are not resolved. In this case the fibers fall much faster than they otherwise should and remain farther apart from each other for a longer time, destroying any semblance of accuracy in the trajectories (errors are order 1). 

We can rescue the original dynamics by still using $N_u=100$ oversampling points for the \emph{nonlocal} parts of the mobility, but using special quadrature for the local parts, i.e., using the mobility \eqref{eq:MFirst}. This mobility gives roughly 2--3 digits of accuracy in the fiber position and other statistics, but it is untenable for Brownian fluctuations because it is not SPD. For this reason, we consider the fat-corrected mobility \eqref{eq:MNew} with $\epsRSh=10^{-2}$, using $N_u=100$ oversampling points. This mobility (cyan) gives differences from the reference that are slightly higher than, but still on the same order as, the special quadrature mobility (green) and oversampled mobility with $N_u=1000$ oversampling points (yellow). This demonstrates that the fat-corrected mobility is a viable option for simulations: it reproduces the accuracy of other reasonable mobility options, while at the same time being the sum of two SPD pieces.

\subsubsection{Two Brownian filaments}
Since the fat-corrected mobility allows us to inexpensively include Brownian fluctuations with nonlocal hydrodynamics, we now consider the role of Brownian hydrodynamics in the suspension of two falling filaments. We set $k_BT =\kappa/L $, so that the persistence length $\ell_p=L$, then repeat the simulations of Fig.\ \ref{fig:TwoFallDet} with four different hydrodynamic models: deterministic local, deterministic nonlocal, Brownian local, and Brownian nonlocal. We use $N=16$ and the fat-corrected nonlocal mobility \eqref{eq:MNew} for all simulations, with time step sizes $\Delta t=\left(2 \times 10^{-6}\right) \bar t$ for local hydrodynamics and $\Delta t=10^{-6} \bar t$ for nonlocal hydrodynamics \rev{(these are required to reproduce the correct end-to-end distance distribution when $g=0$)}. By comparing the hydrodynamic models, we determine if nonlocal entrainment flows play a role even in the presence of Brownian motion.

\begin{figure}
\centering
\includegraphics[width=\textwidth]{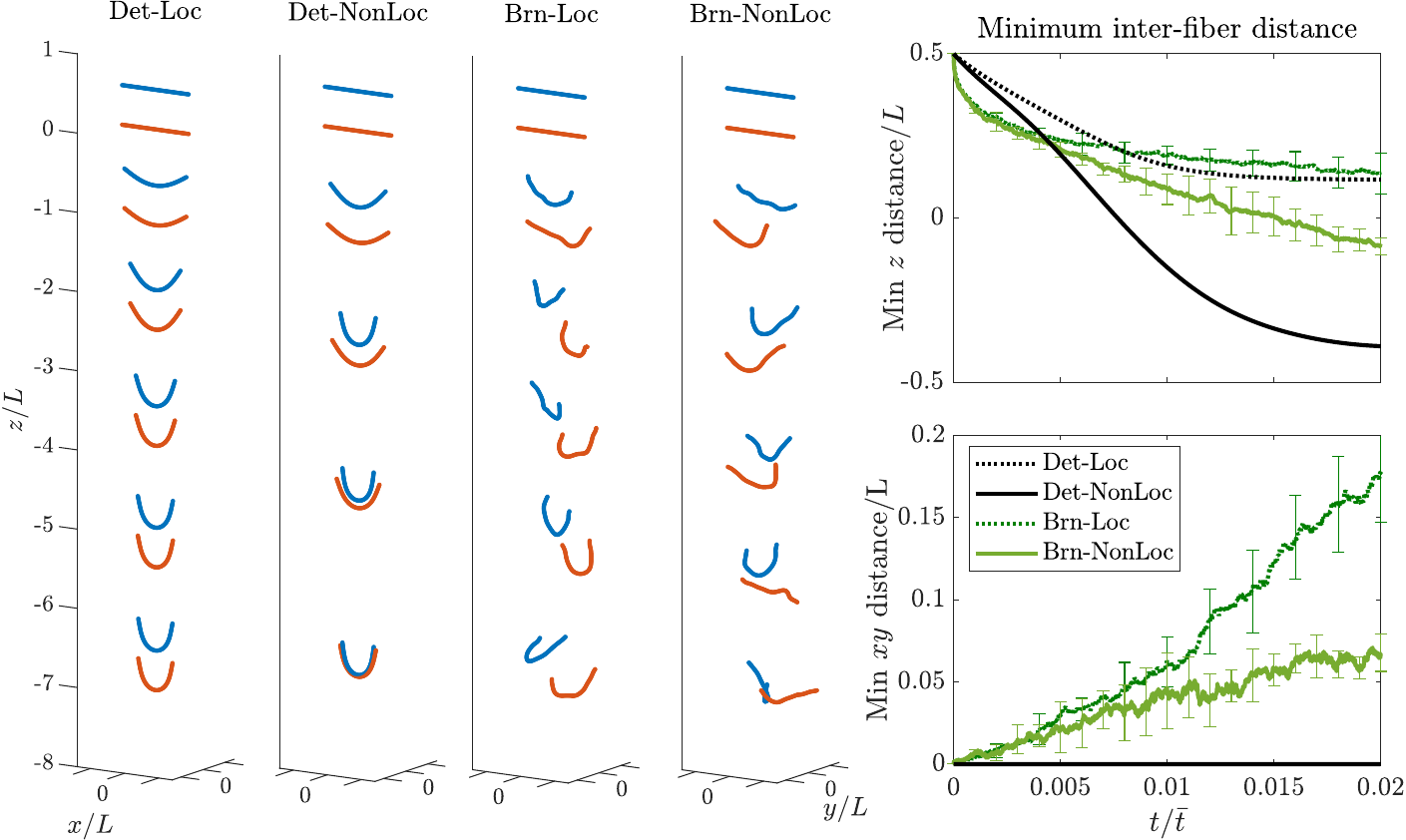}
\caption{\label{fig:TwoFallBrn} Hydrodynamics traps fluctuations in sedimenting Brownian filaments. We repeat the simulations of Fig.\ \ref{fig:TwoFallDet} with four different hydrodynamic models: deterministic local (left), deterministic nonlocal (second from left), Brownian local (third from left), and Brownian nonlocal (fourth from left). At right, we show the mean inter-fiber distance $\ind{\V{X}}{1}-\ind{\V{X}}{2}$, where $\ind{\V{X}}{1}$ denotes the top (blue filament), in the $z$ direction (parallel to gravity) and $xy$ direction (perpendicular to gravity). The Brownian examples shown at left are representative of the mean. With nonlocal hydrodynamics, fibers approach each other faster in the $z$ direction, and are prevented from diffusing away from each other in the $xy$ direction.}
\end{figure}

The results in Fig.\ \ref{fig:TwoFallBrn} show that, indeed, nonlocal flows attract the two filaments to each other, impeding translational diffusion in the $xy$ plane. The representative examples at left first show \emph{deterministic} motion; in this case, fibers simulated with local hydrodynamics fall until they reach a steady state shape, with the top fiber being shifted by $z/L=0.5$ units relative to the bottom. With nonlocal hydrodynamics, the fibers fall faster and are attracted to each other, as we have already seen (this is the same simulation as Fig.\ \ref{fig:TwoFallDet}). With Brownian motion, the fibers simulated with local hydrodynamics appear to escape each other in the $xy$ plane, and never come into close contract in the $z$ plane. By contrast, Brownian fibers simulated with nonlocal hydrodynamics stay in closer contact, appearing to be entrained in each other's flow fields. The minimum inter-fiber distance, projected onto the $z$ axis ($\text{min}_s\left(\ind{\V{X}}{1}_3(s)-\ind{\V{X}}{2}_3(s)\right)$; can be negative) and $xy$ plane ($\text{min}_s\norm{\ind{\V{X}}{1}_{1:2}(s)-\ind{\V{X}}{2}_{1:2}(s)}$; positive), shows that the differences between local and nonlocal hydrodynamics are quite significant, demonstrating that nonlocal fluid flows can confine translational and rotational diffusion of the filaments. 

\subsubsection{Brownian suspensions with nonlocal hydrodynamics and steric interactions}
\begin{figure}
\centering
\includegraphics[width=\textwidth]{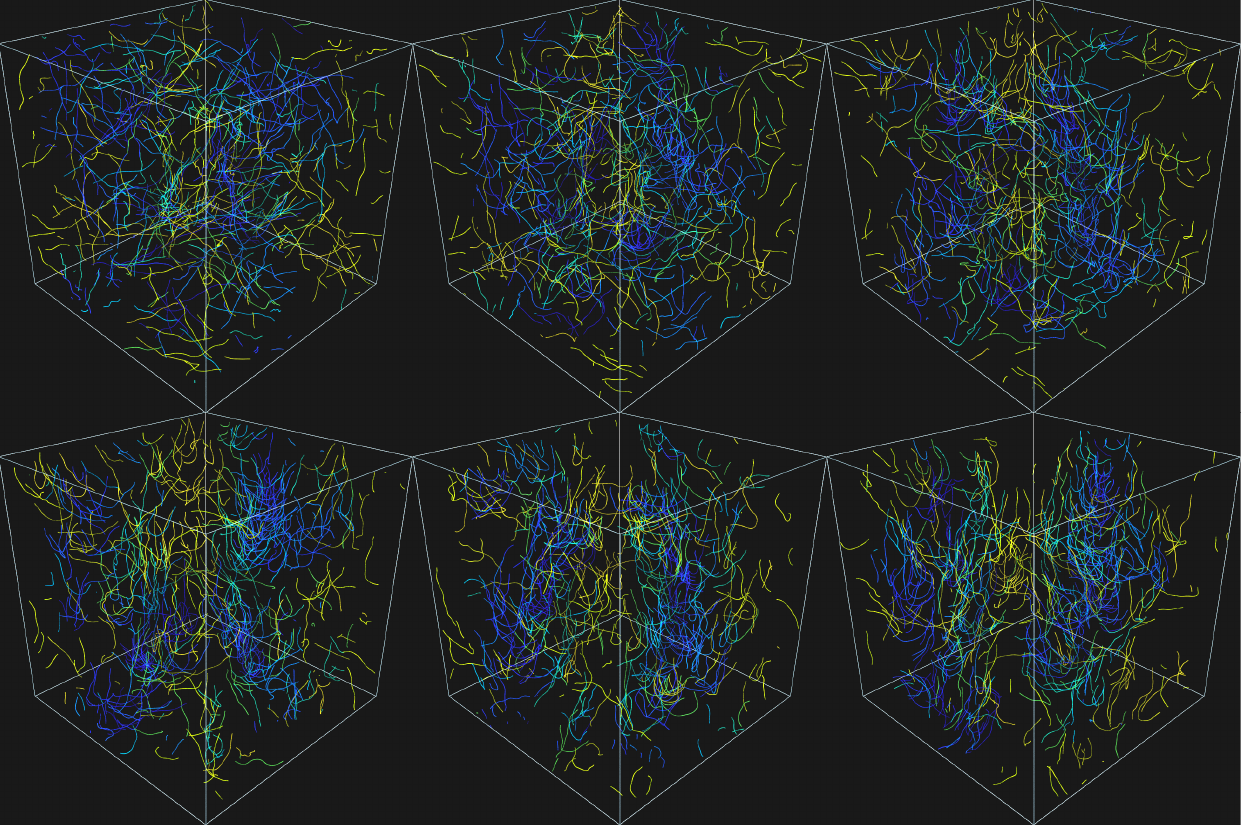}
\caption{\label{fig:BrownianSed}Simulation of 400 sedimenting Brownian fibers with $L=1$ on a domain of size $4 \times 4 \times 4$. Snapshots are shown at equally-spaced intervals from $t=0$ to $t/\bar{t}=0.02$. With gravity oriented downward, the fibers are color-coded by the vertical distance traveled over the simulation time, with blue having the most negative vertical displacement and yellow the most positive.}
\end{figure}

Finally, we demonstrate the capability of the simulation platform by scaling up to $F=400$ semiflexible filaments (still with length $L=1$) in a domain of size $L_d=4$. All parameters are unchanged from previous sections, except that we switch $\epsc=4 \times 10^{-3}$ so fibers have larger volume and steric effects are more important. We simulate sedimentation with Brownian fluctuations and steric interactions (using steric force constant $\Pot_0=1000k_B T$ in \eqref{eq:PotSter}), showing some snapshots in Fig.\ \ref{fig:BrownianSed}. At initial times, the fibers can be seen to be fluctuating, and gravity can hardly be detected. At later times, the nonlocal fluid flows build up, and the fibers change shape and separate into two visible groups, with one falling downwards and the other moving upwards. These qualitative observations agree with previous studies on sedimentation \cite{gustavsson2009gravity, manikantan2016effect}. As shown there \cite{gustavsson2009gravity}, a more detailed estimation of the sedimentation velocity, as well as how the ``flocks'' of fibers change over time, requires ensemble averaging over long times (even for deterministic simulations). While these simulations are certainly possible in our framework, we leave them for future work, which could focus also on the competition between Brownian fluctuations and hydrodynamic interactions in setting the fiber orientation \cite{manikantan2016effect, cunha2022settling}.

\section{Cross-linked fiber networks \label{sec:Bundling}}
The original motivation for this work was to examine the dynamics of the cytoskeleton, in particular to mimic an \emph{in vitro} reconstitution of actin fibers with cross-linking proteins \cite{schmoller2009structural, lieleg2009structural, falzone2012assembly, foffano2016dynamics}. In the experimental system \cite{falzone2012assembly}, actin and various concentrations of the cross linker $\alpha$-actinin are added to a solution, and steady state architectures taken some time later reveal various degrees of fiber agglomeration and bundling. By reconstituting this system \emph{in silico}, we were able to study the role of various parameters and physical phenomena, such as cross-linker kinetics and filament bending stiffness \cite{maxian2021simulations, maxian2021bundling}, in shaping the dynamics of the bundling process. Here we complete these studies by examining the role of steric and hydrodynamic interactions in the bundling process, noting in particular the role of sterics vs.\ cross linking in the kinetic arrest of bundling fibers \cite{falzone2012assembly, foffano2016dynamics}.

\subsection{Model set up \label{sec:CLs}}
In our model of a cross-linked fiber network, the cross linkers are virtual entities which exert time-dependent external forces $\FCL(t)$ on the filaments. The details of the algorithm for transient cross linking are discussed in \cite{maxian2021simulations, maxian2021bundling}; here we give a brief summary in the interest of completeness. At each time, the ``cross-linked network'' is defined by a list of pairs of points (on distinct filaments) which are connected by a cross linker (CL). To define the list of potential points, we sample each filament $\ind{\XPoly}{i}(s)$ at $N_\text{cl}$ uniformly-spaced points $s^*_i$ with spacing $\Delta s_\text{cl} = L/(N_\text{cl}-1)$. At a given time, a CL can have one end bound to a binding site (``singly-bound''), or both ends bound to binding sites on distinct filaments (``doubly-bound''). Each singly- or doubly-bound CL has a rate of unbinding: the rate at which a singly-bound CL unbinds is $\koff$, while a doubly bound CL can unbind from either binding site with rate $\koffb$ (making the total rate $2\koffb$). 

We model the process of CL diffusion and binding by an effective attachment rate for one CL end, given at each site by $\kon \Delta s_\text{cl}$. If a binding site on a nearby filament $\ind{\XPoly}{j}(s)$ is a distance $\ell=\norm{\ind{\XPoly}{i}\left(s_p^*\right)-\ind{\XPoly}{j}(s_q^*)}$ away from the docking site of the first end, the energy associated with binding to that site is given by 
\begin{equation}
\label{eq:ECL}
\mathcal{E}^\text{(CL)} = \frac{K_c}{2}\left(\norm{\ind{\XPoly}{i}\left(s_p^*\right)-\ind{\XPoly}{j}\left(s_q^*\right)}-\ell_c \right)^2 =\frac{K_c}{2}\left(\ell-\ell_c\right)^2,
\end{equation}
and the rate of binding of the second end to the nearby binding site is given by
\begin{equation}
\label{eq:koneqn}
\konb\left(\ell \right)= \konb^0 \exp{\left(-\frac{\mathcal{E}^\text{(CL)}}{k_B T}\right)},
\end{equation}
where $\ell_c$ is the rest length of the CL, and $K_c$ is its stiffness. In \cite{maxian2021bundling}, we derived the relationship\ \eqref{eq:koneqn} by assuming that the fluctuations in CL length are in equilibrium with respect to the energy \eqref{eq:ECL}, thus ensuring that the links are in detailed balance and therefore passive. Because the binding rate \eqref{eq:koneqn} is compactly-supported (to finite precision), we perform a neighbor search to obtain a list of all uniform point pairs within two standard deviations of the Gaussian \eqref{eq:koneqn}; that is, all uniform point pairs separated by a distance $\ell_c+2\sqrt{k_B T/K_c}$. This defines a list of the possible pairs of points to which new CLs can bind, and completes a list of four possible reactions which we simulate using a standard stochastic simulation / Gillespie algorithm \cite{gillespie2007stochastic,  anderson2007modified,maxian2021simulations}.

We use a first-order splitting algorithm to update the cross linkers and filaments. At each time step, we first take a step of size $\Delta t$ in the Markov chain. This produces a list of fibers and uniform points which are connected by cross linkers. To compute forces, we denote the discrete uniform points as $\Xuni = \M{R}^{(u)} \V{X}$, so that the $p$th uniform point on fiber $i$ can be written as $\ind{\Xuni}{i}_{\Bind{p}} = \M{R}^{(u)}_{\Bind{p,:}} \V{X}$, and the energy \eqref{eq:ECL} can be written in terms of $\V{\ell}=\ind{\Xuni}{i}_{\Bind{p}}-\ind{\Xuni}{j}_{\Bind{q}}$. \rev{In a similar way to the steric force calculation \eqref{eq:sumForceUp}, the} CL force on each (Chebyshev) fiber point $\ind{\V X}{i}_{\Bind{a}}$ can be obtained by differentiating the energy \eqref{eq:ECL}, 
\begin{align}
\label{eq:ForceCL}
\ind{\V{F}}{i}_{\Bind{a}}&=-\frac{\partial \Pot^\text{(CL)}}{\partial \ell}\widehat{\V{\ell}}\M{R}^{(u)}_{\Bind{p,a} } \qquad \qquad \ind{\V{F}}{j}_{\Bind{b}}= \frac{\partial \Pot^\text{(CL)}}{\partial \ell}\widehat{\V{\ell}}\M{R}^{(u)}_{\Bind{q,b} }.
\end{align}
Adding up the forces \eqref{eq:ForceCL} for all pairs of points gives the external force $\FCL$ which is added as an external force to the right-hand side of the saddle-point system \eqref{eq:SPmidp}. We note that the CL forces can be nonsmooth; in the special case when the two uniform points are Chebyshev points, we get $\M{R}^{(u)}_{\Bind{p,a}} =\delta_{ap}$ and $\M{R}^{(u)}_{\Bind{q,b}}=\delta_{qb}$, and the resulting force is a delta function between the two points on the Chebyshev grid.

\subsection{Basic dynamics}
\begin{figure}
\centering
\includegraphics[width=\textwidth]{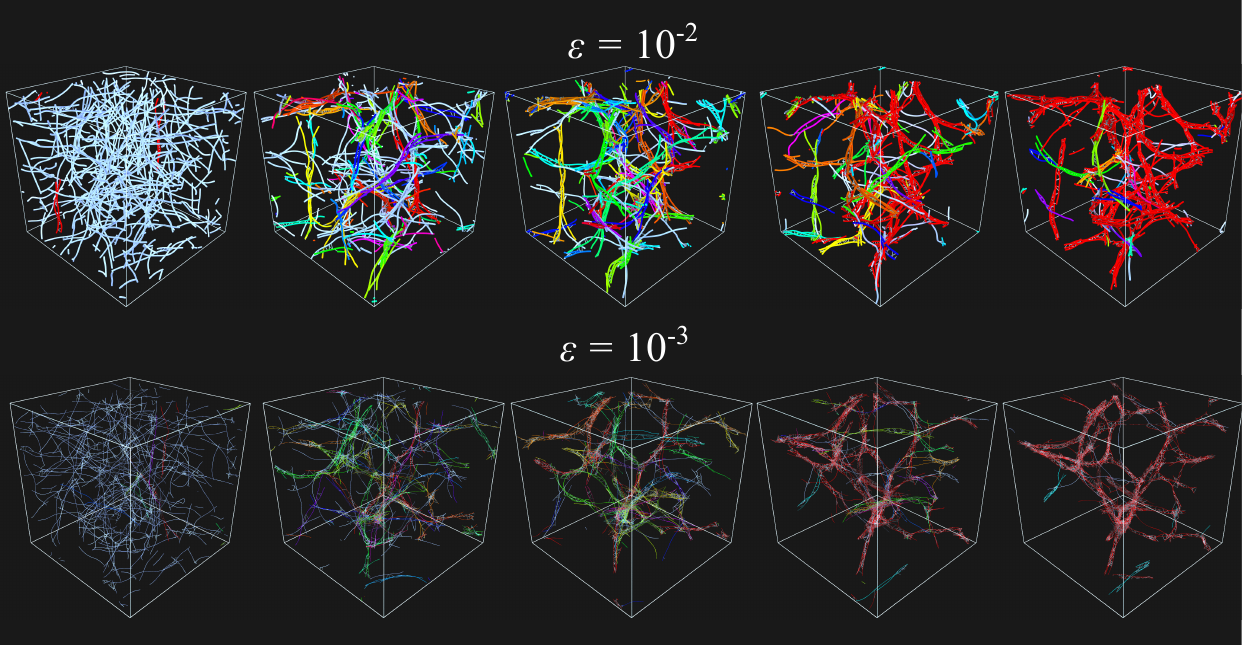}

\vspace{1 cm}

\includegraphics[width=0.9\textwidth]{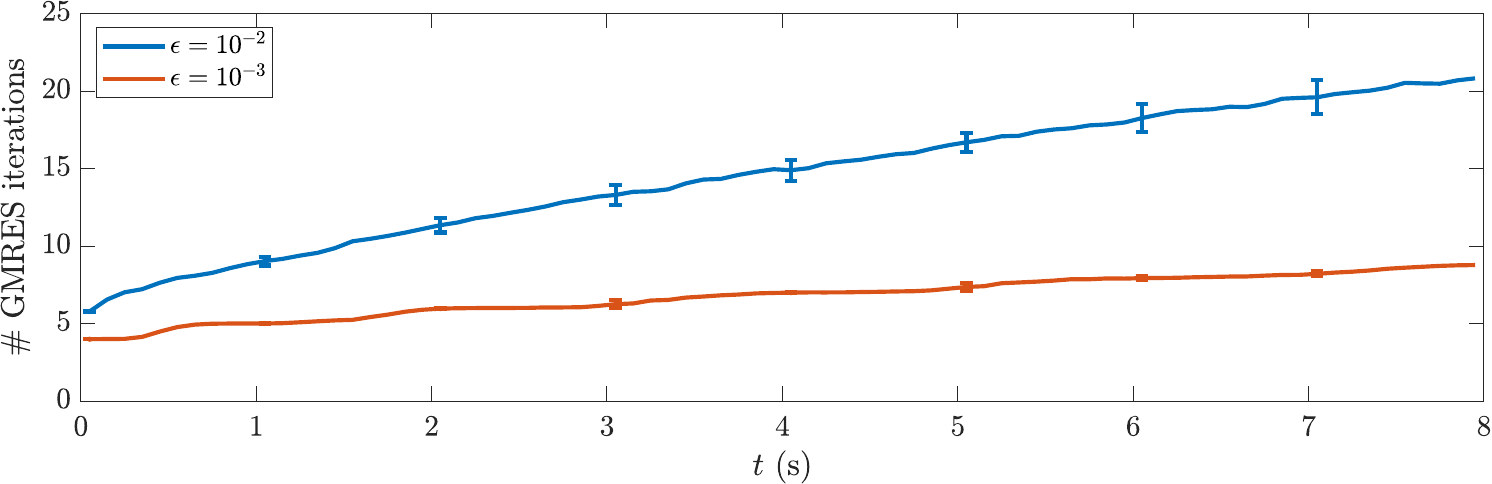}
\caption{\label{fig:BundSeq} Sequence of bundling with two different aspect ratios. Top figure: Snapshots are shown every 2 seconds from $t=0$ to $t=8$ with $\epsc=10^{-2}$ (top) and $\epsc=10^{-3}$ (bottom) in simulations with local hydrodynamics and no steric interactions. Fibers are color-coded by bundle; those not in bundles are shown in white. Thinner white filaments are cross linkers. Bottom: average number of GMRES iterations required to solve the system over time for $\epsc=10^{-2}$ (blue) and $\epsc=10^{-3}$ (red). Data points are collected by averaging over all timesteps in 0.1 s intervals. }
\end{figure}

Since the effect of hydrodynamics and sterics will presumably depend upon the fiber aspect ratio, we consider bundling fiber networks with both $\epsc=10^{-2}$ and $\epsc=10^{-3}$. To get a sense of the dynamics, we first simulate 8 seconds of bundling with local hydrodynamics, neglecting steric interactions. The complete list of parameters for these simulations can be found in \cite[Table~1]{maxian2021bundling}; here we will only list deviations from these parameters. For example, we maintain $L=1$ $\mu$m, and vary the fiber radius $\rc$ in accordance with the aspect ratio of interest. As we previously showed \cite[Sec.~6]{maxian2023bending} that bending fluctuations are only of importance when $\ell_p/L \sim \mathcal{O}(1)$, in these simulations we set $\ell_p=L=2$ $\mu$m, making $\kappa=\ell_p k_B T = 0.0082$ pN$\cdot$ $\mu$m$^2$. We use $N=12$ tangent vectors to represent the filaments, which have initial mesh size $\ell_m=\sqrt{L_d^3/(FL)}=0.2$ $\mu$m. This corresponds to $F=675$ filaments in a domain of extent $L_d=3$ $\mu$m on each side, and $F=200$ filaments in a domain of extent $L_d=2$. For this paper, all images are shown on the smaller domain, while statistics are collected on the larger domain. Unless otherwise noted, we use $n=5$ simulations to generate error bars (two standard errors in the mean). The steric force constant in \eqref{eq:PotSter} is $\Pot_0=1000k_B T$ pN. 

We previously showed \cite{maxian2021bundling, maxian2023bending}, that the flexibility of the CLs, encoded via the energy \eqref{eq:ECL}, drives bundling via a thermal zippering mechanism. When a CL binds two filaments in a stretched configuration, it subsequently relaxes, pulling filaments together through the forcing \eqref{eq:ForceCL}. The approach of the two filaments allows for further CL binding, which creates a zippering mechanism that tends to align filaments in parallel bundles. Subsequent stages of bundling bring multiple bundles together, forming huge agglomerates which span the entire simulation domain. 

Figure \ref{fig:BundSeq} shows this process playing out in our simulations with two different aspect ratios.  Snapshots of the agglomeration/bundling process are shown at timepoints $t=0,2, \dots 8$ s, and fibers are colored by ``bundle.'' Here bundles are defined by mapping the filaments to a graph, where a connection exists between two filaments if they are cross linked in two locations $L/4$ apart. Each ``bundle'' is a connected region in this graph. Figure \ref{fig:BundSeq} shows that the agglomeration process works in two stages, whereby filaments first organize themselves into many bundles of a few filaments, which then coalesce into larger connected structures. By the end of the simulation ($t=8$), one large (red) bundle spans the entire simulation box, and the filaments are decorated with (small white) CLs. 

From a numerical perspective, the bundling process provides an interesting test for our nonlocal hydrodynamic solvers. How robust are they to the network architecture? To explore this, we repeat the simulations of Fig.\ \ref{fig:BundSeq} with nonlocal hydrodynamics. For the larger aspect ratio $\epsc=10^{-2}$, we use the oversampled mobility \eqref{eq:Mref} with $N_u=100$, while the smaller aspect ratio $\epsc=10^{-3}$ requires the fat-corrected mobility \eqref{eq:MNew} with $\epsRSh=10^{-2}$ (and 100 oversampled points). The bottom panel of Fig.\ \ref{fig:BundSeq} shows that, indeed, the number of GMRES iterations required to solve the system increases as networks get more bundled. This is no surprise; when filaments come closer together, hydrodynamic interactions become more important, and the local dynamics becomes less dominant. The increase in iterations is therefore particularly egregious when $\epsc=10^{-2}$, where we see about 4 times as many iterations required at the end of the simulation than at the beginning. For the more slender $\epsc=10^{-3}$, the local dynamics still dominates even in a highly-bundled state, and so the average number of iterations is only about nine at the end of the simulation, an increase from about five near $t=0$.

\subsection{Role of nonlocal hydrodynamics and sterics}
\begin{figure}
\centering
\includegraphics[width=\textwidth]{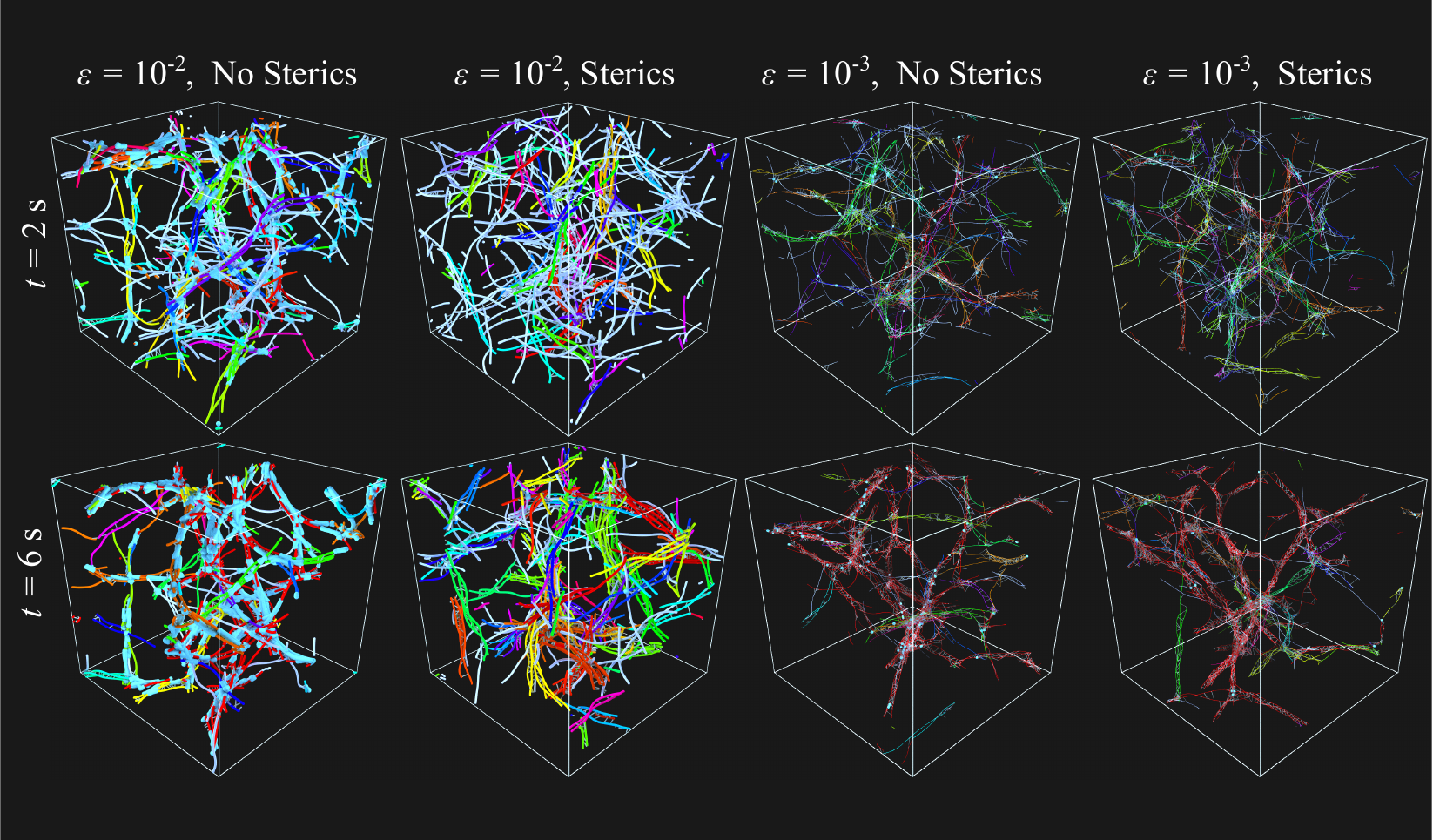}
\includegraphics[width=\textwidth]{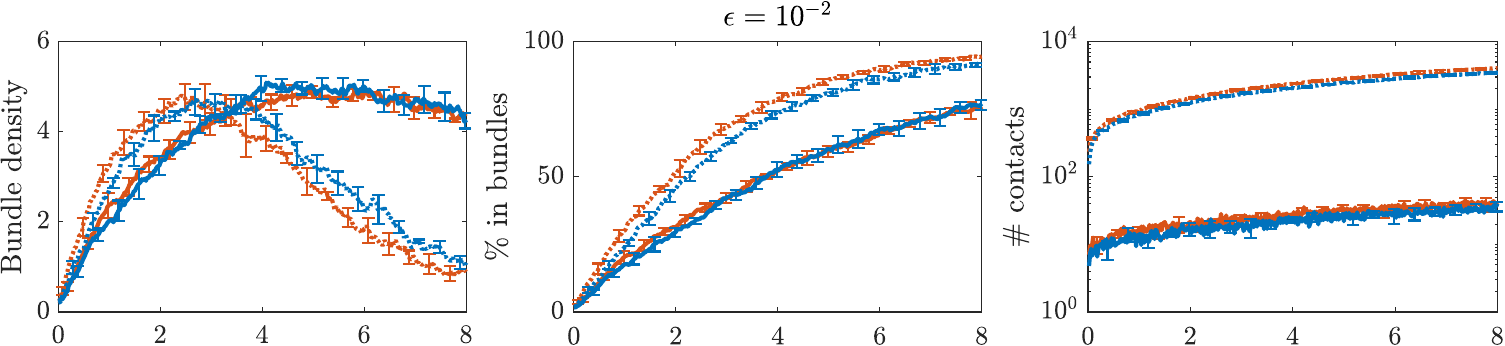}
\includegraphics[width=\textwidth]{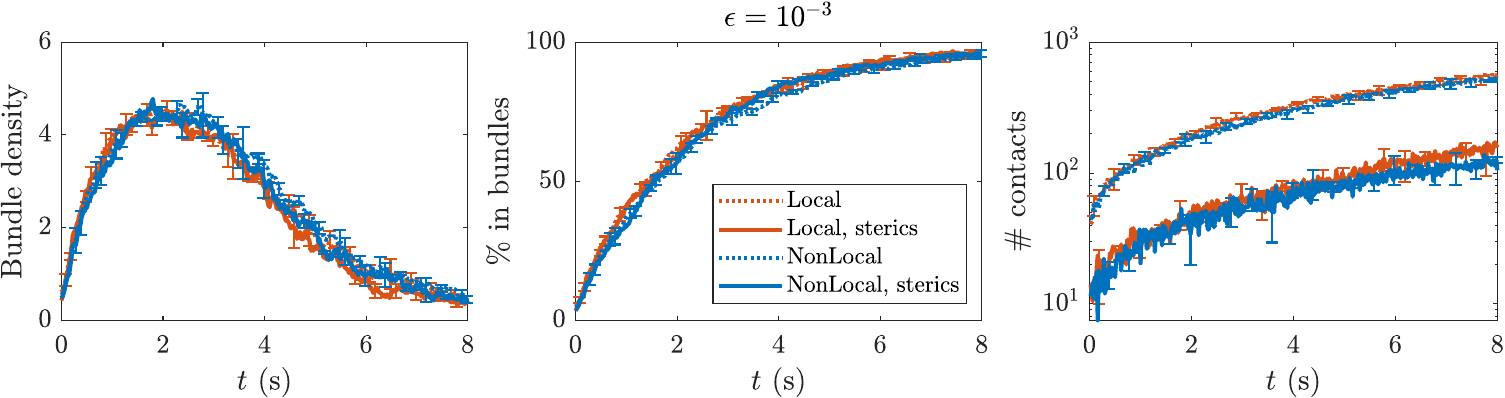}
\caption{\label{fig:BundSnaps}Role of hydrodynamics and steric interactions in bundling. (Top images) Snapshots of bundled fiber configurations with four different simulation conditions. From left to right: $\epsc=10^{-2}$ without sterics, $\epsc=10^{-2}$ with sterics, $\epsc=10^{-3}$ without sterics, and $\epsc=10^{-3}$ with sterics. In all cases we show the results of \emph{local} hydrodynamics only, as those with nonlocal hydrodynamics are indistinguishable by eye. Snapshots are shown at $t=2$ s (top) and $t=6$ s (bottom). Fibers are color-coded by bundle; those not in bundles are shown in white. Thinner white filaments are cross linkers. (Bottom graphs) Statistics of bundling as a function of aspect ratio, hydrodynamic model, and inclusion of steric interactions. For $\epsc=10^{-2}$, the time step sizes are $10^{-4}$ s for non-sterics, $5 \times 10^{-5}$ s for sterics.  For $\epsc=10^{-3}$, the time step sizes are $10^{-4}$ s for local non-sterics, $5 \times 10^{-5}$ s for nonlocal non-sterics, and $10^{-5}$ s for sterics (only 2 trials). }
\end{figure}

Given the degree of agglomeration in the later stages of the bundling process, it seems obvious that steric interactions and hydrodynamic interactions should play a role in determining the timescale of the bundling process. To study this, we perform a systematic study where we vary the hydrodynamic model and steric force potential for both $\epsc=10^{-2}$ and $\epsc=10^{-3}$. Qualitative results and quantitative comparisons are shown in Fig.\ \ref{fig:BundSnaps}, where we display snapshots of the networks at $t=2$ and $t=6$ s (using local hydrodynamics), and compare the bundle density (number of bundles per unit volume), percentage of fibers in bundles, and number of contacts, across the two different aspect ratios and four different simulation conditions.

Our first result can be predicted from the initial snapshots in Fig.\ \ref{fig:BundSeq}. There we clearly see the larger aspect ratio $\epsc=10^{-2}$ taking up much more of the volume than $\epsc=10^{-3}$. It is therefore easy to speculate that steric interactions might substantially change this trajectory. Figure\ \ref{fig:BundSnaps} shows that this is clearly the case; when steric interactions are included, the bundling process is delayed by about a factor of three, with the peak bundle density occuring around $t=5$ s, rather than $t=2$ s. Contrary to previous speculation \cite[Note~23]{muller2014PRL} that steric repulsion would slow down the \emph{initial} part of the bundling process more than the later parts (where fibers are constrained by the CLs), our simulations show that the later stages, where bundles form large agglomerates, are sharply constrained by steric interactions, while the first second or so of bundling appears unchanged. In the case of $\epsc=10^{-2}$, the larger aggregates form structures which are difficult to reorganize into larger connected structures without fiber overlaps, shown in the left column of Fig.\ \ref{fig:BundSnaps} using blue spheres. For $\epsc=10^{-2}$, our steric algorithm is effective at reducing these overlaps without a strong constraint on the time step size, as using a time step size two times smaller than the base parameters removes 99\% of the contacts between fibers.

The role of hydrodynamic interactions in the higher aspect ratio suspensions is intertwined with the role of steric interactions. If we examine the number of contacts in the simulation, there appears to be a slight but noticeable reduction when nonlocal hydrodynamics is included (compare blue and red curves in right panels of Fig.\ \ref{fig:BundSnaps}). Without steric interactions, the reduction in contacts translates into a slower bundling process; in other words, hydrodynamic interactions keep the fibers apart, slowing down bundling (under some simulation conditions, we find that this effect increases the allowable time step size with steric interactions). This effect is natural given the nature of flows: when two fibers move towards each other, they create flows which tend to push them apart, or make it more difficult for them to come together. These flows are extremely short-ranged, however, as the inclusion of steric interactions keeps the fibers far enough apart so as to eliminate their effects (compare solid blue and red curves in the top panel of Fig.\ \ref{fig:BundSnaps}). This demonstrates that nonlocal hydrodynamics has a minor effect on the bundling of filaments. Since there are no long-range flows in the suspension, local steric repulsion can reproduce the effect of nonlocal hydrodynamics.

In the case of $\epsc=10^{-2}$, it is fair to say that the kinetic arrest of the filaments is driven by a combination of steric effects and cross linking. But what happens when we drop to $\epsc=10^{-3}$? The initial snapshots in Fig.\ \ref{fig:BundSeq} gave us a picture of how the fibers look; they appear quite slender, barely taking up any volume. We might consequently expect that cross linking will play more of a role than sterics in constraining their motion. Our more detailed simulations in Fig.\ \ref{fig:BundSnaps} confirm this. While the simulations without sterics have about an order of magnitude more contacts than those with sterics, the number of contacts is still far less than when $\epsc=10^{-2}$, and the inclusion of steric interactions has neither a qualitative (Fig.\ \ref{fig:BundSeq}) nor quantitative (Fig.\ \ref{fig:BundSnaps}) effect on the overall bundling trajectory. Thus, when when the fibers are sufficiently slender, simple simulations based on local hydrodynamics without steric interactions give accurate estimates of the bundling dynamics.

\begin{figure}
\centering
\includegraphics[width=0.49\textwidth]{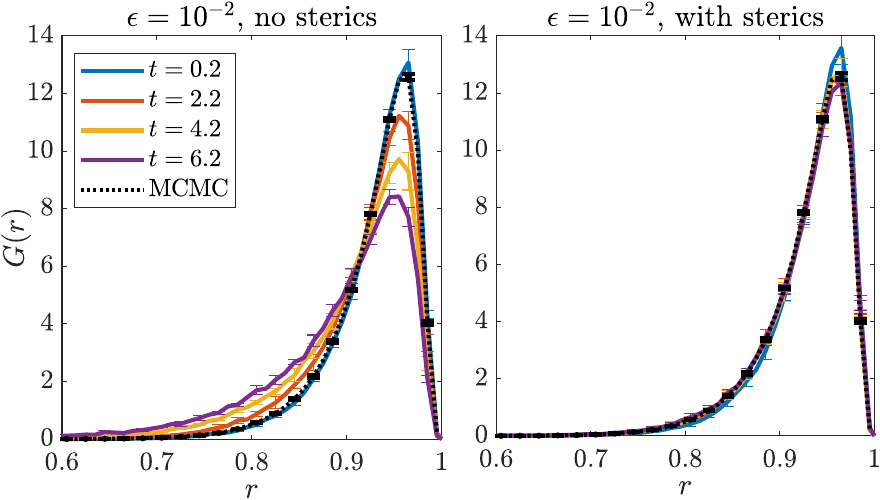}
\includegraphics[width=0.49\textwidth]{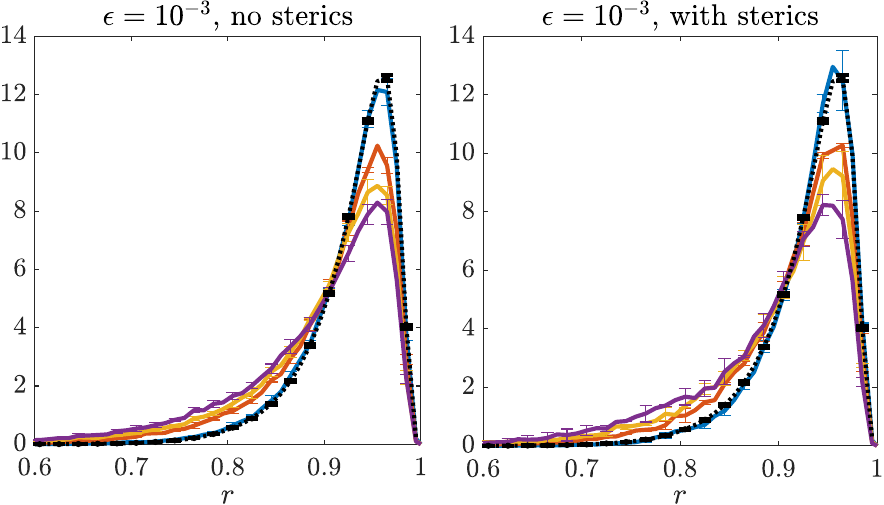}
\caption{\label{fig:EEDistCL}\rev{End-to-end distance distributions for fibers in bundled cross-linked networks. For each of the simulations shown in Fig.\ \ref{fig:BundSnaps} ($\epsc=10^{-2}$ and $\epsc=10^{-3}$, with and without sterics), we plot the filament end-to-end distance distribution over times $t \in (0.1,0.3)$ s (blue), $(2.1,2.3)$ s (red), $(4.1,4.3)$ s (yellow), and $(6.1,6.3)$ s (purple), compared to the equilibrium distribution obtained from MCMC (black, see Section \ref{sec:eqfl}). These plots are generated using simulations with nonlocal hydrodynamics; results for local hydrodynamics are similar.}}
\end{figure}
 
\subsection{\rev{Effect of cross linking on end-to-end distribution}}
\rev{Because the cross linkers generate additional force on the filaments, the filament shapes will no longer be samples from the equilibrium distribution \eqref{eq:GBDist}. To quantify the differences from equilibrium in the presence of strong cross linking, we plot the end-to-end fiber distance over time in Fig.\ \ref{fig:EEDistCL}, considering each of the simulation conditions in Fig.\ \ref{fig:BundSnaps} ($\epsc=10^{-2}$ and $\epsc=10^{-3}$, with and without sterics). As might be expected from the bunding trajectories, a characteristic behavior emerges for all simulation conditions \emph{except} $\epsc=10^{-2}$ with steric repulsion. This behavior is characterized by a gradual widening of the equilibrium distribution, beginning around $t \approx 1$ s (between 0.2 and 2.2), and continuing throughout the simulation. The addition of more filaments with shorter end-to-end distances is the result of strong cross linking effectively ``locking'' large-wavelength deformations into place \cite{mulla2019origin}. Even though the links are transient, the local density of links is sufficiently large as to prevent these deformations from relaxing, and non-equilibrium configurations result.}

\rev{The exception to these dynamics occurs with relatively thick aspect ratios ($\epsc=10^{-2}$) when steric interactions are considered, in which case the fibers retain their equilibrium distribution of end-to-end distances throughout the simulation. While the overall course of bundling is slower with sterics (see Fig.\ \ref{fig:BundSnaps}), comparing simulations at the time of peak bundle density ($t \approx 3.2$ without sterics, $t \approx 4.1$ with sterics) shows that the end-to-end distribution with sterics is still closer to the equilibrium one. Thus, for fixed link density (correlated with bundle density), cross links are less effective in deforming the filaments, since steric interactions keep the thicker filaments sufficiently far apart to allow for the relaxation of \emph{all} bending modes. Since the thicker filaments have diameter (0.02 $\mu$m) on the same order as the cross-linker rest length (0.05 $\mu$m), steric interactions restrict the angle of cross linker binding to be more perpendicular, which induces less deformation.}

\section{Conclusion}
This paper represents the end of a series of works \cite{maxian2021integral,maxian2022hydrodynamics,maxian2023bending} on slender filament hydrodynamics. Our purpose here was two fold: first, we summarized the (often tinkered-with) governing equations and spatial discretization that came out of previous work, so that the reader can have one reference instead of four. This summary revealed an important missing piece that represents the key contribution of this paper: formulating a \emph{nonlocal} mobility that is symmetric positive definite and has cost independent of fiber aspect ratio. Combining this with a novel method for handling steric interactions between filaments, we presented a simulation package which can account for filament slenderness, inextensibility, hydrodynamics, and Brownian motion. We used this package to simulate fibers under gravity (a common example where nonlocal flows are integral to the dynamics) and cross-linked fiber networks (an example from biology which is dominated more by \rev{local} interactions). 

Designing an SPD nonlocal mobility turned out to be a winding journey which alternated between discrete and continuum pictures. At first glance, a continuum mobility based on RPY integrals \eqref{eq:Ucont} ought to generate a discrete SPD mobility, and indeed it does if all of the integrals are resolved with the same number of points. This is where the main problem arises: if all of the integrals must use the same number of points, then the number of points is defined by the self integral, which requires $\mathcal{O}(1/\epsc)$ discretization points to resolve accurately. Attempts to circumvent this difficulty (by using a special quadrature scheme for the self term) destroy the SPD property of the mobility. Using previous methodology, we therefore are able to simulate \emph{either} Brownian motion \emph{or} nonlocal interactions, but not both. 

In this paper, we circumvented this difficulty by making two key observations: first, for a given numerical scheme (oversampling or special quadrature), the eigenvalues of the mobility decrease as the \rev{fibers get fatter}. Second, for a given aspect ratio, the eigenvalues of the mobility decrease as the accuracy of integration increases. Putting these together, we coupled a more accurate quadrature scheme for the self term with \rev{an artificially-fattened fiber radius}. This generated an SPD mobility, which we split into an oversampling step \rev{(for the fattened aspect ratio)}, plus a local correction, which is the special quadrature scheme on the correct aspect ratio minus the special quadrature scheme on the thicker aspect ratio. While this approximation, which comes out of numerical and not theoretical necessity, introduces an error into the dynamics, we demonstrated in both static (Appendix \ref{sec:fatNL}) and dynamic (Sec.\ \ref{sec:grav}) examples that this error is about the same order of magnitude as using $1/\epsc$ points in the self integrals, making it relatively small. In fact, for well-separated filaments the RPY kernel takes the form of a Stokeslet + $\mathcal{O}(\eps^2) \times$ the doublet, and so changing from $\eps$ to $\epsh$ introduces an asymptotic error of order $(\epsRSh)^2$. The thickened aspect ratio $\epsRSh$ represents a degree of freedom in the dynamics; we fixed $\epsRSh=10^{-2}$, believing that $1/\epsRSh=100$ oversampling points per filament would be a good ceiling, but it can drop even further if computationally feasible. Of course, if the filaments are not too slender ($\epsc \approx 10^{-2}$), then the oversampled mobility is sufficient to resolve the dynamics without incurring a large oversampling expense.

In addition to a new formulation of the mobility, we also designed a novel algorithm to treat steric (excluded volume) interactions in filaments discretized with Chebyshev collocation points. The key was to formulate the energy as a double integral, similar to what was done in \cite{popov2016medyan} for straight cylinders. However, in \cite{popov2016medyan}, a singular energy density is chosen, which allows for analytic integration over straight segments. Because our simulations include Brownian fluctuations, fibers might overlap regardless of the strength of the steric potential, and so we chose an energy density that is Gaussian (nonsingular and smooth) in nature, and consequently not possible to integrate analytically. As such, we needed to design efficient quadrature schemes to integrate the density in the limit $\epsc \rightarrow 0$. 

We considered two different quadrature schemes to integrate the steric interaction potential and obtain the steric forces at the Chebyshev collocation points. In the first scheme, we simply resampled the fiber at a set of uniform points, then performed a neighbor search to find pairs of interactions that contribute to the sum. This scheme does not scale well in the limit $\epsc \rightarrow 0$, so we developed a second scheme based on breaking the fibers into $\mathcal{O}(1)$ pieces, then finding and progressively zeroing in on pieces that could be interacting. Similar to the mobility, we validated our more efficient segment-based algorithm by comparing the steric forces it generated to those obtained using oversampling with $1/\epsc$ points (Appendix \ref{sec:StCompare}). 

We chose to treat steric forces explicitly in our temporal update. This leads to a reduction of a factor of roughly 5 to 10 in the allowable time step size, with our original time step leading to unstable simulations when sterics is included. It is tempting in this case to consider an implicit approach to achieve stability at larger time step sizes, such as the complementarity formulation used in \cite{broms2023barrier, yan2019computing, yan2022toward, ferguson2021intersection}, whereby an optimization problem is solved to obtain a configuration that is contact-free while minimizing the contact forces between particles. While these methods do indeed increase the time step size, the method used to solve the optimization problem requires multiple applications of the hydrodynamic mobility, and as such makes each time step at least 10 times more expensive than when contacts are not resolved. Thus employing a more complex contact algorithm would make our simulations just as expensive as treating the steric force explicitly. Furthermore, the complementarity formulation introduces new constraints into the problem, which leads to new drift terms in the Ito Langevin equation that describes the fiber dynamics. Formulating and accounting for these drift terms will be a challenge for any method that attempts to combine Brownian dynamics with contact resolution.

Explicit time integration is more efficient in our case because our volume-exclusion potential is ``soft,'' in the sense that we allow some particle overlaps to keep the time step relatively large. If we were to stiffen up the potential to prevent more overlaps, then we likely would arrive at a point where contact-based algorithms \emph{would} be more efficient. It remains unclear, however, how to extend these algorithms from rigid particles \cite{yan2019computing, ferguson2021intersection, broms2023barrier} to flexible fibers, as the most recent works of Yan et al.\ \cite{yan2019computing,  yan2022toward} treat flexible fibers via a series of rigid rods. The most straightforward way would be to keep the segment-based fiber representation, but replace our special quadrature scheme with segment-based quadrature, \rev{such as that previously developed for regularized Stokeslets} \cite{cortez2018regularized,hall2019efficient}. However, as discussed in the introduction, regularized Stokeslets suffer from a lack of symmetry which destroys the SPD-property of the mobility, so a different regularized singularity (e.g., RPY or FCM) would be required.

We ended this paper by applying our simulation algorithm to two very different examples of fiber dynamics. In the example of fibers under gravity, there were nonlocal fluid flows which entrained Brownian motion; the local and nonlocal mobilities generated qualitatively different dynamics. For cross-linked fiber networks, the only effect of hydrodynamics was to keep the fibers apart, something that we could also account for via steric repulsion. Since our interests are in the latter system, it was a disappointment to see flows making such little difference, but it is no surprise given that the \emph{in vitro} experimental dynamics never demonstrated any large-scale flows \cite{falzone2012assembly}. Indeed, in biology, examples of large-scale fluid flows typically occur when molecular motors are around to actively contract or expand filament networks \cite{wu2024laser,chakrabarti2024cytoplasmic,dutta2024self}. It would be straightforward to expand the cross-linking algorithm presented here (Sec.\ \ref{sec:CLs}) to account for motors, which are essentially moving cross linkers \cite{alberts,quintanilla2023non}. Doing so will allow for a more thorough study of the hydrodynamics of the cytoskeleton, including how actin architecture affects fluid flows \cite{muresan2022f}.

\rev{The platform presented here could also be used to study stress relaxation in cross-linked networks, where previous theories have posited a relationship between cross link activity and stress relaxtion timescales \cite{mulla2019origin, muller2014PRL, broedersz2010PRL,  yao2013PRL}. However, the system stress, which is the quantity of interest in rheology, includes a contribution from Brownian motion \cite[Eq~(3.169)]{doi1988theory}. Classical work showed that the Brownian stress can be computed as a divergence of the matrix mapping velocity to stress, which, when done naively, would require additional resistance solves \cite{brady1993rheological}, significantly increasing the cost of the algorithm presented here. Developing numerical methods for efficient calculation of Brownian stress is consequently an important future direction that would enable large-scale rheological studies.}

\rev{Regarding fiber sedimentation, an avenue for future exploration is to study how a bottom wall impacts the hydrodynamic interactions between fibers. Previous studies on colloidal particles \cite{sprinkle2020driven} have shown how proximity to the wall can impact agglomeration and layering of active colloids (magnetic rollers). Extending such a study to filaments would require efficient evaluation of the hydrodynamic mobility for filaments above a wall. While the corrections to the RPY kernel in this case are known \cite{swan2007simulation}, fast (linear-scaling) methods to compute them remain under development \cite{yan2020scalable, bagge2023fast, hashemi2023computing}. Extending these methods to filaments would again require the development of special quadrature schemes that preserve the overall SPD properties of the RPY kernel. These, and extensions to other non-periodic geometries, represent a challenge for future work.}

\section*{Acknowledgments}
Ondrej Maxian thanks Leslie Greengard for useful discussions about regularization. Ondrej Maxian was supported by the NYU Dean's Dissertation Fellowship and the NSF via GRFP/DGE-1342536, and is currently supported by the Yen and Chicago fellow programs at University of Chicago. This work was also supported by the National Science Foundation through Division of Mathematical Sciences award DMS-2052515, and through a Research and Training
Group in Modeling and Simulation under award RTG/DMS-1646339 (both to A.D.).

\subsection*{Data availability}
All of the python/C++ codes and corresponding input files to reproduce our results
can be found at \url{https://github.com/stochasticHydroTools/SlenderBody}.

\begin{appendices}
\section{RPY Mobility \label{sec:SQ}}
In this appendix, we give more details on the RPY mobility discussed in Section \ref{sec:Mob}. We first discuss the special quadrature scheme for the self term, which was previously presented elsewhere \cite{maxian2022hydrodynamics}, but is repeated here in an effort to make this paper self contained. We then give some more details of the error analysis for the fat-corrected mobility \eqref{eq:MNew}. 

\subsection{Special quadrature scheme}
The self part of the RPY integral mobility \eqref{eq:Ucont} can be rewritten as 
\begin{align}
 \label{eq:transmob}
\V{U}(s) &= 
\int_{D(s)} \left(\Slet{\XPoly(s),\XPoly(s')} +\frac{2\eps^2}{3}\Dlet{\XPoly(s),\XPoly(s')} \right)\V{f}\left(s'\right) \, ds' \\[2 pt] \nonumber &+ \EPMI \int_{D^c(s)} \left(\left(\dfrac{4}{3\eps}-\dfrac{3R\left(\XPoly(s),\XPoly(s')\right) }{8\eps^2}\right)\M{I}+\dfrac{\V{R} \Rhat\left(\XPoly(s),\XPoly(s')\right) }{8\eps^2} \right)\V{f}(s') \, ds',
\end{align}
where $\M{S}$ and $\M{D}$ represent the Stokeslet and doublet for Stokes flow, respectively, and we have used the definition\ \eqref{eq:MbttRPY} to split the integral into a region $D(s)$ for \mbox{$R > 2\eps$} and $D^c(s)$, which uses the RPY kernel for $R \leq 2\eps$. As discussed in Section \ref{sec:SpatDisc}, the fibers as polynomials do not satisfy $\ds \XPoly \cdot \ds \XPoly = 1$ everywhere, and so we say they are \emph{approximately} inextensible. In our quadrature schemes, we therefore do not assume $\norm{\ds \XPoly}=1$, but we still make the approximation $R \approx |s'-s|$ when $R \lesssim 2\eps$, so that 
\begin{equation}
\label{eq:Ddom}
D(s) =\begin{cases} \left(0,s-2\eps\right) \cup \left(s+2\eps,L\right) & 2\eps \leq s \leq L-2\eps \\
\left(s+2\eps,L\right) & s < 2\eps \\
 \left(0,s-2\eps\right) & s > L-2\eps
\end{cases},
\end{equation}
with the complement $D^c(s) = [0,L] \setminus D(s)$. We now discuss each of the integrals in\ \eqref{eq:transmob} in more detail, separating out the Stokeslet and doublet into two separate integrals.

\subsubsection{Stokeslet integral}
For the integral of the Stokeslet in\ \eqref{eq:transmob}, we use a singularity subtraction technique which is closely tied with the asymptotics of the Stokeslet. In particular, we subtract from the integrand the leading order singular behavior and perform that integral separately, which gives
\begin{align}
\label{eq:Stsubtr}
\V{U}^{(S)}=& \int_{D(s)} \Slet{\XPoly(s),\XPoly(s')} \V{f}\left(s'\right)\, ds'\\[2 pt] \nonumber
= \EPMI &\int_{D(s)}\left(\frac{\M{I}+\widehat{\ds{\XPoly}}(s)\widehat{\ds{\XPoly}}(s)}{\norm{\ds{\XPoly}(s)}|s-s'|}\right)\V{f}(s) \, ds'\\ \nonumber+ 
& \int_{D(s)} \left(\Slet{\XPoly(s),\XPoly(s')} \V{f}\left(s'\right) -  \EPMI \left(\frac{\M{I}+\widehat{\ds{\XPoly}}(s)\widehat{\ds{\XPoly}}(s)}{\norm{\ds{\XPoly}(s)}|s-s'|}\right)\V{f}(s)\right) \, ds' \\[2 pt] \nonumber
:=&\V{U}^{(\text{inner, S})}(s)+\V{U}^\text{(int, S)}(s)\\ \nonumber
\label{eq:UinnerSt}
\text{where} \quad 
&\V{U}^{(\text{inner,S})}(s) =
\frac{a_L(s)}{8 \pi \mu}\frac{\left(\M{I}+\widehat{\ds{\XPoly}}(s)\widehat{\ds{\XPoly}}(s)\right)\V{f}(s) }{\norm{\ds{\XPoly}(s)}},
\end{align}
where $a_L(s)$ is given by 
\begin{equation}
a_S(s) = \begin{cases}
\ln{\left(\dfrac{(L-s)s}{4\eps^2}\right)} & 2\eps < s < L-2\eps\\[4 pt]
\ln{\left(\dfrac{(L-s)}{2\eps}\right)}& s \leq 2\eps \\[4 pt]
\ln{\left(\dfrac{s}{2\eps}\right)} & s \geq L-2\eps
\end{cases}
\end{equation}
This leaves a nearly-singular integral that represents the Stokeslet minus the leading order singularity, 
\begin{equation}
\label{eq:StNSint}
\V{U}^\text{(int, S)}(s)=\int_{D(s)} \left(\Slet{\XPoly(s),\XPoly(s')} \V{f}\left(s'\right) -  \EPMI \left(\frac{\M{I}+\widehat{\ds{\XPoly}}(s)\widehat{\ds{\XPoly}}(s)}{\norm{\ds{\XPoly}(s)}|s-s'|}\right)\V{f}(s)\right) \, ds' ,
\end{equation}
which has the same integrand as the finite part integral in slender body theory \cite{ts04}, but with the different domain of integration $D(s)$. Because of our singularity subtraction scheme, the second term in the integral\ \eqref{eq:StNSint} cancels the leading order $1/|s'-s|$ singularity in the first. The next singularity is $\text{sign}(s'-s)$, which means that the near singular integral\ \eqref{eq:StNSint} can be written as
\begin{gather} 
\label{eq:FPre}
\V{U}^\text{(int, S)}(s) =\int_{D(s)} \V{g}_\text{stok}(s,s') \frac{s'-s}{|s'-s|} \, ds'
= \frac{L}{2}\int_{D'(\eta)} \V{g}_\text{stok}(\eta, \eta') \frac{\eta'-\eta}{|\eta'-\eta|} \, d\eta', \qquad \text{where}\\ 
\nonumber
\V{g}_\text{stok}(s,s') = \left(\Slet{\XPoly(s),\XPoly(s')} \V{f}\left(s'\right)  |s'-s| -\frac{1}{8 \pi \mu}\frac{\left(\M{I}+\widehat{\ds{\XPoly}}(s)\widehat{\ds{\XPoly}}(s)\right)\V{f}(s) }{\norm{\ds{\XPoly}(s)}}\right)\frac{1}{s'-s}
\end{gather}
is a smooth function, $\eta=-1+2s/L$ is a rescaled arclength coordinate on $[-1,1]$, and $D'(\eta)=[0,\eta_\ell] \cup [\eta_h,L]$ is defined from $D(s)$. The function $\V{g}_\text{stok}$ is nonsingular at $s=s'$ with the finite limit
\begin{gather}
\label{eq:gttlimit}
\lim_{s' \to s} \V{g}_\text{stok}(s,s') = \EPMI  \Bigg{[} \frac{1}{2\norm{\ds{\XPoly}(s)}^3} \bigg{(} \ds{\XPoly}(s)\ds^2{\XPoly}(s)+\ds^2{\XPoly}(s)\ds{\XPoly}(s)\\ \nonumber
-\left(\ds{\XPoly}(s) \cdot \ds^2{\XPoly}(s)\right) \left(\M{I}+3\widehat{\ds{\XPoly}}(s)\widehat{\ds{\XPoly}}(s)\right)\bigg{)}\V{f}(s) + \left(\frac{\M{I}+\widehat{\ds{\XPoly}}(s)\widehat{\ds{\XPoly}}(s)}{\norm{\ds{\XPoly}(s)}}\right) \ds{\V{f}}(s) \Bigg{]}
\end{gather}
Furthermore, $\V{g}_\text{stok}$ is smooth, so we can express it in a truncated Chebyshev series on $[-1,1]$, 
\begin{equation}
\label{eq:gmono}
\frac{L}{2}\V{g}_\text{stok}(\eta,\eta') \approx \sum_{k=0}^{N-1} \V{c}_k(\eta) T_k(\eta'), 
\end{equation}
where $\V{c}_k$ is a vector of 3 coefficients for each $\eta$. Substituting the Chebyshev  expansion\ \eqref{eq:gmono} into the integrand\ \eqref{eq:FPre}, we obtain
\begin{align}
\label{eq:expandmono}
 \V{U}^\text{(int, S)}(\eta) &=\sum_{k=0}^{N-1} \V{c}_k(\eta) \int_{D'(\eta)} T_k\left(\eta^\prime \right) \frac{\eta'-\eta}{|\eta'-\eta|} \, d\eta' = \sum_{k=0}^{N-1} \V{c}_k(\eta) q^{(S)}_k(\eta)=\V{c}^T(\eta) \V{q}^{(S)}(\eta),\\[6 pt]
\label{eq:qRPY}
\text{where} \qquad 
q_k^{(S)}(\eta) & = \int_{D'(\eta)} T_k(\eta') \frac{\eta'-\eta}{|\eta'-\eta|} \, d\eta'= -\int_{-1}^{\eta_\ell} T_k(\eta') \, d\eta' +\int_{\eta_h}^L T_k(\eta') \, d\eta',
\end{align}
are integrals that can be precomputed to high accuracy for each $\eta$ on the Chebyshev collocation grid. The precomputed integrals in\ \eqref{eq:qRPY} can then be used in an adjoint method to accelerate the repeated evaluation of\ \eqref{eq:expandmono}, so that at each time step only the values of $ \V{g}_\text{stok}(\eta,\eta')$ on the collocation grid need to be determined  \cite{tornquad,maxian2022hydrodynamics}. 

\subsubsection{Doublet integral \label{sec:DbNS}}
The next part of the translational mobility\ \eqref{eq:transmob} is the integral of the doublet kernel. We begin by applying the same singularity subtraction technique we used in\ \eqref{eq:Stsubtr} for the Stokeslet integral
\begin{align}
\label{eq:DbSS}
{\V{U}}^{(D)}(s)=&\int_{D(s)} \Dlet{\V{X}(s),\V{X}(s')} \V{f}\left(s'\right)\, ds'\\[2 pt] \nonumber
=  \EPMI & \int_{D(s)} \left(\frac{\M{I}-3\widehat{\ds{\XPoly}}(s)\widehat{\ds{\XPoly}}(s)}{\norm{\ds{\XPoly}(s)}^3|s-s'|^3}\right)\V{f}(s)\\ \nonumber + &\int_{D(s)} \left(\Dlet{\V{X}(s),\V{X}(s')} \V{f}\left(s'\right) -\ \EPMI \left(\frac{\M{I}-3\widehat{\ds{\XPoly}}(s)\widehat{\ds{\XPoly}}(s)}{\norm{\ds{\XPoly}(s)}^3|s-s'|^3}\right)\V{f}(s)\right) \, ds'\\[2 pt] \nonumber
 =& {\V{U}}^{(\text{inner, D})}(s)+{\V{U}}^\text{(int, D)}(s).
\end{align}
The term ${\V{U}}^{(\text{inner, D})}(s)$ is given for all $s$ by
\begin{gather}
\label{eq:UinnerDb}
{\V{U}}^{(\text{inner,D})}(s) =  
\frac{\left(\M{I}-3\Xs(s)\Xs(s)\right)a_D(s)\V{f}(s)  }{8 \pi \mu \norm{\ds{\XPoly}(s)}^3},
\end{gather}
where $a_D(s)$ is defined by
\begin{equation}
a_D(s) = \begin{cases}
\dfrac{1}{4 \eps^2}-\dfrac{1}{2s^2}-\dfrac{1}{2(L-s)^2}& 2\eps < s < L-2\eps\\[4 pt]
\dfrac{1}{8 \eps^2}-\dfrac{1}{2(L-s)^2}& s \leq 2\eps \\[4 pt]
\dfrac{1}{8 \eps^2}-\dfrac{1}{2s^2}& s \geq L-2\eps
\end{cases}\\[4 pt] \nonumber
\end{equation}

The nearly-singular doublet integral is handled in exactly the same way as the Stokeslet one, but the singular behavior is different, namely,
\begin{gather}
 {\V{U}}^\text{(int, D)}(s) =\int_{D(s)} \frac{(s'-s)}{|s'-s|^3} \, \V{g}_D(s',s) \, ds'= \frac{2}{L}\int_{D'(\eta)} \V{g}_\text{D}(\eta, \eta') \frac{(\eta'-\eta)}{|\eta'-\eta|^3} \, d\eta',\\[4 pt] 
\label{eq:dint} \nonumber
\text{where} \qquad \V{g}_D(s',s)=\left(\Dlet{\V{X}(s),\V{X}(s')} \V{f}\left(s'\right) -\EPMI \left(\frac{\M{I}-3\widehat{\ds{\XPoly}}(s)\widehat{\ds{\XPoly}}(s)}{\norm{\ds{\XPoly}(s)}^3|s-s'|^3}\right)\V{f}(s)\right) \frac{|s'-s|^3}{(s'-s)}
\end{gather}
is a smooth function with the finite limit
\begin{gather*}
\lim_{s'\rightarrow s} \V{g}_D(s',s) = \EPMI  \Bigg{[} \frac{1}{2\norm{\ds{\XPoly}(s)}^5} \bigg{(} -3\left(\ds{\XPoly}(s)\ds^2{\XPoly}(s)+\ds^2{\XPoly}(s)\ds{\XPoly}(s)\right)\\ -\left(\ds{\XPoly}(s) \cdot \ds^2{\XPoly}(s)\right) \left(3\M{I}-15\widehat{\ds{\XPoly}}(s)\widehat{\ds{\XPoly}}(s)\right)\bigg{)}\V{f}(s) + \left(\frac{\M{I}-3\widehat{\ds{\XPoly}}(s)\widehat{\ds{\XPoly}}(s)}{\norm{\ds{\XPoly}(s)}^3}\right) \ds{\V{f}}(s) \Bigg{]}.
\end{gather*}
If we expand $\V{g}_D$ in terms of Chebyshev polynomials, the efficient evaluation of\ \eqref{eq:DbSS} requires precomputing integrals of the form
\begin{equation}
\label{eq:intMod}
q_k^\text{(D)}(\eta) =\int_{D'(\eta)} \frac{(\eta'-\eta)}{|\eta'-\eta|^3} T_k\left(\eta'\right) \, d\eta',
\end{equation}
for all $\eta$ on the Chebyshev grid. Notice that the integrals\ \eqref{eq:intMod} are not defined if $\eta'=\eta$ is included in the integration domain, but $D'(\eta)$ does \emph{not} contain $\eta$.

\subsubsection{Integral on $R < 2\eps$}
After discretizing the Stokeslet and doublet integrals, we are left with the integral over $D^c(s)$ in the third line of\ \eqref{eq:transmob}. The integrand is nonsingular, but behaves like $|s-s'|$, and so for each $s$ we split the domain $D^c$ into $(s-2\eps, s)$ and $(s,s+2\eps)$ (with appropriate modifications at the endpoints), and use $N_2/2$ Gauss-Legendre quadrature points to sample the fiber and force density (i.e., sample the Chebyshev interpolant of each) and evaluate the integral on each of the two subdomains. We use $N_2/2$ points for each of these integrals so that there are a total of $N_2$ additional (local) quadrature nodes per collocation point. We previously found $N_2=4$ to be sufficient when $\epsc < 10^{-3}$; otherwise we set $N_2=8$ \cite{maxian2022hydrodynamics}.

\subsection{Empirical error analysis for the fat-corrected mobility \label{sec:fatNL}}
\begin{figure}
\centering
\includegraphics[width=0.5\textwidth]{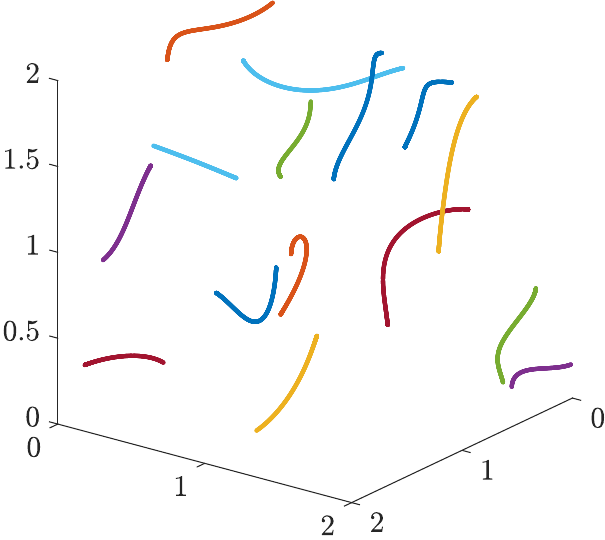}
\caption{\label{fig:FibShap}Fiber shapes for the deterministic test of the mobility. 15 fibers are shown.}
\end{figure}

To perform a convergence study and quantify errors in the mobility, smooth fiber shapes and forces are necessary. At the same time, we want this test to be realistic, so that the fiber shapes and forces roughly match those generated from the equilibrium distribution \eqref{eq:GBDist}. To obtain fibers that satisfy both conditions, we perform MCMC sampling (as previously described \cite{maxian2023bending}) to generate semiflexible fibers with $N_x=5$ Chebyshev nodes from the Gibbs-Boltzmann equilibrium distribution \eqref{eq:GBDist} (Fig.\ \ref{fig:FibShap} shows some examples). For each fiber, we compute the elastic force density $\V{f}=\Wt^{-1}\V{F}^{(\kappa)}$ (see \eqref{eq:FKappa}) using the grid with $N_x=5$ nodes, which defines a fourth order Chebyshev polynomial $\V{f}(s)$ for the continuum force density in the mobility \eqref{eq:Ucont}. Because there are only $N_x=5$ Chebyshev nodes, both the fiber shapes and forces are well resolved, and the only error left is in computing the RPY integrals with the near-singular RPY kernel. 

\subsubsection{Accuracy of oversampled mobility for the self term}
\begin{figure}
\centering
\subfigure[Self velocity]{\includegraphics[width=0.5\textwidth]{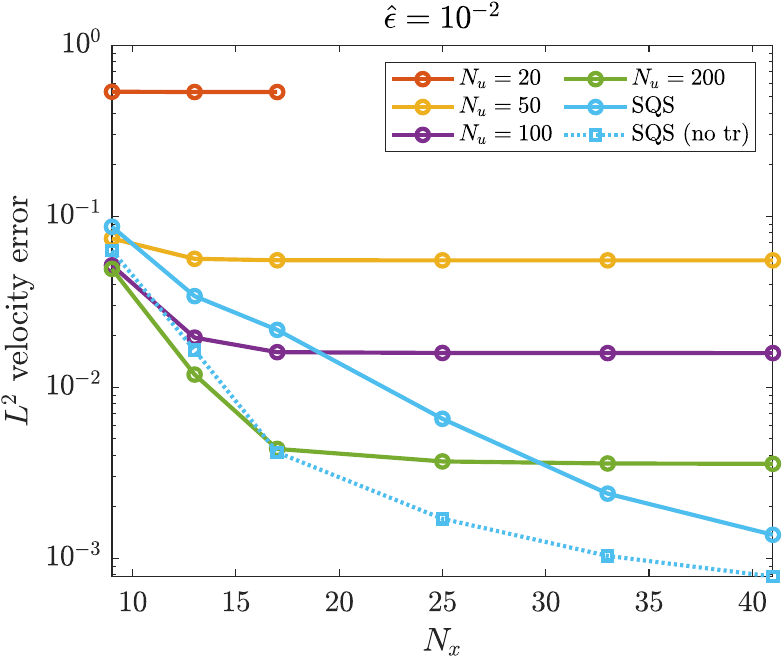}
\includegraphics[width=0.48\textwidth]{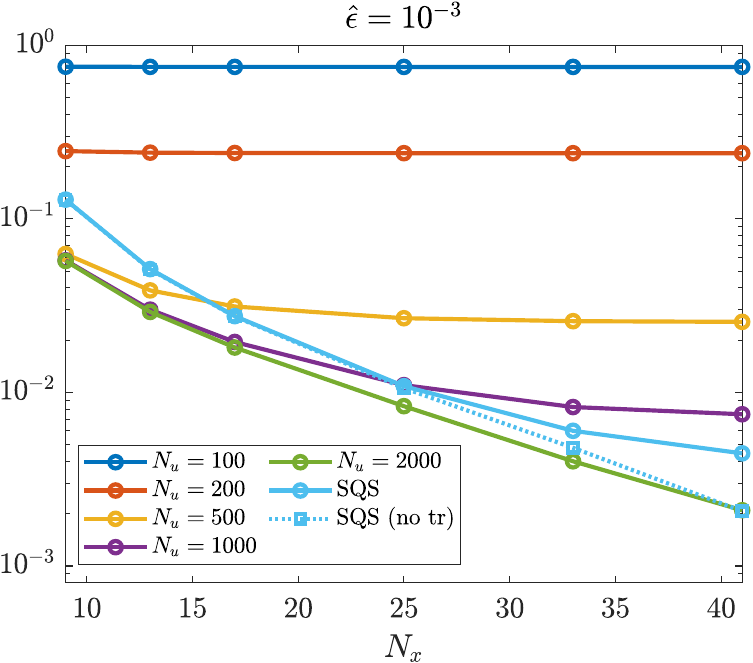}}
\subfigure[Nonlocal velocity]{\includegraphics[width=0.49\textwidth]{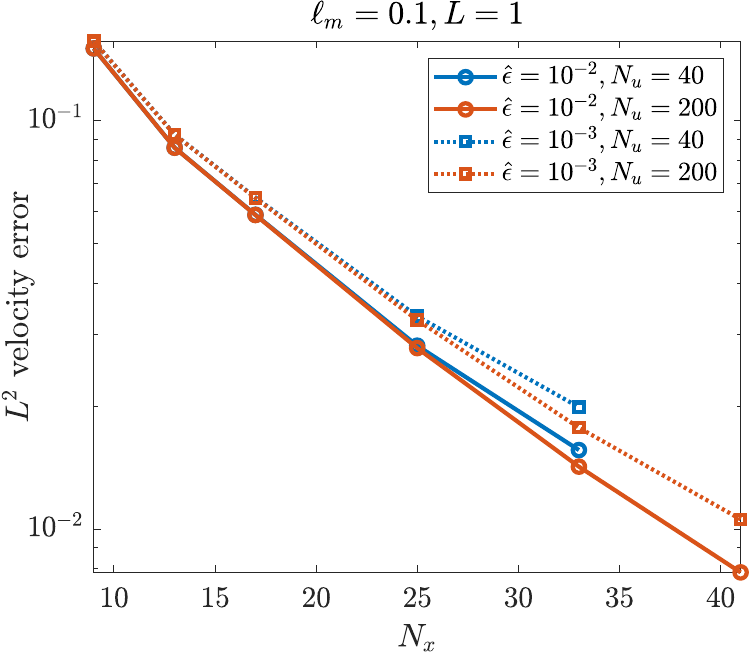}
\includegraphics[width=0.49\textwidth]{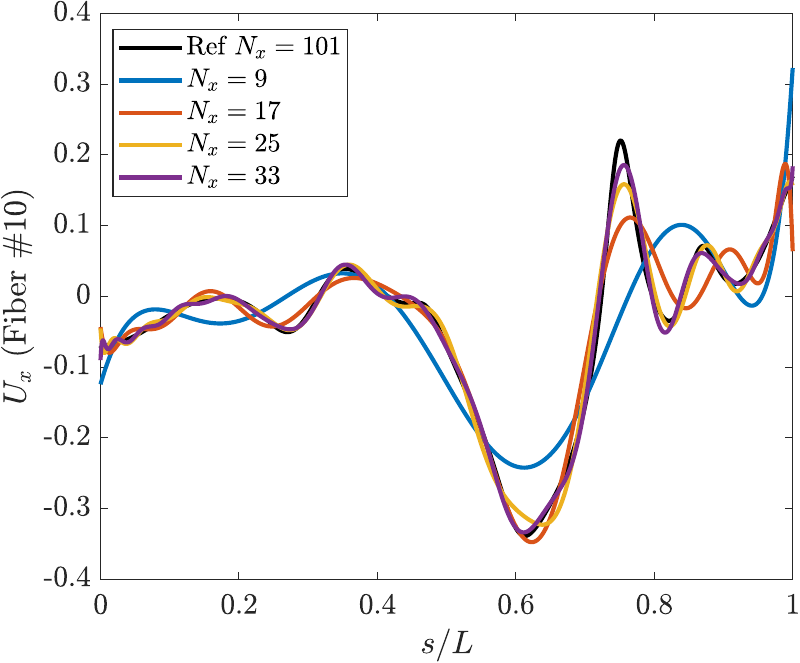}}
\caption{\label{fig:SelfEr}Errors in the fiber velocity using the upsampled mobility \eqref{eq:Mref}. (a) Relative $L^2$ errors in the \emph{local} (self) velocity on a single fiber (first term in \eqref{eq:Ucont}), averaged over ten filaments, as a function of the number of Chebyshev nodes $N_x$. Blue ($N_u=0.1/\epsRS$), red (0.2/$\epsRS$), yellow (0.5/$\epsRS$), purple (1/$\epsRS$), and green (2/$\epsRS$) lines show the errors using mobility\ \eqref{eq:Mref} with increasing numbers of oversampled points; cyan lines show the result with special quadrature (dotted lines show the corresponding result without eigenvalue truncation). (b) Relative errors in nonlocal velocity (second term in \eqref{eq:Ucont}) in a mesh size $\ell_m/L=0.1$. Left panel: $L^2$ velocity errors for $N_u=40$ (blue) and $N_u=200$ (red) for $\epsRS=10^{-2}$ (solid lines with circles) and $\epsRS=10^{-3}$ (dotted lines with squares). Right panel:  the \emph{nonlocal} $x$ velocity field on fiber \#10. The reference is shown with the result computed with varying $N_x$ and $N_u=200$.   }
\end{figure}

We begin by studying the errors in the self velocity. We construct a reference velocity solution using $N_x=81,201$ Chebyshev nodes and $N_u=1000,4000$ oversampling points (for $\epsRS=10^{-2}$ and $\epsRS=10^{-3}$, respectively). Then, we compute the velocity using the reference mobility \eqref{eq:Mref} with varying numbers of Chebyshev nodes and oversampled points. We use ten fibers and compute the mean $L^2$ error in the velocity (on an upsampled grid), divided by the $L^2$ norm of the true velocity. Figure \ref{fig:SelfEr}(a) shows the errors in the self velocity field with the oversampled mobility are a strong function of the number of oversampling points. In fact, the errors in the mobility are \emph{dominated} by the number of oversampling points $N_u$, and 2 digits of accuracy can only be obtained when $N_u=1/\epsRS$. 

\subsubsection{Accuracy of special quadrature for the self term}
We now examine the errors in the self-velocity using the special quadrature scheme (cyan lines in Fig.\ \ref{fig:SelfEr}(a)). Here we set ${\lambda^*}$ as the smallest eigenvalue of the reference mobility \eqref{eq:Mref} for the same number of collocation nodes and $1/\epsRS$ upsampling points. Unlike the errors in oversampled quadrature, which for any $N_u$ eventually saturate as $N_x$ increases, the errors using special quadrature consistently drop to about 3 digits of accuracy for $N_x=41$, demonstrating spectral convergence. Thus, the spectral quadrature scheme gives an accurate self velocity field with cost independent of $\epsc$. \rev{The corresponding errors without eigenvalue truncation are even smaller (dotted lines), though the eigenvalue truncation is clearly a secondary effect.}

\subsubsection{Errors in the nonlocal velocity}
We contrast the behavior of the self mobility, which is dominated by the number of oversampling points $N_u$, to that of the nonlocal mobility, which is not sensitive to $N_u$. To perform a convergence study, we set up 800 non-overlapping fibers of length $L=1$ (we use a radius $\eps= 10^{-2}$ for the no overlap condition) in a box of size $L_d=2$ on all sides (the mesh size is therefore $\ell_m=0.1$). For $\epsRS=10^{-2}$ and $10^{-3}$, we first compute the nonlocal flows (subtracting the self terms) with $N_x=101$ and $N_u=500$ to establish a reference solution, then vary $N_u$ and $N_x$ and plot the $L^2$ velocity error (the mean over the first ten filaments, normalized by the self velocity of those filaments). So that the normalization is the same across both radii, we normalize by the average self velocity for $\epsRS=10^{-2}$. 

Figure \ref{fig:SelfEr}(b) shows the results. Unlike for the self term, the errors in the nonlocal terms are roughly independent of aspect ratio and oversampling factor, and depend only on the number of collocation points, which must be large to resolve the nonsmooth nonlocal velocity field (Fig.\ \ref{fig:SelfEr}(b)). We thus have a contradiction in the mobilities: for the local mobility, accuracy is controlled by the number of upsampled points proportional to $1/\epsRS$, while for the nonlocal mobility, the cost to obtain a given accuracy is independent of $\epsRS$ and $N_u$ and a function of the number of collocation points alone. Thus, in the reference mobility \eqref{eq:Mref}, the self term necessarily dominates the computational cost by setting a high floor on the number of oversampled points $N_u$. 

\subsubsection{Errors using the fat corrected mobility}
\begin{figure}
\centering
\includegraphics[width=\textwidth]{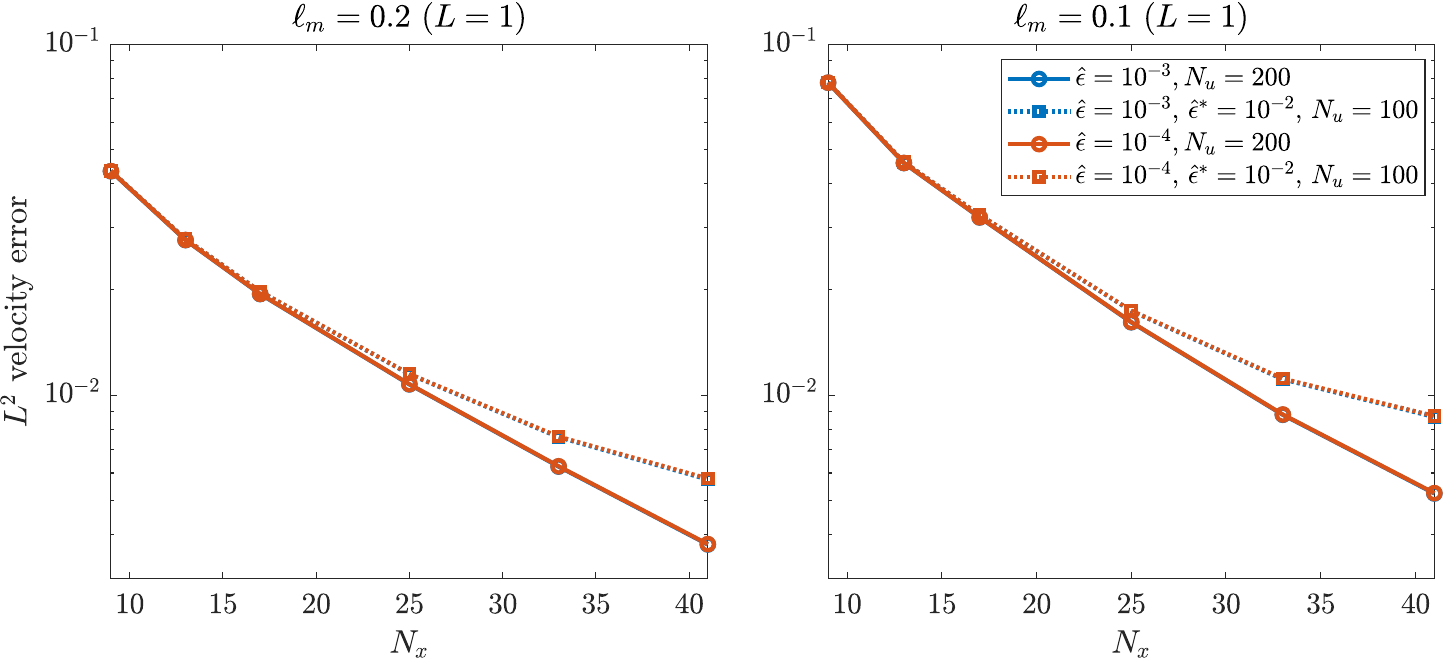}
\caption{\label{fig:NLStar}Errors in the nonlocal velocity with the fat-corrected mobility \eqref{eq:MNew}. We compare the error with the true $\epsRS$ (using various numbers of collocation points and $N_u=200$) with the error using a larger value of $\epsRSh$ and $N_u=100$. Regardless of the mesh size, the error from using a larger aspect ratio is secondary to discretization error, and both errors are independent of aspect ratio (because the $\mathcal{O}(\epsRS^2)$ term is negligible at these aspect ratios, these errors are approximately equal to the errors in integrating the line of Stokeslets). }
\end{figure}

 In Fig.\ \ref{fig:NLStar}, we use the fat-corrected mobility \eqref{eq:MNew} and examine the nonlocal velocity errors for two different mesh sizes $\ell_m=0.2$ and $\ell_m=0.1$. We fix $\epsRS=10^{-3}$ or $10^{-4}$ and compute reference nonlocal flows (subtracting the self terms) with $N_x=101$ and $N_u=500$. We then recompute these flows with the \emph{true} aspect ratio with smaller values of $N_x$, setting $N_u=200$ (recall from Fig.\ \ref{fig:SelfEr}(b) that $N_x$ is the dominant contributor to the error), and plot the errors using solid lines in Fig.\ \ref{fig:NLStar} (the normalization is the mean local velocity when $\epsRS=10^{-3}$; it is convenient to use the same normalization for both quantities). We then repeat the calculation of nonlocal velocity, but with ``fattened'' filaments, i.e., we use $\epsRSh=10^{-2}$ to compute the nonlocal velocity, but still compare to the true answer with the correct value of $\epsRS$. Figure \ref{fig:NLStar} shows that the additional error incurred in using this approximation is small relative to the discretization errors; it does not become noticeable until $N_x=33$, when we are already achieving 2 digits of accuracy. As might be intuitively expected, the saturated errors in this approximation increase when the mesh size decreases. 

\begin{figure}
\centering
\includegraphics[width=0.5\textwidth]{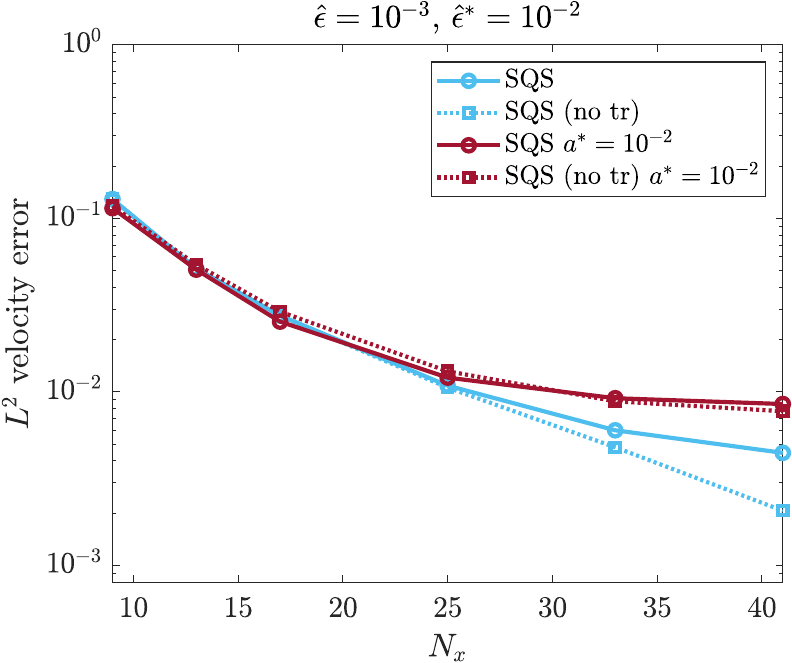}
\caption{\label{fig:SelfErFat}Self errors in the mobility using $\epsRSh=10^{-2}$ for a true aspect ratio $\epsRS=10^{-3}$. We repeat the test of Fig.\ \ref{fig:SelfEr} (with special quadrature results plotted in cyan), but this time also plot the errors using the new mobility \eqref{eq:MNew} in dark red. This adds the difference between oversampled and special quadratures for the larger $\epsRSh=10^{-2}$ to the self error term. We do not see any effect until we drop below two digits of accuracy.  }
\end{figure}

\rev{Finally, Fig.\ \ref{fig:SelfErFat} reveals that the errors in the self velocity are not strongly affected by using the fattened mobility \eqref{eq:MNew}. As with the nonlocal calculation, the errors in the self velocity (which come from adding and subtracting two different calculations of the self velocity integral in \eqref{eq:Ucont}) do not appear noticeable until $N_x \approx 30$, at which we already obtain two digits of accuracy.}

\section{Steric interactions \label{sec:StericAppen}}
In this appendix, we give more details on the steric interaction algorithm we presented in Section \ref{sec:SterSeg}. There are four main points to clarify: the errors in approximating curved fiber pieces as straight segments, the Newton solver for curved segments, the approximation of the fibers as quadratic curves to obtain the integration domains, and the union of integration domains on nearly interacting fibers. We conclude by studying the accuracy of the steric forces in the case of bundling fibers, finding that our segment-based algorithm gives the same accuracy as oversampling in estimating the steric forces.

\subsection{Errors from straight segments \label{sec:ClosePts}}

One difficult part of the segment-based algorithm is step 3, where we need to estimate the curvature of the fibers and use Newton's method to obtain the closest points on the two segments. This step is actually strictly necessary to preserve independence from $\epsc$. If we use segments (dashed lines in Fig.\ \ref{fig:SegSch}(b)), we will have some error in estimating the distance between curved fiber pieces. Since $\rmax \sim \rc$, this error will matter more when $\epsc$ is smaller, and so there will be more erroneous pairs of fiber pieces as $\epsc$ decreases (assuming the curvature of the fibers does not change with $\epsc$). Thus we need to eliminate errors by performing the nonlinear optimization\ \eqref{eq:NewtonMin}.

This section deals with how we estimate the curvature of the fiber segments relative to their straight backbones. Figure\ \ref{fig:SegSch}(b) shows what we mean by this: the pieces of the fibers $\XPoly(s)$ that we are concerned with are shown as solid lines, while the straight segment approximations are shown as dotted lines. We define a distance $\Delta h(s)$ that is simply the distance between the straight segment and the curved fiber. Our goal is to get an estimate for the maximum $\Delta h$ along the segment. We expect this distance to depend on the number of segments we put in each persistence length, which is the ratio $\ell_p/\Lseg$. 

\begin{figure}
\centering
\includegraphics[width=\textwidth]{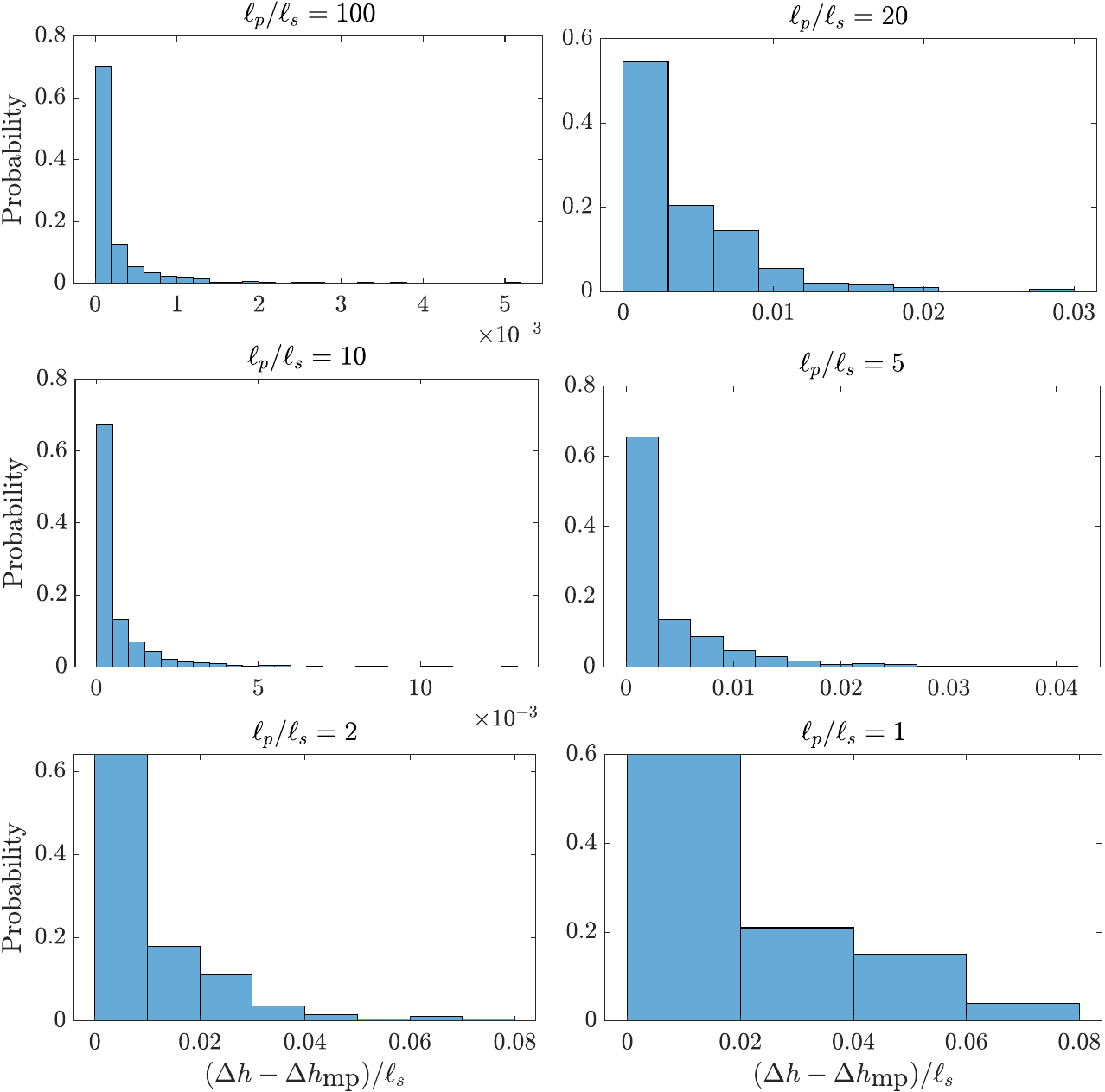}
\caption[Error in estimating deviation from straight segments using segment midpoint]{\label{fig:StraightErMP}Error in estimating $\Delta h$ using its value at the midpoint of the segment. When there are at least 2 panels per persistence length the error is typically 5\% of the segment length or less.}
\end{figure} 

There are two possible ways to do this. First, we can pretabulate, for each $\ell_p/\Lseg$, the expected deviation from straightness, and use that as $\Delta h$. This is a good idea for equilibrium simulations, but in reality other external forces could bend the fibers, so we should really do something with the actual fiber shape. The simplest thing we can do is sample the fiber at the segment midpoint, and then use this to estimate the deviation from straight. For fibers sampled from the equilibrium distribution, Fig.\ \ref{fig:StraightErMP} shows that we can use segments of size $\ell_p$ and still make an error of no more than 10\% in doing this. Therefore, we just sample the segment at the midpoint, compute $\Delta h_\text{mp}$ (at the fiber midpoint), and set $\Delta h=\Delta h_\text{mp}$ in step 3. 

\subsection{Newton solve \label{sec:Newton}}
In this section, we discuss the modified Newton method we use to solve the minimization problem\ \eqref{eq:NewtonMin}. We denote the objective function as $D\left(\ind{s}{i},\ind{s}{j}\right)$ and define the gradient and Hessian matrix
\begin{gather}
\label{eq:NewtG}
\V{g}=\V{\nabla} D = \begin{pmatrix} 
2\ind{\XPoly}{i} \cdot \ds{\ind{\XPoly}{i}}-2 \ds{\ind{\XPoly}{i}} \cdot \ind{\XPoly}{j}\\
2\ind{\XPoly}{j} \cdot \ds{\ind{\XPoly}{j}}-2 \ds{\ind{\XPoly}{j}}\cdot \ind{\XPoly}{i}
\end{pmatrix}\left(\ind{s}{i},\ind{s}{j}\right) \\ \nonumber
\M{H} = \V{\nabla}^2 D.
\end{gather}
Starting with initial guess $\ind{s}{i}=\ind{s}{i}_\text{seg}$ and $\ind{s}{j}=\ind{s}{j}_\text{seg}$ (the closest points between the straight segments), we then perform the following projected Newton iteration. Denoting $\V{s}=\left(\ind{s}{i},\ind{s}{j}\right)$, to perform step $k$ we \cite{bertsekas1982projected, kim2010tackling}
\begin{enumerate}
\item Evaluate $\V{g}_k=\V{g}\left(\V{s}_{k}\right)$ and $\M{H}_k=\M{H}\left(\V{s}_k\right)$.
\item Compute the eigenvalue decomposition of $\M{H}_k=\M{\Psi} \M{\Lambda}\M{\Psi}^{-1}$ to verify positive definiteness. Letting $\M{\Lambda}=\text{diag}\left(\lambda_{1},\lambda_{2}\right)$, we abort if both $\lambda_{1}$  and $\lambda_{2}$ are negative (and report a failure of Newton's method). If $\lambda_{2} < \lambda_{1}/\kappa_\text{max}$, then we set $\lambda_{2}^* \leftarrow \kappa_\text{max}\lambda_1$; otherwise $\lambda_{2}^*=\lambda_{2}$. Here $\kappa_\text{max}$ represents the maximum condition number, which is 1 for gradient descent; we find that a larger $\kappa_\text{max}$ gives faster convergence, so we use $\kappa_\text{max}=10^6$. We compute the modified Hessian $\widetilde{\M{H}}_k$ which is symmetric positive definite, and compute the modified inverse Hessian $\widetilde{\M{H}}_k^{-1}=\M{\Psi}\text{diag}\left(\lambda_1^{-1},\left(\lambda_2^*\right)^{-1}\right)\M{\Psi}^{-1}$.
\item We define a boundary set of indices 
$$B = \left\{q: \left(\ind{s}{q} \leq \varepsilon \quad \text{and} \quad \ind{\V{g}}{q}> 0\right)\quad  \text{or}  \quad \left(\ind{s}{q} \geq L-\varepsilon \quad \text{and} \quad \ind{\V{g}}{q} < 0\right)\right\},$$
which are the indices whose descent directions would leave the domain $[0,L] \times [0,L]$ (we use $\varepsilon=10^{-10}$ for a double precision Newton solver). For all indices $q$ in $B$, we set the $q$th column and $q$th row of $\widetilde{\M{H}}^{-1}$ to zero, and set the $\ind{\V{g}}{q}=0$ as well. 
\item Compute the descent direction $\V{p}_k =-\widetilde{\M{H}}_k^{-1}\V{g}_k$. If $\norm{\V{p}_k} > 0.1L$, renormalize so that $\norm{\V{p}_k}=0.1L$ (this step is necessary to prevent a long line search). 
\item Take a projected step along the free variables only
$${\V s}_{k+1} = \M{P}_{\left[0,L\right]} \left(\V{s}_{k}+\alpha_k \V{p}_k\right) $$
where the relative step size $\alpha_k$ is the largest number $2^{-n}$ that satisfies the Armijo condition
$$D\left(\V{s}_{k}\right)-D\left(\V{s}_{k+1}\right) \geq -\frac{1}{2}\alpha_k \V{g}_k^T \V{p}_k,$$
and $\M{P}_{\left[0,L\right]}$ projects the step onto $[0,L] \times [0,L]$.
\end{enumerate}
This iteration repeats until $\norm{\V g}$ is less than the tolerance, set to $0.01\left(\delta/L\right)$ in the results that follow. Because we zero out the gradient in the constrained directions in step 3, the Armijo and convergence criteria are one-dimensional when one of the directions is constrained.

\subsection{Quadratic approximation to obtain integration domain \label{sec:quadApprox}}
Once we identify the closest points on the two fibers, which are at coordinates $\ind{s}{i}_*$ and $\ind{s}{j}_*$, we identify $\V{d}_*=\ind{\XPoly}{i}_*-\ind{\XPoly}{j}_*$, where the notation $\ind{\XPoly}{i}_*$ means ${\ind{\XPoly}{i}}\left(\ind{s}{i}_*\right)$. If the minimum distance $d_*=\norm{\V{d}_*}$ is less than $\rmax$, we identify the domains on the fibers over which they could be interacting within distance $\rmax$. We do this by approximating the fibers via a quadratic interpolant, $\ind{\XPoly}{i}(s) \approx \ind{\XPoly}{i}\left(\ind{s}{i}_*\right)+\ds{\ind{\XPoly}{i}}\left(\ind{s}{i}_*\right)\ind{\D s}{i} + \ds^2{\ind{\XPoly}{i}}\left(\ind{s}{i}_*\right)\left(\ind{\D s}{i}\right)^2/2$, so that the slack distance between the fibers $\Delta d^2 = d^2 - \rmax^2$ is roughly (dropping terms of order $\D s^3$ or larger), 
\begin{gather}
\label{eq:quadDist}
\Delta d^2\left(\ind{s}{i}, \ind{s}{j}\right) \approx a \left(\ind{\D s}{i}\right)^2 + b \left(\ind{\D s}{j}\right)^2+2c\ind{\D s}{i} +2d\ind{\D s}{j}+2e\ind{\D s}{i}\ind{\D s}{j}+f\\  \nonumber
a = \left(\ds{\ind{\XPoly}{i}}_* \cdot \ds{\ind{\XPoly}{i}}_*\right)+\V{d}_* \cdot \ds^2\ind{\XPoly}{i}_*,  \quad 
b = \left(\ds{\ind{\XPoly}{j}}_* \cdot \ds{\ind{\XPoly}{j}}_*\right)-\V{d}_* \cdot \ds^2\ind{\XPoly}{j}_* \\ \nonumber
c = \V{d}_* \cdot  \ds{\ind{\XPoly}{i}}_*, \quad d = -\V{d}_* \cdot \ds\ind{\XPoly}{j}_*, \quad e = -\ds{\ind{\XPoly}{i}}_* \cdot \ds{\ind{\XPoly}{j}}_*, \quad f = d_*^2 -\rmax^2 < 0.
\end{gather}

The solution of $\Delta d^2 =0$ (equivalently, $d^2=\rmax^2$) for $\D \ind{s}{i}$ and $\D \ind{s}{j}$ is an ellipse. The non-parametric curves that bound the ellipse are given by
\begin{gather*}
\D \ind{s}{j}_{\pm}\left(\D \ind{s}{i}\right)=\frac{1}{2b}\left(-2d-2e\D \ind{s}{i}\pm \sqrt{\left(2d+2e\D \ind{s}{i}\right)^2-4\left(a\left(\D \ind{s}{i}\right)^2+2c\D \ind{s}{i}+f\right)b}\right) \\ \nonumber
\D \ind{s}{i}_{\pm}\left(\D \ind{s}{j}\right)=\frac{1}{2a}\left(-2c-2e\D \ind{s}{j} \pm 
\sqrt{\left(2c+2e\D \ind{s}{j}\right)^2-4\left(b\left(\D \ind{s}{j}\right)^2+2d\D \ind{s}{j}+f\right)a}\right) 
\end{gather*} 
We are interested in the bounding box that we can draw around the ellipse in $\left(\D \ind{s}{i},\D \ind{s}{j}\right)$ space. The points where the bounding box touches the ellipse are where the discriminants are zero, and are given by
\begin{gather}
\label{eq:s1roots} \nonumber
\left(2d+2e\D \ind{s}{i}\right)^2-4\left(a\left(\D \ind{s}{i}\right)^2+2c\D \ind{s}{i}+f\right)b=0 \\ 
\D \ind{s}{i}_{\pm} = \frac{8cb-8de\pm\sqrt{(8de-8cb)^2-4(4e^2-4ab)(4d^2-4bf)}}{8e^2-8ab}\\ \nonumber
\left(2c+2e\D \ind{s}{j}\right)^2-4\left(b\left(\D \ind{s}{j}\right)^2+2d\D \ind{s}{j}+f\right)a=0 \\ \nonumber
\D \ind{s}{j}_{\pm} = \frac{8ad-8ce\pm\sqrt{(8ce-8ad)^2-4(4e^2-4ab)(4c^2-4af)}}{8e^2-8ab}
\end{gather}
We then set $\D \ind{s}{i} = \text{max}\left(|\D \ind{s}{i}_+|, |\D \ind{s}{i}_-|\right)$, and likewise for $\D \ind{s}{j}$, then form an interval $\ind{S}{i}=\left(\max{\left(0,\ind{s}{i}_*-\D \ind{s}{i}\right)}, \min{\left(L,\ind{s}{i}_*+\D \ind{s}{i}\right)}\right)$, and likewise for $\ind{S}{j}$. The pieces of the fiber $\ind{\XPoly}{i}\left(\ind{S}{i}\right)$ and $\ind{\XPoly}{j}\left(\ind{S}{j}\right)$ could interact sterically (see Fig.\ \ref{fig:SegSch}(c)), and we need to compute the double sum\ \eqref{eq:sumForceUp} using Gauss-Legendre nodes over $\ind{S}{i} \times \ind{S}{j}$ (see \hyperref[step7]{Step 7 in Algorithm 1}).

\subsection{Combining integration domains \label{sec:Unions}}
\begin{figure}
\centering
\includegraphics[width=0.5\textwidth]{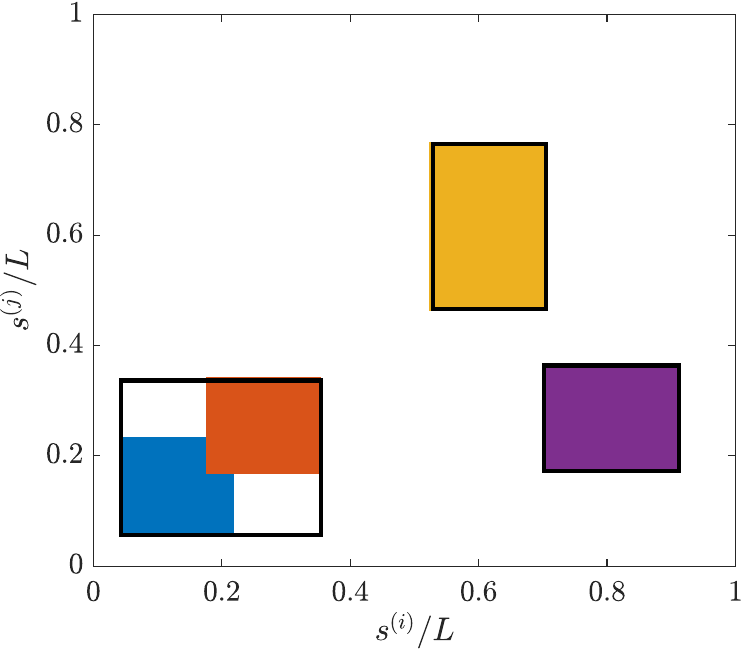}
\caption[Combining integration domains in two dimensions]{\label{fig:DomComb}Combining integration domains in two dimensions for the steric energy integral\ \eqref{eq:EDCInt}. The colored rectangles represent the integration domains returned by the quadratic approximation algorithm (Section\ \ref{sec:quadApprox}). We combine these four intervals into three (the black rectangles) for Gauss-Legendre integration.}
\end{figure}

This section describes how we combine the integration domains generated in the previous section to efficiently integrate\ \eqref{eq:EDCInt} over disjoint two-dimensional intervals. Consider the picture in Fig.\ \ref{fig:DomComb}. We identify the colored rectangles as intervals that were chosen via \eqref{eq:s1roots} on which the fibers have a nonzero steric interaction energy. Our procedure is to first sort the intervals in order of increasing starting $\ind{s}{i}$ location. Then, for each following interval, we check if there is overlap with the previous interval (the red and blue rectangles) in the $\left(\ind{s}{i},\ind{s}{j}\right)$ plane. If there is, then we extend the integration domain to encompass both domains by taking a sort of union, as shown in the black rectangle that outlines both the red and blue regions. If there is no overlap (as for the yellow and purple rectangles), then we add the interval to the integration list as a separate interval. This algorithm prevents having to integrate over the whole fiber if only the two endpoints are interacting (for example), while at the same time being simple to implement.

\subsection{Comparing segment-based algorithm to uniform points \label{sec:StCompare}}
Because the uniform point implementation and segment-based implementation are supposed to evaluate the same integrals, they ought to identify the same intervals on pairs of fibers, and consequently evaluate the same integrals. Thus, as the number of uniform points gets larger and $N_\delta$ in\ \eqref{eq:NGL} gets larger, we should converge to the same forces on the Chebyshev points. While we cannot guarantee a certain accuracy in our segment-based algorithm, we can set a goal for \emph{the segment-based algorithm to give the same accuracy as uniform points spaced $1/\epsc$ apart} for ``most'' configurations.

\begin{figure}
\centering
\includegraphics[width=\textwidth]{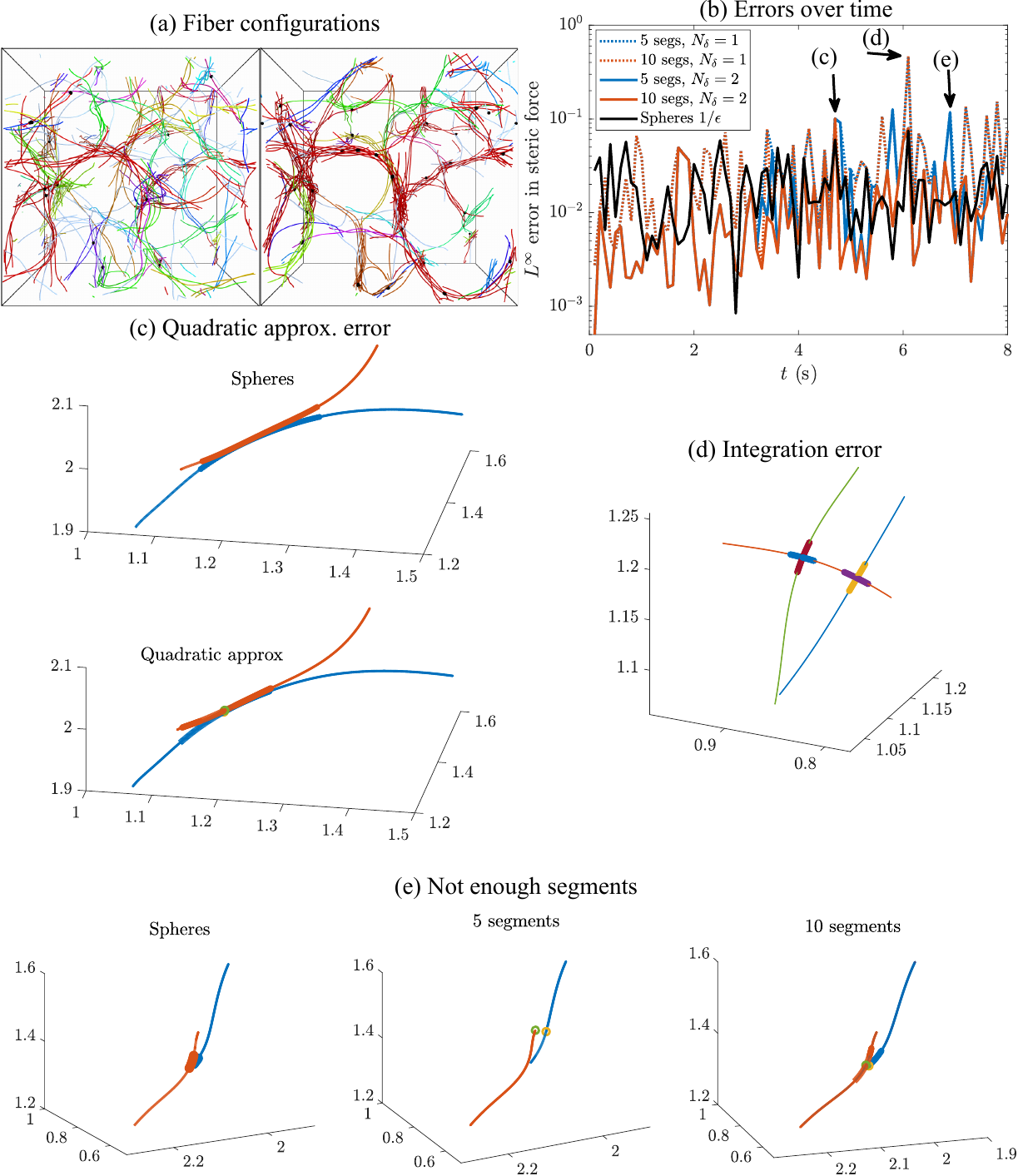}
\caption{\label{fig:StericEr}Using bundled configurations to test the segment-based algorithm and identify problem cases. (a) Samples of the fiber configurations we use, which come from bundling (Section \ref{sec:Bundling}) with $\epsc=4 \times 10^{-3}$. These snapshots are at $t=4$, and 8 s, with black circles showing points of overlap. (b) $L^\infty$ errors in the steric forces, computed using the segment-based algorithm (colored lines) and the uniform point algorithm with $1/\epsc$ points (black lines). In all cases, the reference solution uses $8/\epsc$ uniform points. For the segments, we consider 5 (blue) and 10 segments (red), and the dotted lines show $N_\delta=1$ in\ \eqref{eq:NGL}, while the solid lines show $N_\delta=2$. We isolate three cases with high errors (indicated by the arrows) for further study in (c--e).}
\end{figure}

To test this assumption, we consider configurations generated from simulations of fiber bundling with steric interactions (Section\ \ref{sec:Bundling}), in particular those shown in Fig.\ \ref{fig:StericEr}(a), which were simulated with $\epsc=4 \times 10^{-3}$. For each configuration, we compute steric forces using $8/\epsc$ uniform points (to obtain a reference force), then look at the errors with various $N_\delta$ and numbers of segments. The $L^\infty$ errors in the steric forces are shown in Fig.\ \ref{fig:StericEr}(b). We observe that using $1/\epsc$ uniform points gives 1--2 digits of accuracy in the forces. Typically, using $N_\delta=1$ for the Gauss-Legendre intervals gives the same accuracy as $1/\epsc$ uniform points, independent of how many segments we use. This is expected because the average point spacing ($\epsc L = \rc =\delta$) is the same in both cases. Increasing to $N_\delta=2$ gives an extra digit of accuracy, but at the cost of more expensive quadratures. 

There are a few cases in Fig.\ \ref{fig:StericEr} where we obtain larger than expected errors. The sources of error are as follows:
\begin{enumerate}
\item Errors in using the quadratic approximation of Section\ \ref{sec:quadApprox} to approximate the filament near the local minimum. For large $\Delta s$, this approximation can underestimate the size of the region where the filaments are closer than $\rmax$, thus giving an error in the steric force (see Fig.\ \ref{fig:StericEr}(c), where the region of interaction is shown using thicker lines). Unlike cases (d) and (e), in this case there is no way for us to reduce the error; as shown in arrow (c) in Fig.\ \ref{fig:StericEr}, the error does not change when we change the number of segments or $N_\delta$.
\item Inaccuracies in the quadrature scheme, which can happen when $N_\delta$ is too small. The arrow (d) in Fig.\ \ref{fig:StericEr} is an example of this, as the error decreases when we increase $N_\delta$. The specific fibers that are causing the problem are isolated in Fig.\ \ref{fig:StericEr}(d), where the regions of overlap ($r \leq \rmax$) are shown using thick lines. We see that the fibers intersect each other at a 90 degree angle (distance between the lines $< 10^{-3}$ for both intersections), which means that the integrand is sharply peaked and decays rapidly. This is what causes the inaccuracies in the quadrature, and why we need larger $N_\delta$ to get a more accurate result.
\item Not finding all the local minima, which can happen when there are too few segments. The arrow (e) in Fig.\ \ref{fig:StericEr} is an example of this, as the error decreases when we add more segments. Figure\ \ref{fig:StericEr}(e) shows the pair of fibers that causes the problem here. The closest point of contact is near the fiber endpoint, and when there are five segments the initial guess for Newton's method \emph{is} the fiber endpoint. Because of the way the fiber curves at the endpoint, the gradient there is positive (descent direction is negative, off the fiber), and so the projected Newton's method forces us to stay at the endpoint. This is fixed in the 10 segment case, where we have a better initial guess and we find the true local minimum. For this reason, we use 10 segments in all simulations in this work.
\end{enumerate}

\end{appendices}

\end{document}